\documentclass[a4paper]{amsart}
\usepackage{amssymb}
\usepackage{mathtools}
\usepackage{mathrsfs}
\usepackage[hidelinks]{hyperref}
\usepackage{tikz}
\usepackage{graphicx}
\usepackage{todonotes}
\usepackage[font=small]{caption}
\usetikzlibrary{decorations.markings}
\usepackage{enumitem}

\tikzset{
	wing/.style={
		decoration={
			markings, mark=at position 0.8 with {\arrow{latex}}
		},
		postaction={decorate},
	},
}

\newcommand{\Z}{\mathbb{Z}}
\newcommand{\R}{\mathbb{R}}
\newcommand{\C}{\mathbb{C}}
\newcommand{\T}{\mathbb{T}}
\newcommand{\ii}{{\rm i}}
\newcommand{\cent}[1]{#1^\circ}

\newcommand{\graph}{\mathsf}
\newcommand{\gV}{\graph{V}}
\newcommand{\gE}{\graph{E}}

\newcommand{\gF}{\graph{F}}
\newcommand{\gH}{\graph{H}}
\newcommand{\gG}{\graph{G}}

\newcommand{\gB}{\graph{B}}
\newcommand{\gC}{\graph{C}}

\newcommand{\surface}{\mathscr}
\newcommand{\sM}{\surface{M}}
\newcommand{\sS}{\surface{S}}

\newcommand{\Fhor}{F^\text{hor}}
\newcommand{\Phor}{P^\text{hor}}
\newcommand{\Fver}{F^\text{ver}}
\newcommand{\Pver}{P^\text{ver}}

\newcommand{\cH}{\mathcal{H}}
\newcommand{\cS}{\mathcal{S}}
\newcommand{\cA}{\mathcal{A}}
\newcommand{\cB}{\mathcal{B}}
\newcommand{\cC}{\mathcal{C}}
\newcommand{\cV}{\mathcal{V}}

\newcommand{\cPhor}{\mathcal{P}^\text{hor}}
\newcommand{\cFhor}{\mathcal{F}^\text{hor}}
\newcommand{\cPver}{\mathcal{P}^\text{ver}}
\newcommand{\cFver}{\mathcal{F}^\text{ver}}

\newcommand{\rotate}{\varsigma}
\newcommand{\involute}{\iota}

\newcommand{\grad}{\operatorname{grad}}
\renewcommand{\div}{\operatorname{div}}
\newcommand{\curl}{\operatorname{curl}}
\newcommand{\mdiv}{\operatorname{mdiv}}

\newcommand{\re}{\operatorname{Re}}
\newcommand{\im}{\operatorname{Im}}
\newcommand{\Res}[2]{\operatorname{Res}\left(#1,#2\right)}

\newcommand{\rank}{\operatorname{rank}}

\theoremstyle{plain} 
\newtheorem{theorem}{Theorem}[section]
\newtheorem{lemma}[theorem]{Lemma}
\newtheorem{corollary}[theorem]{Corollary}
\newtheorem{proposition}[theorem]{Proposition}

\theoremstyle{definition} 
\newtheorem{definition}[theorem]{Definition}
\newtheorem{assumption}[theorem]{Assumption}
\newtheorem{example}[theorem]{Example}

\theoremstyle{remark}
\newtheorem{remark}[theorem]{Remark}

\begin{document}

\title[Gluing saddle towers I: TPMS]{Gluing Karcher--Scherk saddle towers I:\\Triply periodic minimal surfaces}

\author{Hao Chen}
\address[Chen]{Georg-August-Universit\"at G\"ottingen, Institut f\"ur Numerische und Angewandte Mathematik}
\email{h.chen@math.uni-goettingen.de}
\thanks{H.\ Chen is supported by Individual Research Grant from Deutsche Forschungsgemeinschaft within the project ``Defects in Triply Periodic Minimal Surfaces'', Projektnummer 398759432.
M.\ Traizet is supported by the ANR project Min-Max (ANR-19-CE40-0014)}

\author{Martin Traizet}
\address[Traizet]{Institut Denis Poisson, CNRS UMR 7350, Faculté des Sciences et Techniques, Université de Tours }
\email{martin.traizet@univ-tours.fr}


\date{\today}

\begin{abstract}
	We construct minimal surfaces by gluing simply periodic Karcher--Scherk
	saddle towers along their wings.  Such constructions were previously
	implemented assuming a horizontal reflection plane.  We break this symmetry
	by prescribing phase differences between the saddle towers.  It turns out
	that, in addition to the previously known horizontal balancing condition, the
	saddle towers must also be balanced under a subtle vertical interaction.
	This interaction vanishes in the presence of a horizontal reflection plane,
	hence was not perceived in previous works.

	Our construction will be presented in a series of papers.  In this first
	paper of the series, we will explain the background of the project and
	establish the graph theoretical setup that will be useful for all papers in
	the series.  The main task of the current paper is to glue saddle towers into
	triply periodic minimal surfaces (TPMSs).  Our construction expands many
	previously known TPMSs into new 5-parameter families, therefore significantly
	advances our knowledge on the space of TPMSs.
\end{abstract}

\maketitle

\section{Background} \label{sec:intro}

In the last decades, the node-opening technique has been very successful in
gluing catenoids into minimal surfaces of finite or infinite topology in
Euclidean space forms~\cite{traizet2002, traizet2002b, traizet2008,
morabito2012, chen2021b}.

In fact, the technique was first developed to glue Karcher--Scherk saddle
towers. More specifically, the second named author desingularized arrangements
of vertical planes into minimal surfaces by replacing the intersection lines
with Scherk surfaces (saddle towers with four wings).  This was first done by
solving non-linear PDEs~\cite{traizet1996} and later using the node-opening
technique~\cite{traizet2001}.  In his thesis~\cite{younes2009}, Rami Younes
desingularized $\gG \times \R$, where $\gG$ is a graph embedded in a flat
2-torus, into triply periodic minimal surfaces by placing saddle towers over
the vertices of $\gG$.  These constructions work only when all saddle towers
are balanced under a horizontal interaction or, equivalently, when the graph
embedding is a non-degenerate critical point of the length functional.

\begin{figure}[hbt]
\includegraphics[height=0.4\textwidth]{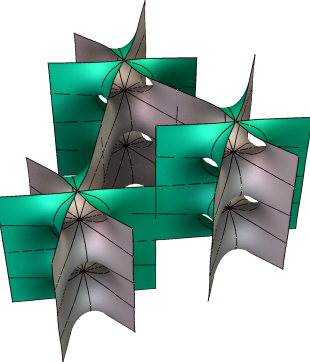}
\hspace{1cm}
\includegraphics[height=0.4\textwidth]{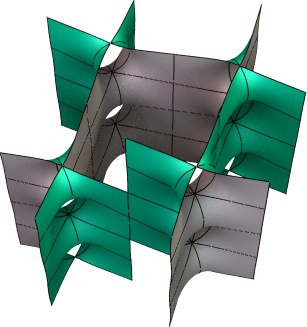}
\caption{
	A rhombohedral (left) and a tetragonal (right) deformation of the Gyroid near
	the saddle tower limits.  They do not have any horizontal symmetry plane.
}
\label{fig:tGrGL}
\end{figure}

All these previous works assumed, however, that the surfaces have a horizontal
reflection plane.  It was wondered in~\cite{traizet1996} whether this symmetry
is necessary for simply periodic minimal surfaces (SPMSs) with ends of Scherk
type (i.e.\ asymptotic to vertical planes).  A positive answer was given for
SPMSs of genus 0~\cite{perez2007}, but Mart\'in and Ramos Batista constructed
the Karcher--Costa towers~\cite{martin2006} that provide examples of genus 1
without horizontal symmetry.  The same question can be asked for triply
periodic minimal surfaces (TPMSs). But the deformation families of the Gyroid,
which were recently confirmed by the first named author~\cite{chen2019} admit
saddle tower limits but do not have any horizontal symmetry plane (see Figure
\ref{fig:tGrGL}).

The goal of this series of papers is to glue Karcher--Scherk saddle towers into
minimal surfaces, without the assumption of a horizontal reflection plane.

Let $\gG$ be a graph in the complex plane.  Formal definitions will be given in
Section~\ref{sec:graphs}.  At the moment, it suffices to understand that $\gG$
is a set of straight segments (edges) and half-lines (rays) that intersect only
at their endpoints (vertices).  There are different possible setups, depending
on weather the graph $\gG$ is finite, periodic with finite quotient, or
infinite aperiodic.

For sufficiently small $\varepsilon>0$, we want to construct a 1-parameter
family $(\sM_\varepsilon)_{\varepsilon>0}$ of embedded minimal surfaces of
vertical period $2\pi\varepsilon^2$ that tends to $\gG \times \R$ as
$\varepsilon \to 0$.  For this purpose, we place suitably scaled saddle towers
over the vertices of $\gG$ and glue their wings along the edges of the graph.

Since the fluxes along the wings of a saddle tower sum up to $0$, the graph
$\gG$ must be balanced in the sense that for each vertex $v$, the sum of the
unit vectors in the directions of the outgoing edges adjacent to $v$ is zero.
We call the sum the horizontal force.

Orientability of $\sM_{\varepsilon}$ requires the graph $\gG$ to be orientable,
in the sense that its faces can be labelled with $+$ or $-$ signs so that
adjacent faces have opposite signs. The prototype result in the symmetric case
is the following:

\begin{theorem}[Informal statement] \label{metatheorem1}
	Let $\gG$ be a balanced, rigid and orientable graph. Then $\gG\times\R$ can
	be desingularized into a family of minimal surfaces $\sM_\varepsilon$ that
	have vertical period $2\pi\varepsilon^2$ and are symmetric with respect to a
	horizontal plane.
\end{theorem}
Rigidity concerns the invertibility of the Jacobian of horizontal forces. Its
precise formulation depends on the setup under consideration.

Theorem \ref{metatheorem1} was first proved in \cite{traizet2001} when $\gG$
appears as a line arrangement, yielding a family of SPMSs with Scherk-type
ends.  It was then proved in \cite{younes2009} when $\gG$ is a doubly-periodic
graph with finite quotient, yielding a family of TPMSs.  See Theorem
\ref{thm:rami} for a precise statement in this case.

\medskip

To break the horizontal symmetry, we must move the saddle towers vertically so
that they have different reflection planes.  The ``height differences'' between
the reflection planes, known as the phase differences, are prescribed on the
edges of the graph.

Technically, the phase differences correspond to the complex arguments of the
node-opening parameters; see Section~\ref{ssec:opennode} for details.  Complex
node-opening parameters were considered in previous works~\cite{traizet2008},
but never did the complex argument play such an important role.

This paper reveals that the phase differences must satisfy a balancing
condition, known as vertical balancing, which is more subtle and delicate than
horizontal balancing of the graph $\gG$.  The formulation of vertical balancing
requires quite a lot of preliminaries, but we can give an idea of what it looks
like.

Assume that two saddle towers $\sS$ and $\sS'$ are glued along their wings.
Let $\ell$ be the length of the corresponding edge in the graph, and $\phi$ be
the prescribed phase difference of $\sS'$ over $\sS$.  Then, as $\varepsilon
\to 0$, the vertical force exerted by $\sS'$ on $\sS$ is asymptotically of the
form
\[
	K \sin(\phi) \exp(-\ell/\varepsilon^2).
\]
where the coefficient $K$ is a positive real number that depends on the
undulation of the saddle towers and a first-order deformation of the graph.
Because of the factor $\sin(\phi)$, the force vanishes whenever $\phi =0$ or
$\pi$.  This explains why this interaction was not perceived in previous works
that assumed a reflection symmetry.

The vertical balancing condition requires, for every set of saddle towers, that
the forces sum up to $0$ over the edges that separate the set from other saddle
towers.  Because of the exponential factor, the interaction dominates on
shorter edges of the graph for small $\varepsilon$.  So in the limit
$\varepsilon \to 0$, it suffices to take the sum of $K\sin(\phi)$ over the
shortest separating edges.

The prototype result in the non-symmetric case is the following:

\begin{theorem}[Informal statement] \label{metatheorem2}
	Let $\gG$ be a balanced, rigid and orientable graph with phase differences
	$\phi$ prescribed on the edges. Assume that $\phi$ is balanced and rigid.
	Then $\gG\times\R$ can be desingularized into a family of minimal surfaces
	$\sM_\varepsilon$ which have vertical period $2\pi\varepsilon^2$.  Moreover,
	for each vertex $v$, $\varepsilon^{-2}\sM_{\varepsilon}$ converges after
	suitable horizontal translation to a saddle tower $\sS_v$, and the phase
	differences between adjacent saddle towers are given by $\phi$.
\end{theorem}

Rigidity of $\phi$ concerns the invertibility of the Jacobian of vertical
forces.  Again, its precise formulation depends on the given setup.

\medskip

\begin{figure}[bht]
	\small
	\begin{tabular}[c]{ccc}
		\includegraphics[height=0.3\textwidth]{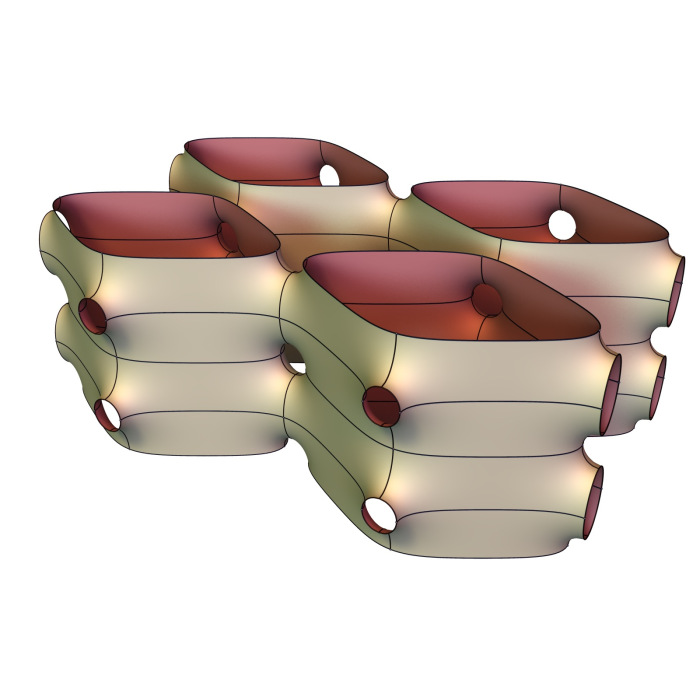}&
		\includegraphics[height=0.3\textwidth]{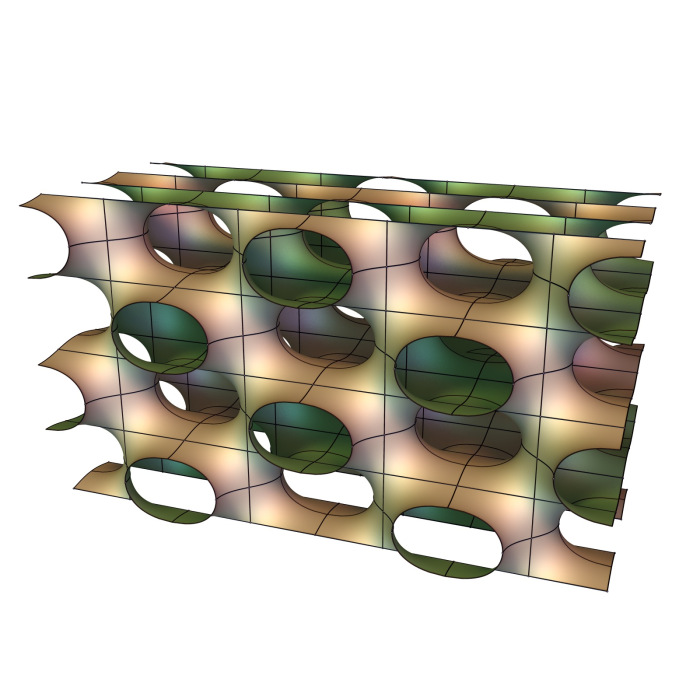}&
		\includegraphics[height=0.3\textwidth]{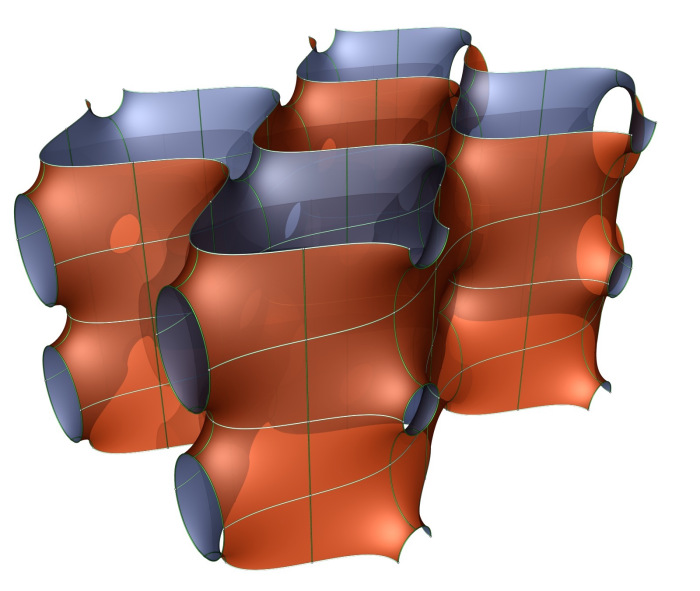}\\
		(a) oPa & (b) oPb & (c) oCLP\\
		\includegraphics[height=0.3\textwidth]{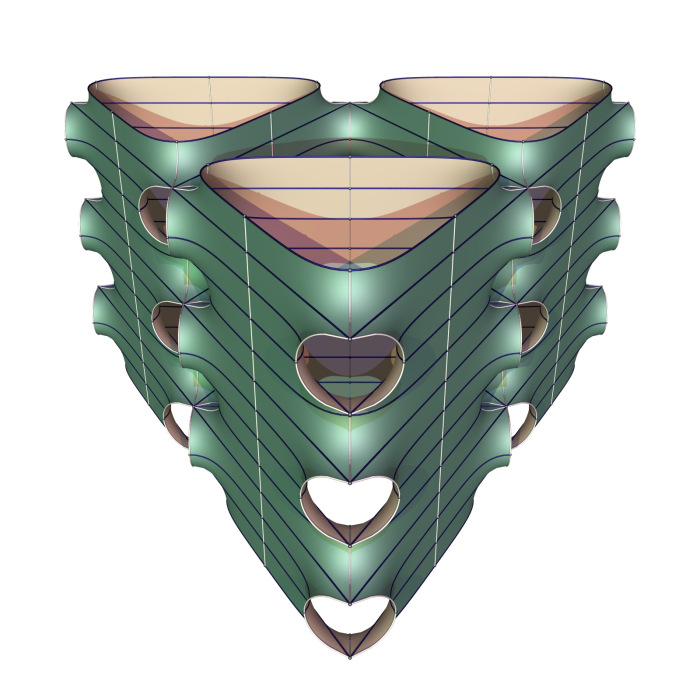}&
		\includegraphics[height=0.3\textwidth]{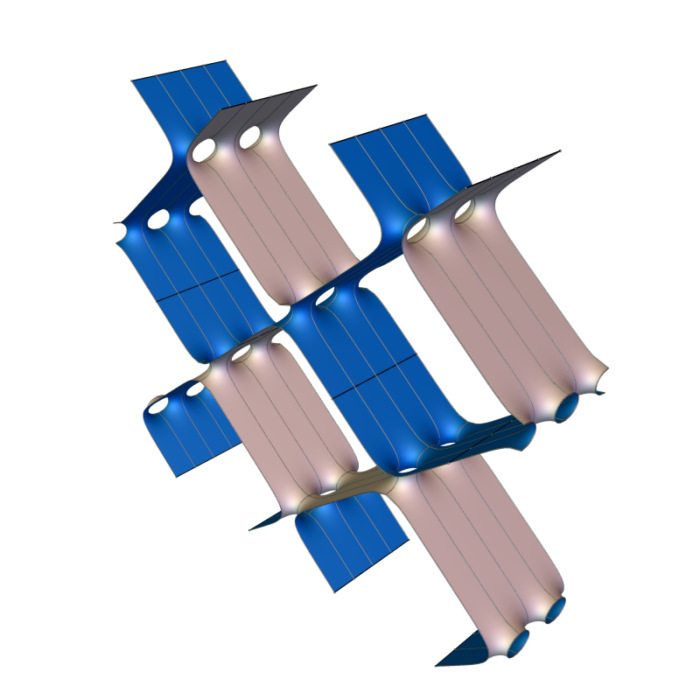}&
		\includegraphics[height=0.3\textwidth]{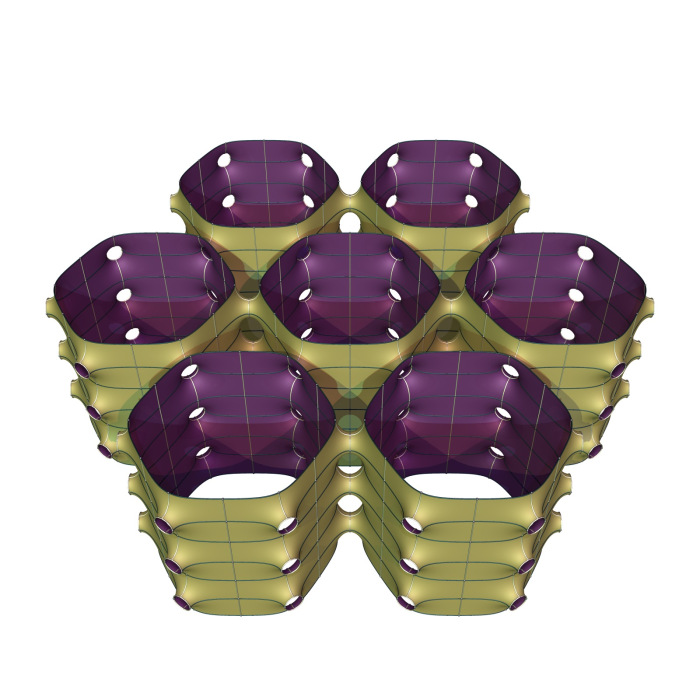}\\
		(d) H & (e) oH & (f) H'--T\\
		\includegraphics[height=0.3\textwidth]{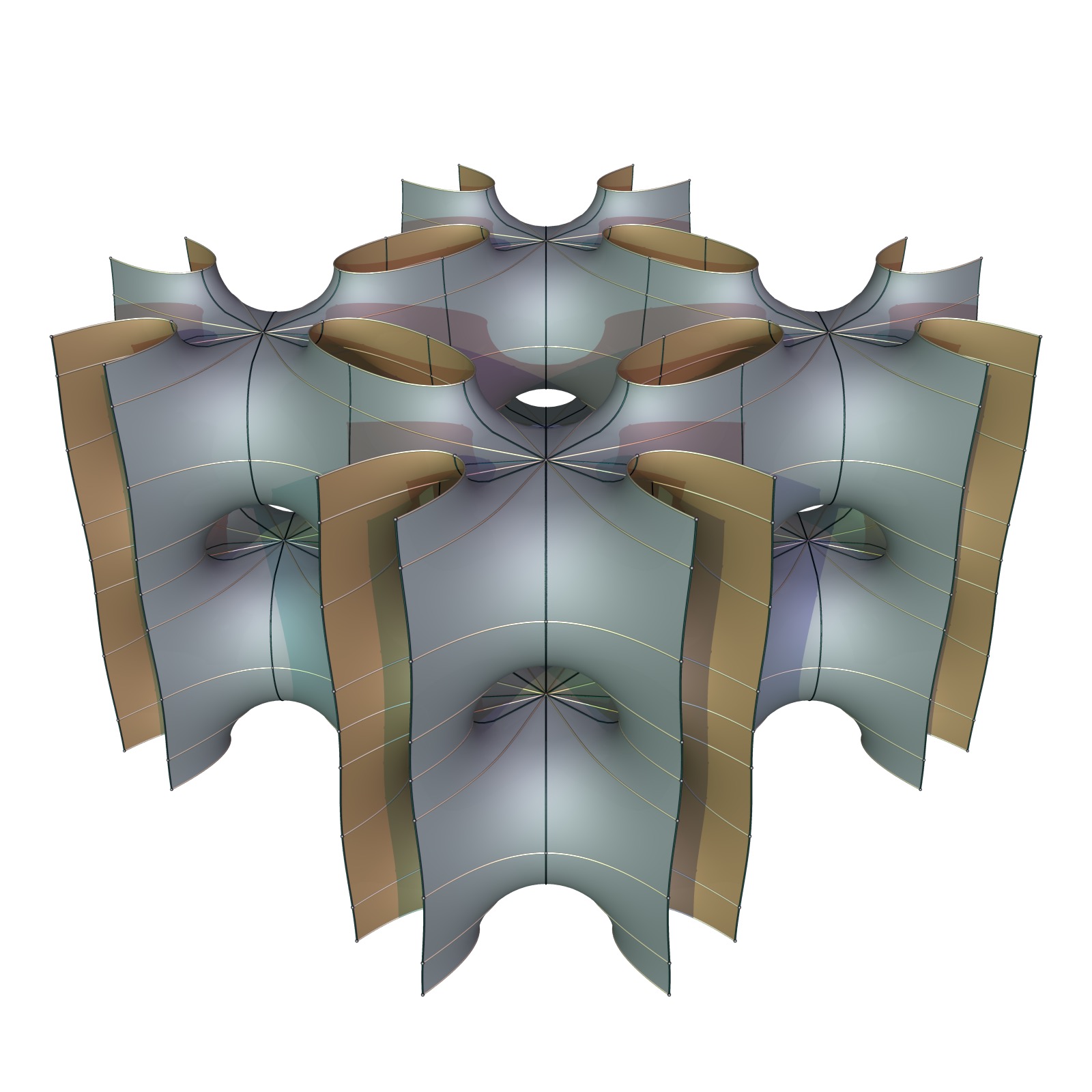}&
		\includegraphics[height=0.3\textwidth]{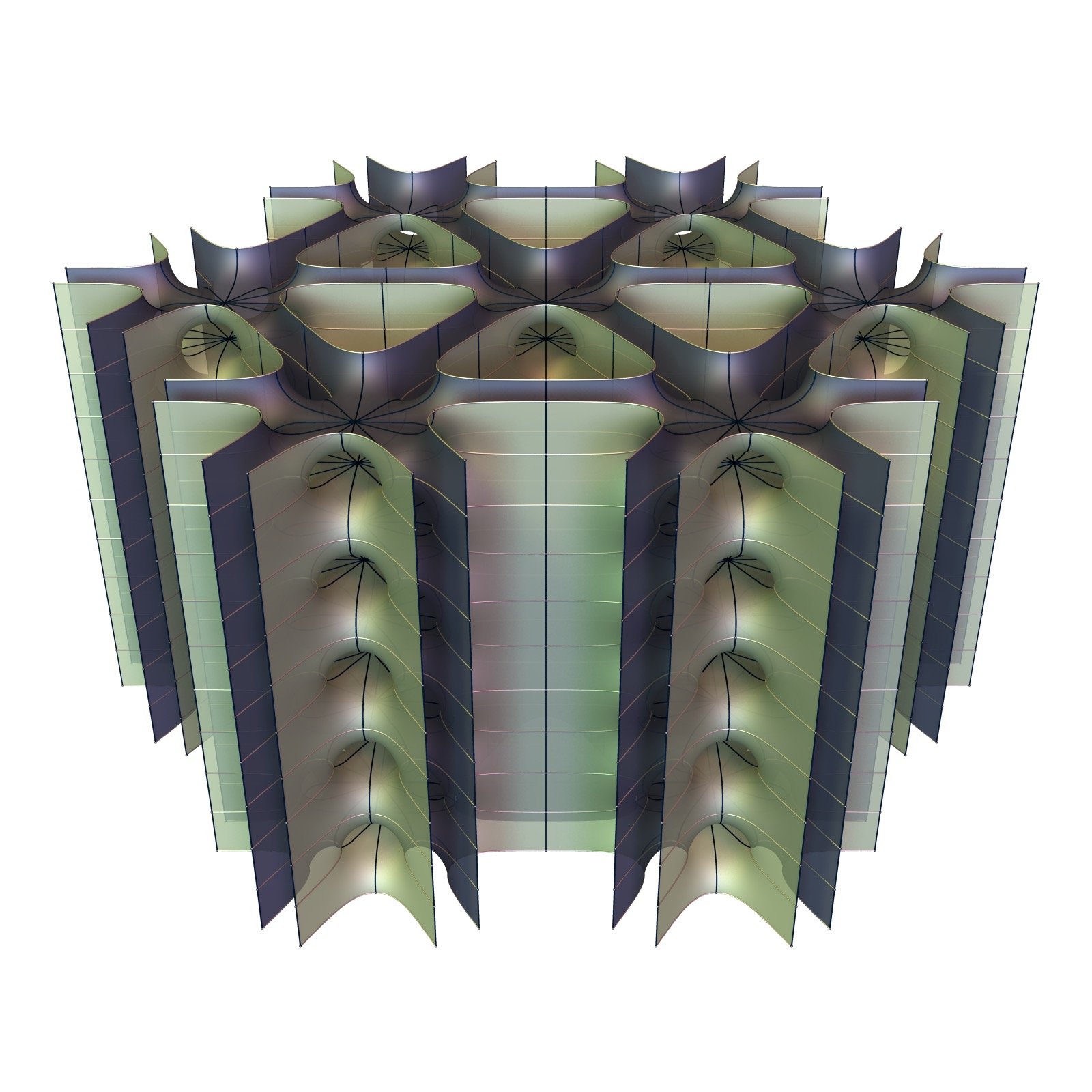}&
		\includegraphics[height=0.3\textwidth]{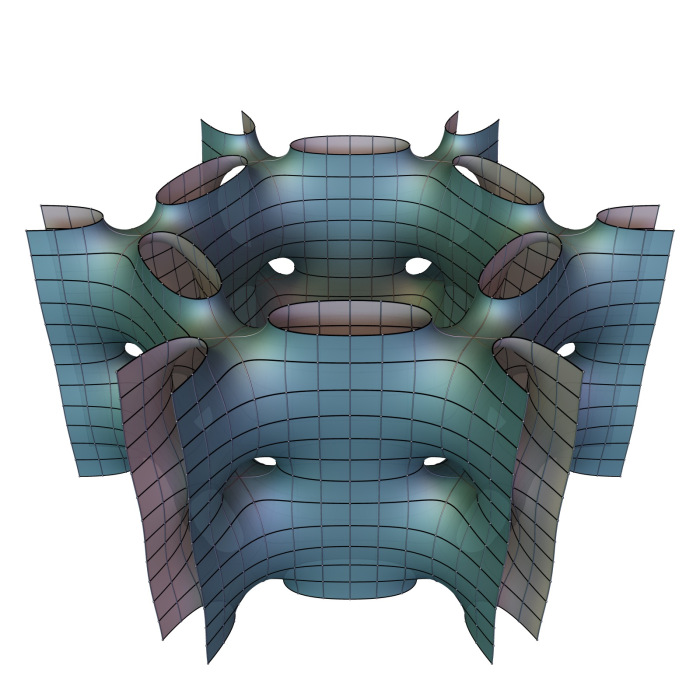}\\
		(g) S'--S'' & (h) T'--R' & (i) H''-R
	\end{tabular}
	\caption{
		A gallery of known examples of saddle tower limits of triply periodic
		minimal surfaces (Source: Matthias Weber).  Generalizations of them will be
		constructed in this paper.  \label{fig:gallery}
	}
\end{figure}

The goal of this first paper in the series is to prove Theorem
\ref{metatheorem2} when $\gG$ is a doubly-periodic graph with finite quotient,
resulting in TPMSs.  See Theorem~\ref{thm:main} for the precise formulation of
our main result.  TPMSs are given the privilege because many known families of
TPMSs admit saddle tower limits, but most of them are symmetric in a horizontal
plane; see Figure~\ref{fig:gallery}.  So we expect to find many interesting new
examples after breaking this symmetry.

\medskip

The paper is organized as follows.

In Section~\ref{sec:jenkins-serrin}, we recall some basic facts about
Karcher--Scherk saddle towers and give a geometric definition of phases.
Section~\ref{sec:graphs} sets up the graph theoretical foundation, which is
crucial for formal statements of the main results.  Then in
Section~\ref{sec:horizontal}, we define horizontal balance and rigidity, and
recall Younes' result about symmetric TPMSs.

Section~\ref{sec:towers} is dedicated to an in-depth investigation on the
shapes of the wings.  We will define new quantities in terms of Weierstrass
data.  They will be useful in Section~\ref{sec:main} for defining vertical
balance and rigidity, and for announcing our main Theorem~\ref{thm:main} about
non-symmetric TPMSs.  Then we present in Section~\ref{sec:example} examples of
TPMSs that can be produced from our construction, including many new TPMSs of
genus 3.  The construction is finally proved in Section~\ref{sec:construction}.

Note that the sections~\ref{sec:jenkins-serrin}, \ref{sec:graphs},
\ref{sec:towers} are preparatory, meant to be relevant not only for the current
manuscript but also for future papers.  

\subsection*{Acknowledgement}

The authors thank Matthias Weber for an unbelievable video of triply periodic
minimal surfaces that motivated this project.

\section{First look on Karcher--Scherk saddle towers}
\label{sec:jenkins-serrin}

Recall that a \emph{saddle tower} is an embedded SPMS with $n$ Scherk-type
ends, where $n\ge 4$ is an even integer.  Saddle towers with $n=4$ were
discovered by Scherk~\cite{scherk1835} in the 19th century.  All other saddle
towers were constructed by Karcher~\cite{karcher1988}.  Saddle towers are
classified in \cite{perez2007} as the only embedded SPMS with $n$ Scherk-type
ends and genus zero in the quotient by the period.

We sketch Karcher's construction as follows: Let $P$ be a convex polygon with
$n$ sides of length one, where $n\geq 4$ is even.  We allow $P$ to be
non-strictly convex but exclude the \emph{degenerate} case where $P$ appears as
a line segment of length $n/2$, and the \emph{special} case where $n\geq 6$ and
$P$ appears as a parallelogram with 2 sides of length one~\cite{perez2007}.
Then the Jenkins--Serrin theorem~\cite{jenkins1966} guarantees the existence of
a minimal graph over the interior of $P$ that takes the values $\pm\infty$,
alternately, on the edges of $P$, unique up to vertical translation.  This
graph is bounded by vertical lines over the vertices of $P$.  Its conjugate
minimal surface is bounded by $n$ curvature lines that lie, alternately, in two
horizontal planes at distance 1 from each other.  Reflections in these planes
extend the conjugate surface into an SPMS called a saddle tower, which we
denote by $\sS$.

The edges of $P$ are $n$ unit vectors that sum up to $0$.  To be consistent
with later notations, we label the edges by $\gH = \{0,\cdots,n-1\}$ in the
counterclockwise order, and define $\rotate: \gH\to\gH$ by $\rotate(h) = h+1
\pmod{n}$.  Let the corresponding unit vectors be $u_h = \exp(\ii \theta_h) \in
\C \simeq \R^2$, $h \in \gH$.  The corresponding ends of $\sS$ are called
\emph{wings}.  Each wing is asymptotic to a vertical half-plane that is
parallel to the corresponding edge.  Wings extending in the same direction
(same $\theta_h$) are called \emph{parallel}.

We now make a distinction between the two horizontal reflection planes.  Fix an
orientation of $\sS$.  It is known that the conjugate minimal surface of the
Jenkins--Serrin graph is a graph over an unbounded concave domain $\Omega$.
The complement of $\Omega$ has $n$ unbounded convex components, each bounded by
the projection of a curvature line (see \cite{karcher1988}).  We say that a
curvature line $\gamma$ is a $0$-arc if the Gauss map on $\gamma$ points into
the corresponding convex component.  All $0$-arcs lie on the same horizontal
symmetry plane, which we call the $0$-plane.  We scale all our saddle towers so
that their vertical period is $2\pi$.

\begin{definition}[Phase]\label{def:phase}
	Let $\sS$ be an oriented saddle tower of vertical period $2\pi$.  We say that
	the \emph{phase} of $\sS$ is $\phi \in \R/2\pi\Z$ if the horizontal plane
	$x_3 = \phi$ is a $0$-plane of $\sS$.
\end{definition}

\begin{figure}
	\includegraphics[width=.8\textwidth]{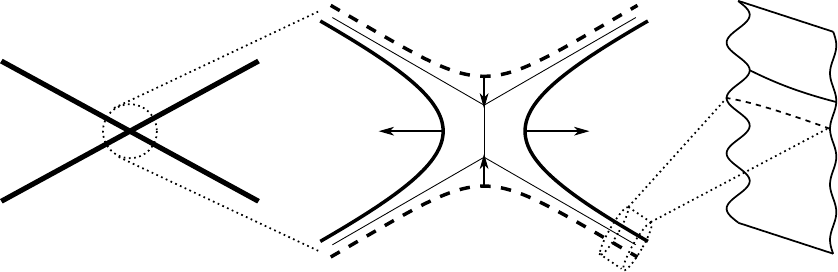}
	\caption{
		A four-wing saddle tower seen from a distance (left), then scaled at the
		``axis'' (middle).  The shape of a wing is illustrated on the
		right, featuring the undulation.\label{fig:tower}
	}
\end{figure}

\begin{remark} \label{remark:tower}
	Seen from a distance, $\sS$ looks like $n$ vertical half-planes sharing a
	vertical boundary; see Figure~\ref{fig:tower} (left).  But a closer look
	reveals that the asymptotic planes of the wings do not intersect at a single
	vertical line; see Figure~\ref{fig:tower} (middle) for an illustration and
	Section~\ref{sec:symtower} for a detailed analysis.
\end{remark}

\begin{figure}
	\includegraphics[height=0.5\textwidth]{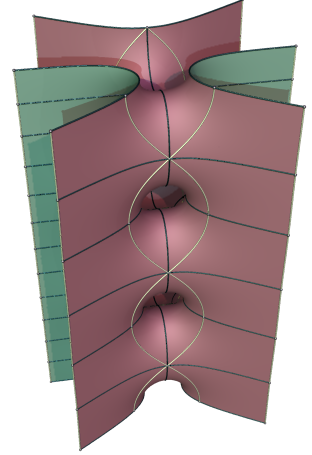}
	\includegraphics[height=0.5\textwidth]{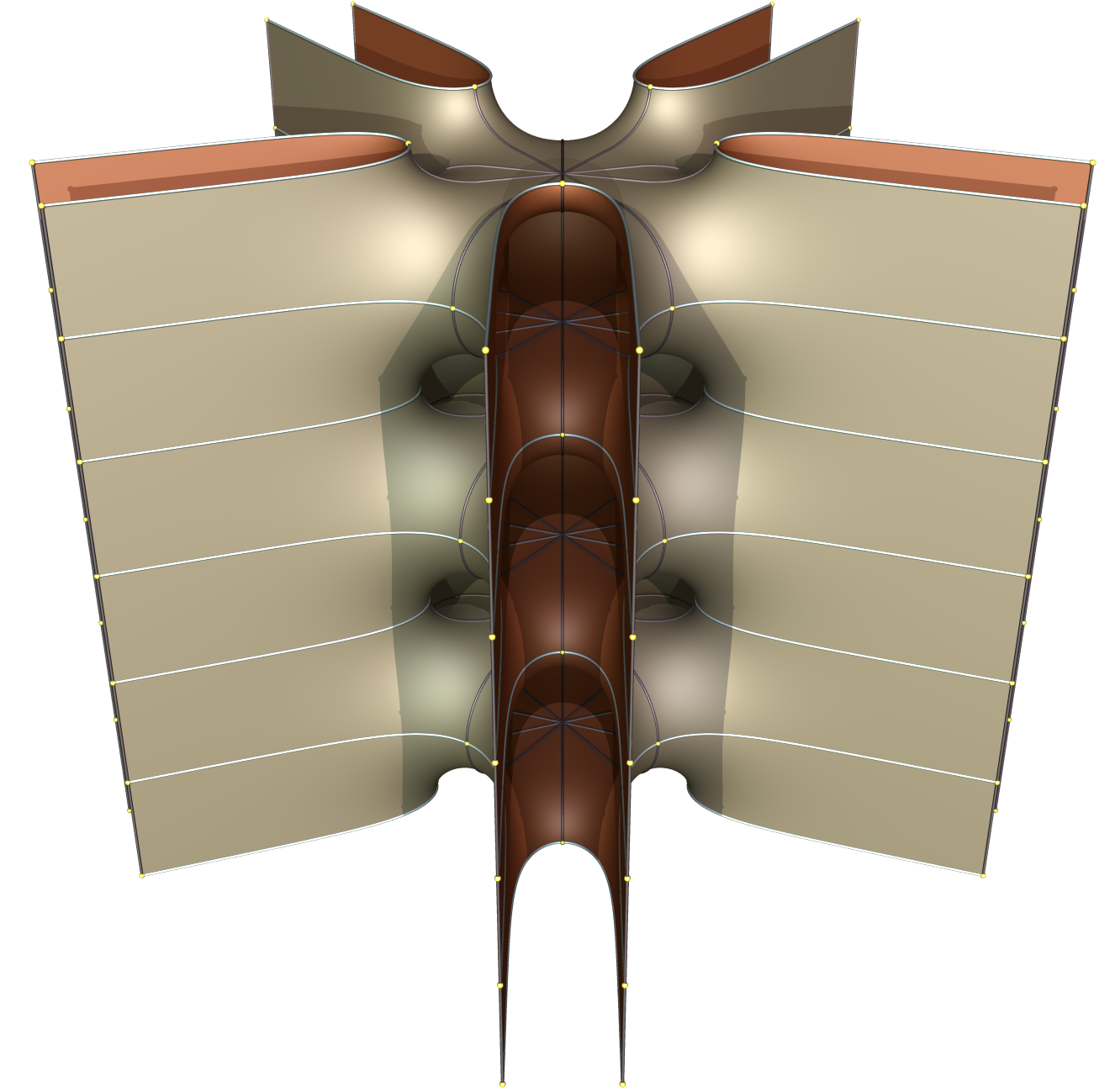}
	\caption{
		Two saddle towers with six and ten wings; see
		Section~\ref{sec:exampletower}.  The one on the right has five pairs of
		parallel wings. (Source: Matthias Weber)\label{fig:Matthias}
	}
\end{figure}

\section{Graphs}\label{sec:graphs}

Graph theory plays a central role in our construction.  Combinatorially, we use
a graph to describe how saddle towers (vertices) are glued along wings (edges)
into a minimal surface.  Geometrically, the degenerate limit of our surface
projects to a horizontal plane as a geometric ``representation'' of the graph.

In this paper, a ``graph'' actually refers to a multigraph, possibly with
multiple edges and loops.  This maximizes the generality of our construction.
Such a multigraph is best described in terms of rotation systems.  The standard
references for this section are~\cite{godsil2001, mohar2001, diestel2018}.

\subsection{Intuitive picture}
\label{sec:graph-intuitive}

For readers not familiar with rotation systems, we give here an intuitive
picture which, although informal, should suffice for understanding our
construction.

A graph $\gG$ in a flat 2-torus $\T^2$ can be understood as a set $\gV$ of
points (vertices) and a set $\gE$ of straight segments (edges) whose endpoints
are (not necessarily distinct) vertices and whose interiors are disjoint.  The
connected components of $\T^2 \setminus \gG$ are called faces.  The set of
faces is denoted $\gF$.  By Euler formula,
\[
	|\gV|-|\gE|+|\gF|=0.
\]
Each edge has two possible orientations, i.e.\ a choice of ``initial vertex''.
The set of oriented edges of $\gG$ is denoted $\gH$, so $|\gH|=2|\gE|$.  For
$h\in \gH$, we denote $e(h)\in \gE$ the corresponding unoriented edge and $-h$
the same edge with opposite orientation.  The initial vertex of an oriented
edge is denoted $v(h)$, so its terminal vertex is $v(-h)$.

\begin{remark}
	To agree with the language of rotation systems, we call elements of $\gH$
	half-edges and the adjacency relation is denoted $h\in v$.  This terminology
	is preferred because vertices correspond to saddle towers and half-edges
	correspond to wings of saddle towers.
\end{remark}

The number of half-edges adjacent to a vertex $v$ is called the degree of $v$
and denoted $\deg(v)$.  The wings of a saddle tower are naturally ordered by a
cyclic permutation.  Accordingly, we define a permutation $\rotate$ on $\gH$ as
follows: if $h\in v$, $\rotate(h)\in v$ is the half-edge which comes after $h$
when traveling around $v$ in the counterclockwise direction.

\subsection{Formal definition}

In the case that a saddle tower has parallel wings, they project to the same
straight segment in the degenerate limit.  Then the intuitive picture above,
which defines edges geometrically, cannot distinguish them properly.  To
include this situation, it is necessary to use the language of rotation system.

A \emph{rotation system} consists of a set $\gH$ and two permutations
$\involute$ and $\rotate$ acting on $\gH$, such that $\involute$ is an
involution without fixed points and the group generated by $\involute$ and
$\rotate$ acts transitively on $\gH$.  The elements of $\gH$ are called
half-edges.  To ease notations and be consistent with the intuitive picture in
Section \ref{sec:graph-intuitive}, we write $-h$ for $\involute(h)$.

\subsubsection{The multigraph}

The rotation system $(\gH,\involute,\rotate)$ defines a connected multigraph
$(\gH, \gV, \gE)$, where the vertex set $\gV$ consists of the orbits of
$\rotate$ and the edge set $\gE$ consists of the orbits of $\involute$.  As we
have mentioned before, the vertices correspond to saddle towers, the half-edges
correspond to wings, and the edges correspond to glued wing pairs.  In this
sense, the graph describes the gluing pattern for our construction.

Note that vertices and edges are identified with subsets of $\gH$, so $h \in v$
means that the half-edge $h$ is adjacent to the vertex $v$.  For a half-edge $h
\in \gH$, we use $v(h)$ and $e(h)$ to denote the unique vertex and edge
associated to $h$.  We say that $e(h)$ is a \emph{loop} if $v(h) = v(-h)$.

The orbits of $\rotate\involute$ are called (combinatorial) \emph{faces} of the
graph, and the set of faces is denoted by $\gF$.

\begin{assumption} \label{assumption-loops}
	All faces have at least two elements.
\end{assumption}

Among the edges, we define an equivalence relation known as ``parallel'', such
that $e(h)$ is \emph{parallel} to $e(h')$ whenever $\{h, h'\}$ is a face.
Because an equivalence relation is transitive, we have $e$ is parallel to $e'$
whenever there is a sequence $h_0, \cdots, h_n$ such that $e = e(h_0)$, $e' =
e(h_n)$, and $\{h_i, h_{i+1}\}$ is a face for all $0 \le i < n$.

\subsubsection{Topological embedding}

The rotation system $(\gH, \involute,\rotate)$ also defines, up to
homeomorphism, a 2-cell embedding of the multigraph on a closed oriented
surface.

Recall that an embedding represents vertices by distinct points and edges by
curves that do not intersect in their interiors.  The half-edges then
correspond to the curves in small neighborhoods of the vertices.  The
permutation $\rotate$ sends a half-edge to the next half-edge around the same
vertex in the counterclockwise direction, and $\involute$ sends a half-edge to
the other half-edge of the same edge.

The connected components of the complement of a 2-cell embedding are all
homeomorphic to an open disk, and are called (topological) \emph{faces} of the
embedding.  They are in correspondence with the combinatorial faces: Half-edges
on the boundary of a topological face in the clockwise direction form a
combinatorial face of the graph.  Note that two edges are parallel if and only
if their representatives are homologous.

If the graph is finite, the genus of the oriented surface can be calculated as
\[
	g = 1-\frac{1}{2}(|\gV|-|\gE|+|\gF|)
\]
and we call this number the \emph{genus} of the graph.

\begin{assumption}
	Graphs in the current paper are of genus $1$, so they are embedded in a
	2-torus $\T^2$.
\end{assumption}

\begin{remark}
	Graphs of genus $0$ arise when constructing SPMS instead of TPMS.  Higher
	genus graphs are certainly interesting and could be used to glue saddle
	towers in $\mathbb{H}^2 \times \R$.
\end{remark}

\subsubsection{Geometric representation}

We have seen that the rotation system $(\gH,\involute,\rotate)$ determines a
homeomorphism class of 2-cell embeddings.  We choose a geometric representation
$\varrho$ from the closure of this homeomorphism class.  More specifically,
$\varrho$ maps vertices to distinct points in a flat torus $\T^2$ and maps each
edge to a segment between the images of its endpoints, so that parallel edges
are mapped to the same segment, and non-parallel edges are mapped to segments
with disjoint interiors.

Such a geometric representation always exists. Indeed, let $\gG'$ be the graph
obtained from $\gG$ by merging parallel edges. We may apply the genus-one
version of Tutte's embedding theorem (see~\cite{gortler2006} for instance)
which implies that $\gG'$ admits a straight-edge representation on a flat
torus.  We then obtain a straight-edge representation of $\gG$ by mapping
parallel edges to the same segment that represents the corresponding edge of
$\gG'$.

\subsubsection{Orientation}

An \emph{orientation} of the graph $\gG$ is a function $\sigma \colon \gH \to
\{\pm 1\}$ such that $\sigma(-h) = -\sigma(h)$ for all $h$.

\begin{definition}
	An orientation $\sigma$ on $\gG$ is \emph{consistent} if
	$\sigma\circ\rotate = -\sigma$.  A graph is \emph{orientable} if it admits a
	consistent orientation.  Once a consistent orientation is fixed, we say that
	the graph is \emph{oriented}.
\end{definition}

Note that an orientable graph only has vertices of even degrees.  Vertices of
degree $2$ are not relevant for us, since there is no saddle tower
with two wings. 

\begin{assumption}\label{ass:degree}
	Graphs in this paper have only vertices of degree at least $4$.
\end{assumption}

\subsection{Vector spaces on graphs} \label{sec:vector-spaces}

\subsubsection{Cycles and cuts}

A (simple) \emph{cycle} is a set of half-edges $c \subset \gH$ that can be
ordered into a sequence $(h_1, \cdots, h_n)$ such that $v(-h_i)=v(h_{i+1})$ for
$1\leq i<n$, $v(-h_n)=v(h_1)$, and $v(h_i) \ne v(h_j)$ whenever $i \ne j$.  We
use $-c$ to denote the reversed cycle $\{-h \colon h \in c\}$.  The set of
cycles is denoted by $\gC$.  In particular, combinatorial faces are all cycles. 

For some partition $\gV=\gV_1 \sqcup \gV_2$ of the vertices, the \emph{cut}
between $\gV_1$ and $\gV_2$ is the set of half-edges $b \subset \gH$ such that
$v(h) \in \gV_1$ and $v(-h) \in \gV_2$ for all $h \in b$.  We use $-b$ to
denote the reversed cut $\{-h \colon h \in b\}$.  The set of cuts is
denoted\footnotemark~by $\gB$. In particular, for any vertex $v$, the set
\[
	b(v) = \{ h \in v \colon v(-h) \ne v \}
\]
is a cut between $\{v\}$ and $\gV \setminus \{v\}$.  We call $b(v)$ the
\emph{vertex cut} at $v$.

\footnotetext{Inclusion-wise minimal cuts are called \emph{bonds}.  Hence the
cut space is sometimes referred to as the bond space, therefore our notation.}

\subsubsection{Functions on half-edges}

Let $\cH$ be the space of functions $f\colon \gH \to \R$.  We say that $f \in
\cH$ is \emph{symmetric} if $f_{-h} = f_h$, and \emph{antisymmetric} if $f_{-h}
= -f_h$.  They can be seen as edge labelings on, respectively, undirected and
directed graphs.  The orientation $\sigma$ is an example of antisymmetric
function.  We use $\cS$ and $\cA$ to denote, respectively, the space of
symmetric and antisymmetric functions.  

We denote $e_h$ the characteristic function of $\{h\}$, so $(e_h)_{h\in\gH}$ is
the canonical basis of $\cH$.  We equip $\cH$ with the inner product
$(\cdot,\cdot)$ defined by $(e_h, e_{h'}) = \delta_{h,h'}$.  Then $\cS$ and
$\cA$ are orthogonal complementary $|\gE|$-dimensional subspaces of $\cH$,
i.e.~$\cH = \cA \oplus \cS$.  More specifically, an orthogonal basis for $\cA$
is given by
\[
	a_h = e_h - e_{-h}
\]
and an orthogonal basis for $\cS$ is given by
\[
	s_h = e_h + e_{-h}.
\]
So any $f \in \cH$ can be decomposed into $f = (f^a+f^s)/2$, where $f^a \in
\cA$ and $f^s \in \cS$ are defined by
\[
	f^a_h = f_h - f_{-h} \quad \text{and} \quad f^s_h = f_h + f_{-h}.
\]

\subsubsection{Cut space and cycle space} \label{sec:cut-cycle}

We use $\cC$ to denote the subspace of $\cA$ generated by the character
functions
\[
	a_c = \sum_{h \in c} a_h\qquad\mbox{ for } c \in \gC,
\]
known as the \emph{cycle space}.  We use $\cB$ to denote the subspace of $\cA$
generated by the character functions
\[
	a_b = \sum_{h \in b} a_h\qquad\mbox{ for } b \in \gB,
\]
known as the \emph{cut space}.  It is well known \cite[Chapter 14]{godsil2001}
that the cut and cycle spaces are orthogonal complementary subspaces of $\cA$,
i.e.\ $\cA = \cB \oplus \cC$.  The dimension of $\cC$ is $|\gE|-|\gV|+1$, and
the dimension of $\cB$ is $|\gV|-1$.

A cut basis is a set $\gB^* \subset \gB$ such that $(a_b)_{b\in \gB^*}$ form a
basis of the cut space $\cB$.  An explicit cut basis is given as follows: Let
$\gV^*$ be the set of all but one vertices.  Then $\gB^*=\{b(v):v\in\gV^*\}$ is
a canonical cut basis.  In the following, we write
\[
	a_v = a_{b(v)}=\sum_{h\in v} a_h,
\]
for $v \in \gV$.

A cycle basis is a set $\gC^*\subset\gC$ such that $(a_c)_{c\in\gC^*}$ form a
basis of the cycle space $\cC$.  In case $\gG$ has genus one, the cycle space
has dimension $|\gF|+1$ by Euler formula.  An explicit cycle basis is given as
follows.  Let $\gF^*$ be the set of all but one faces.  For $i=1,2$, let
$c_i\in\cC$ be a cycle which is homologous in $\T^2$ to the segment $[0,T_i]$.
Then $\gC^*=\gF^*\cup\{c_1,c_2\}$ is a canonical cycle basis.

\subsubsection{Discrete differential operators}

Let $\cV$ be the vector space of functions $f \colon \gV \to \R$ such that
$\sum_{v \in \gV} f_v = 0$.  For $f\in\cV$, define
\[
	\grad(f)=-\sum_{v\in\gV} f_v a_v=\sum_{h\in\gH}(f_{v(-h)}-f_{v(h)})e_h.
\]
Then $\grad:\cV\to\cB$ is an isomorphism.

For $f = (f_h)_{h \in \gH} \in \cA$, define
\begin{align*}
	\div_b(f) &= (f,a_b)/2 = \sum_{h \in b} f_h, & \div(f)&=(\div_b(f))_{b\in\gB};\\
	\curl_c(f) &= (f,a_c)/2 = \sum_{h \in c} f_h, & \curl(f)&=(\curl_c(f))_{c \in \gC}.
\end{align*}
Because of the orthogonality between $\cB$ and $\cC$, we have
\[
	\ker(\div) = \cC \quad \text{and} \quad \ker(\curl) = \cB.
\]
Hence we can identify
\[
	\operatorname{im}(\div) \simeq \cB
	\quad \text{and} \quad
	\operatorname{im}(\curl) \simeq \cC.
\]
\begin{remark} \label{remark:basis}
	Let $\pi_{\gB}:\R^{|\gB|}\to\R^{|\gV|-1}$ be the projection $(x_b)_{b\in\gB}
	\mapsto (x_b)_{b\in\gB^*}$.  Then
	$\pi_{\gB}:\operatorname{im}(\div)\to\R^{|\gV|-1}$ is an isomorphism. Indeed,
	$(a_v)_{v\in\gV^*}$ is a basis of the cut space so the Gram matrix
	$(a_v,a_{v'})_{v,v'\in V^*}$ is invertible.  In the same way, the projection
	$\pi_{\gC}:\R^{|\gC|}\to\R^{|\gF|+1}$, $(x_c)_{c\in\gC}\mapsto
	(x_{c})_{c\in\gC^*}$ restricts to an isomorphism from
	$\operatorname{im}(\curl)$ to $\R^{|\gC|+1}$.
\end{remark}

\subsection{A divergence over shortest edges}

For each half-edge $h$, let $\cent\ell_h$ be the length of the segment
$\varrho(e(h))$.  For $b \in \gB$, define
\[
	\cent\ell_b= \min_{h \in b} \cent\ell_h\quad\text{and}\quad m(b)=\{ h \in b \mid \cent\ell_h = \cent\ell_b\}.
\]
We define the operator $\mdiv:\cA\to\R^{|\gB|}$ by
\[
	\mdiv_b(\phi)=(a_{m(b)},\phi)/2=\sum_{h\in m(b)}\phi_h, \qquad
	\mdiv(\phi) = (\mdiv_b(\phi))_{b\in\gB}.
\]
In general, $\cB_m := \operatorname{im}(\mdiv)$ is different from $\cB$, but
the following proposition asserts that they have the same dimension.

\begin{proposition} \label{proposition:mdiv}
 	$\mdiv$ has rank $|\gV|-1$.  Moreover, there exists a cut basis $\gB_m^*$ such
 	that $(\mdiv_b)_{b\in\gB_m^*}$ has rank $|\gV|-1$.
\end{proposition}

\begin{proof}
	Let $\phi\in\cB$ such that $\mdiv(\phi)=0$. We can write $\phi=\grad(f)$ with
	$f\in\cV$.  Assume that $f$ is not constant.  Let $\gV_1$ be the set of
	vertices where $f$ achieves maximum and let $b$ be the corresponding cut.
	Then $\mdiv(\phi)<0$, a contradiction. So $f$ is constant and $\phi=0$.
	Hence $\mdiv$ is injective on $\cB$ and
	\[
		\rank(\mdiv)\geq |\gV|-1.
	\]
	For the reverse inequality, consider for $\varepsilon>0$ the operator
	$\mdiv^{\varepsilon}:\cA\to\R^{|B|}$ defined by
	\[
		\mdiv_b^{\varepsilon}(\phi)=\sum_{h\in b}e^{\varepsilon^{-2}(\cent\ell_b-\cent\ell_h)}\phi_h,\quad b\in \gB.
	\]
	Then
	\[
		\mdiv=\lim_{\varepsilon\to 0} \mdiv^{\varepsilon}.
	\]
	But we can write for $\varepsilon>0$
	\[
		\mdiv^{\varepsilon}=\Phi^{\varepsilon}\circ\div\circ \Psi^{\varepsilon},
	\]
	where
	\[
		\Psi^{\varepsilon}(\phi)=(e^{-\varepsilon^{-2}\cent\ell_h}\phi_h)_{h\in \gH}
		\quad\text{and}\quad
		\Phi^{\varepsilon}(X)=(e^{\varepsilon^{-2}\cent\ell_b}X_b)_{b\in\gB}.
	\]
	Hence
	\[
		\rank(\mdiv^{\varepsilon})\leq\rank(\div)=|\gV|-1.
	\]
	Since the rank is lower semi-continuous, it follows that
	\[
		\rank(\mdiv)\leq |\gV|-1.
	\]
	Hence $\mdiv$ has rank $|\gV|-1$. There exists a subset
	$\gB_m^*\subset\gB$ with cardinal $|\gV|-1$ such that $(\mdiv_b)_{b\in\gB_m^*}$
	has rank $|\gV|-1$.  Then for $\varepsilon>0$,
	$(\mdiv_b^{\varepsilon})_{b\in\gB_m^*}$ has rank $|\gV|-1$ by continuity, so
	$(\div_b)_{b\in\gB_m^*}$ has rank $|\gV|-1$ and $\gB_m^*$ is a cut basis.
\end{proof}

\section{Horizontal balance and rigidity} \label{sec:horizontal}

In this part, we assume that $\gG = (\gH, \involute,\rotate, \varrho)$ is a
finite graph represented in a flat torus $\T^2 = \C/\langle T_1, T_2 \rangle$.

To each half-edge $h$ is associated the unit tangent vector
$\cent{u}_h=e^{\ii\cent{\theta}_h}$ of the segment $\varrho(e(h))$ at
$\varrho(v(h))$. Recall that $\cent{\ell}_h$ denotes the length of the segment
$\varrho(e(h))$ and set $\cent{x}_h=\cent{\ell}_h\cent{u}_h$.  As a general
rule, we use a superscript $\circ$ to denote quantities associated to the
given graph $\gG$, which are to be perturbed as parameters in the construction.
\begin{remark}
	In the intuitive picture in Section \ref{sec:graph-intuitive}, $\cent{u}_h$
	is simply the unit vector in the direction of the oriented edge $h$ and
	$\cent{\ell}_h$ is its length.
\end{remark}
Obviously,
$\cent{u} = (\cent{u}_h)_{h \in\gH} \in \cA^2$,
$\cent{\ell} = (\cent{\ell}_h)_{h \in \gH} \in \cS$ and
$\cent{x} = (\cent{x}_h)_{h \in\gH} \in \cA^2$.
For $x\in\cA^2$ in a neighborhood of $\cent{x}$, we define
\[
	u_h(x)=\frac{x_h}{\|x_h\|}
	\quad \text{and} \quad
	u(x)=(u_h(x))_{h\in\gH}\in\cA^2.
\]
We define the horizontal forces as the function
\begin{align*}
	\Fhor\colon \cA^2 &\to \cB^2\\
	x &\mapsto \div(u(x)).
\end{align*}
More explicitly, for any cut $b\in\gB$
\[
	\Fhor_b(x)=\sum_{h\in b}\frac{x_h}{\|x_h\|}.
\]

\begin{definition}\label{def:horblance}
	The graph $\gG$ is balanced if $\Fhor(\cent{x}) = 0$.
\end{definition}

\begin{remark} \label{remark:vertex-cuts}
	By remark \ref{remark:basis}, the graph is balanced if and only if
	\[
		\sum_{h\in v}\frac{\cent{x}_h}{\|\cent{x}_h\|}=0
	\]
	for all $v\in V^*$. In other words, it suffices to consider horizontal forces
	on a canonical cut basis consisting of all but one vertex cuts.
\end{remark}

\begin{remark}
	A balanced graph is a weak local minimal network in the sense of Ivanov and
	Tuzhilin~\cite{ivanov1994}; see Proposition~\ref{prop:minlength} in the
	Appendix.
\end{remark}

Now define the horizontal periods as the operator
\begin{align*}
	\Phor\colon \cA^2 &\to \cC^2\\
	x &\mapsto \curl(x).
\end{align*}
More explicitly, for any cycle $c\in\gC$
\[
	\Phor_c(x)=\sum_{h\in c} x_h.
\]
The graph $\gG$ is represented on the torus $\T^2$ so we have on any canonical
cycle basis $\gC^*$
\begin{equation}\label{eq:Hperiod}
	\Phor_c(\cent{x}) = \begin{cases}
		0, &c \in \gF^*;\\
		T_i, &c=c_i,\quad i=1,2.
	\end{cases}
\end{equation}

\begin{remark}\label{rmk:kirchhoff}
	The balance and period equations can be compared to Kirchhoff's current and
	voltage laws of electrical networks.  More specifically, $\cent{u}$,
	$\cent{x}$, and $\cent{\ell}$ play, respectively, the roles of currents,
	voltages and resistance.
\end{remark}

Recall from Section \ref{sec:vector-spaces} that $\operatorname{ker}(\Phor)=\cB^2$.
\begin{definition} \label{definition:horizontal-rigidity}
	The graph $\gG$ is rigid if the differential $D\Fhor(\cent{x})$ restricted to
	$\operatorname{ker}(\Phor)$ is an isomorphism from $\cB^2$ to $\cB^2$.
	Equivalently, the graph is rigid if
	\[
		(D\Fhor(\cent{x}),\Phor):\cA^2\to\cA^2
	\]
	is an isomorphism.
\end{definition}
A computation reveals that
\begin{equation} \label{eq:DFhor}
	D\Fhor_b(\cent{x})\cdot \chi=\sum_{h\in b}\frac{1}{\|\cent{x}_h\|^3}\big(
	\langle\cent{x}_h,\cent{x}_h\rangle \chi_h-\langle\cent{x}_h,\chi_h\rangle\cent{x}_h\big).
\end{equation}

The following rigidity result was proved in \cite[Theorem 10]{younes2009}.
\begin{theorem} \label{theorem:triangles}
	Assume that all faces of $\gG$ have 2 or 3 edges. Then $\gG$ is rigid.
\end{theorem}
As \cite{younes2009} is not published, we include a proof in Appendix
\ref{appendix:triangles}.

Next we need to assign a saddle tower to each vertex.
\begin{definition}
	Let $v \in \gV$. Consider the convex closed polygon $P$ whose edges are the
	unit vectors $\cent{u}_h$ for $h\in v$ in the order given by the permutation
	$\rotate$.  We say that $v$ is \emph{ordinary} if $P$ is neither
	\emph{degenerate} (appears as a line segment) nor \emph{special} (appears as
		a parallelogram with two sides of length 1 and two sides of length $\geq
	2$).
\end{definition}
Note that non-ordinary vertices occur only in the presence of parallel edges.
If $v$ is ordinary, there exists a saddle tower $\sS_v$, unique up to
translations, with $\deg(v)$ wings which are in correspondence with the
half-edges $h\in v$, so that the direction of the wing corresponding to $h$ is
$\cent{u}_h$ and the natural order on the wings is given by the permutation
$\rotate$ (see Section \ref{sec:jenkins-serrin}).

We want to glue the wings of these saddle towers along the edges of $\gG$.
This construction is achieved if the saddle towers share a horizontal
reflection plane, as stated in the following theorem due to
Younes~\cite{younes2009}.

\begin{theorem}[TPMSs with horizontal symmetry]\label{thm:rami}
 	Let $\gG$ be a graph represented in $\T^2 = \C/\langle T_1, T_2 \rangle$, and
 	assign a phase $\phi_v \in \{0, \pi\}$ to each vertex $v \in \gV$.  If $\gG$ is
 	orientable, balanced, rigid, and all vertices are ordinary, then for
 	sufficiently small $\varepsilon>0$, there is a family $\sM_\varepsilon$ of
 	embedded minimal surfaces of genus $|\gF| + 1$ in the flat 3-torus
	\[
		\T_\varepsilon^3=\R^3/\langle (T_1\varepsilon^{-2},0),
		(T_2\varepsilon^{-2},0), (0,0,2\pi)\rangle.
	\]
	They lift to triply periodic minimal surfaces $\widetilde\sM_\varepsilon$ in
	$\R^3$ such that
 	\begin{enumerate}
		\item $\widetilde{\sM}_{\varepsilon}$ converges, after scaling by
			$\varepsilon^2$, to $\widetilde{\gG}\times\R$ as $\varepsilon\to 0$,
			where $\widetilde{\gG}$ is the lift of $\gG$ to $\R^2$.

		\item For each vertex $v$ of $\widetilde\gG$, there exists a horizontal
			vector $X_v(\varepsilon)$ such that
			$\widetilde\sM_\varepsilon-X_v(\varepsilon)$ converges on compact subset
			of $\R^3$ to a saddle tower $\sS_v$ as $\varepsilon \to 0$. Moreover,
			$\varepsilon^2X_v(\varepsilon) \to \widetilde\varrho(v)$ as
			$\varepsilon\to 0$, where $\widetilde\varrho$ is the lift of $\varrho$.

		\item Each limit saddle tower $\sS_v$ has phase $\phi_v$.  Moreover,
			$\widetilde\sM_{\varepsilon}$ is symmetric with respect to a horizontal
			plane.
	\end{enumerate}
\end{theorem}

\begin{figure}
	\begin{center}
		\includegraphics[height=3cm]{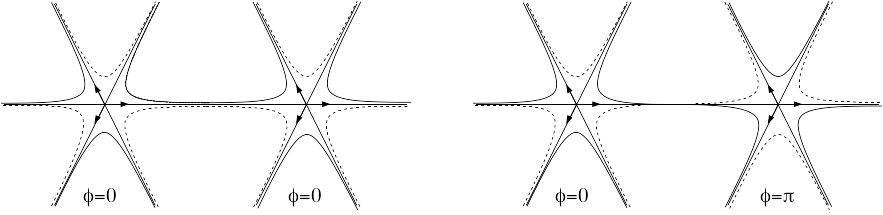}
	\end{center}
	\caption{Left: two saddle towers in phase. Right: two saddle towers in
	opposite phase. The solid and dotted lines represent the level lines $x_3=0$
and $x_3=\pi$, respectively.}
\end{figure}

\section{Closer look on Karcher--Scherk saddle towers}
\label{sec:towers}

We want to break the horizontal symmetry by prescribing phase difference
between adjacent saddle towers. The phase differences must satisfy a balancing
condition that involves higher order shapes of the Scherk ends.  This section
is dedicated to a detailed investigation on the saddle towers.  Quantities
describing the shapes of Scherk ends will be defined in term of Weierstrass
representation, and will be useful in the formulation of the main theorem.

\subsection{Weierstrass parameterization of saddle towers}\label{sec:weierstrass-1tower}

The quotient of a saddle tower by its vertical period is conformally equivalent
to a Riemann sphere $\C \cup \{\infty\}$ with $n$ punctures $p_h$, $h \in \gH$,
corresponding to the $n$ wings.  The punctures must lie on a circle $C$ fixed
by the anti-holomorphic involution $\rho$ corresponding to the horizontal
reflections.  For convenience, $C$ is often taken to be the unit circle or the
real line.

Recall that our saddle tower is scaled so that it has vertical period $2\pi$.
Then its Weierstrass data can be written on the punctured Riemann sphere as
\begin{equation}\label{eq:weierstrassz}
	\Phi_1 = \sum_{h \in \gH} \frac{-\cos\theta_h}{z- p_h} dz,
	\qquad
	\Phi_2 = \sum_{h \in \gH} \frac{-\sin\theta_h}{z- p_h} dz,
	\qquad
	\Phi_3 = \sum_{h \in \gH} \frac{-\ii \sigma_h}{z- p_h} dz,
\end{equation}
where the orientation $\sigma_h=\pm 1$ and $\sigma_{\rotate(h)}=-\sigma_h$ for
$h\in\gH$.  The conformality condition
\[
	\Phi_1^2+\Phi_2^2+\Phi_3^2=0
\]
determines the punctures $p_h$ up to a M\"obius transformation (although doing
this explicitly can be difficult).  The saddle tower is parameterized by the
Weierstrass Representation formula
\begin{equation} \label{eq:weierstrass-representation}
	z \mapsto \re \int_{z_0}^z \Phi =
	\re \int_{z_0}^z (\Phi_1, \Phi_2, \Phi_3).
\end{equation}

The stereographically projected Gauss map $G = - (\Phi_1 + \ii\Phi_2)/\Phi_3$
extends holomorphically to the punctures with $G(p_h) = \ii\sigma_h u_h$.
Then the Gauss map extends at the end $p_h$ with
\[
	N(p_h) = \sigma_h(-\sin(\theta_h), \cos(\theta_h),0).
\]
Consequently, the image of the arc between $p_h$ and $p_{\rotate(h)}$ is a
$0$-arc if and only if $\sigma_h=1$.

\subsection{Shape of wings}\label{sec:undulation}

Let $w_h$ be a local complex coordinate in a neighborhood of $p_h$ with
$w_h(p_h)=0$.  We define
\begin{align*}
	\Upsilon_h &= \sigma_h\Big\langle N(p_h), \Res{\frac{\Phi}{w_h}}{p_h} \Big\rangle_{\!\text{H}}\\
	&= -\sin(\theta_h) \Res{\frac{\Phi_1}{w_h}}{p_h} + \cos(\theta_h) \Res{\frac{\Phi_2}{w_h}}{p_h},
\end{align*}
where $\langle\cdot,\cdot\rangle_\text{H}$ denotes the hermitian scalar product
on $\C^3$ (semi-linear on the left).  Note that $\Upsilon_h$ depends on the
local coordinate $w_h$.

\begin{proposition} \label{proposition:Upsilon}
	\[
		\Upsilon_h=\ii\frac{dG}{G\,dw_h}(p_h).
	\]
\end{proposition}
\begin{proof}
	We may expand $G$ and $\Phi_3$ around $p_h$ as
	\begin{align*}
		G &=  \ii \sigma_h e^{\ii\theta_h} \big( 1 + a w_h + O(w_h^2) \big),\\
 		G^{-1} &= -\ii\sigma_h e^{-\ii\theta_h} \big( 1 - a w_h + O(w_h^2) \big),\\
 		\Phi_3 &= -\ii\sigma_h\big( \frac{1}{w_h} + b + O(w_h) \big)dw_h,
	\end{align*}
	where $a,b\in\C$ and
	\[
 		a=\frac{dG}{G\,dw_h}(p_h).
 	\]
 	This gives
 	\begin{align*}
 	 	\Res{\frac{\Phi_1}{w_h}}{p_h}&= - \cos(\theta_h)b - \ii\sin(\theta_h)a,\\
 	 	\Res{\frac{\Phi_2}{w_h}}{p_h}&= - \sin(\theta_h)b + \ii\cos(\theta_h)a,
 	\end{align*}
 	and
	\[
		\Upsilon_h=-\sin(\theta_h)\left(- \cos(\theta_h)b - \ii\sin(\theta_h)a\right) + \cos(\theta_h)\left(- \sin(\theta_h)b + \ii\cos(\theta_h)a\right)=\ii a.
	\]
\end{proof}

We define the quantities
\begin{align*}
	\mu_h & = \lim_{z \to p_h} \Big(e^{\ii\theta_h} \log |w_h(z)| + \re\int_{z_0}^{z} \Phi_1 + \ii \re\int_{z_0}^{z}\Phi_2 \Big),\\
	\nu_h & = \lim_{z \to p_h} \Big( -\sigma_h \arg(w_h(z)) + \re\int_{z_0}^{z} \Phi_3 \Big).
\end{align*}
It is easy to see that the limits exist, but these quantities depend on the
local coordinate $w_h$ and the choice of the base point $z_0$.

Recall that a saddle tower has a horizontal symmetry plane corresponding to an
anti-holomorphic involution $\rho$ of the Riemann sphere that fixes all
punctures $p_h$ for $h\in\gH$.
\begin{definition}
	The local coordinate $w_h$ is said to be \emph{adapted} if
	$w_h\circ\rho=\overline{w_h}$ and $w_h>0$ near $p_h$ on the arc between $p_h$
	and $p_{\rotate(h)}$.
\end{definition}

\begin{proposition} \label{prop:adapted}
	Assume that the coordinate $w_h$ is adapted and let $\phi$ be the phase of the
	saddle tower.  Then $\Upsilon_h > 0$ and
	\[
		\nu_h=\begin{cases}
			\phi & \pmod{2\pi} \quad \text{if } \sigma_h=+1,\\
			\phi+\pi & \pmod{2\pi} \quad \text{if } \sigma_h=-1.
		\end{cases}
	\]
\end{proposition}
\begin{proof}
	Assume that the coordinate $w_h$ is adapted.  We start by computing $\nu_h$.
	Consider $z$ on the arc between $p_h$ and $p_{\rotate(h)}$.  We have
	$\arg(w_h)=0$ for $z$ sufficiently close to $p_h$.  Then, by the definition
	of the phase $\phi$, we have
	\[
		\nu_h = \re\int_{z_0}^z \Phi_3 = \begin{cases}
			\phi & \pmod{2\pi} \quad \text{if } \sigma_h=+1,\\
			\phi+\pi & \pmod{2\pi} \quad \text{if } \sigma_h=-1.
		\end{cases}
	\]
	Since the coordinate $w_h$ is adapted, we have, for $i=1,2$,
	\[
		\rho^*\Big(\frac{\Phi_i}{w_h}\Big)=\frac{\overline{\Phi_i}}{\overline{w_h}},
	\]
	so $\Res{\Phi_i/w_h}{p_h}\in\R$ and $\Upsilon_h\in\R$.

	To understand its sign we need to go back to the Jenkins-Serrin construction.
	The solution of the Jenkins-Serrin problem is a graph on the convex domain
	bounded by the polygon $P$ so its Gauss map is non-horizontal in the interior
	of $P$. Recalling that the conjugate minimal surface has the same gauss map,
	we have $|G|\neq 1$ on the Riemann sphere minus the circle $C$ fixed by the
	symmetry $\rho$, and $|G|=1$ on $C$.  Hence in a neighborhood of $p_h$, $\log
	G$ is a well-defined holomorphic function and is pure imaginary only on
	$C$, so its zero at $p_h$ is simple by the local behavior of holomorphic
	functions in a neighborhood of a zero.  Hence $\Upsilon_h\neq 0$ by
	Proposition \ref{proposition:Upsilon}.

  Each horizontal symmetry
	curve on a saddle tower is a convex curve. Since the unitary
	vectors $u_h$ are ordered in the counterclockwise order, the argument of $G$
	is non-increasing on the arc between $p_h$ and $p_{\rotate(h)}$. This implies
	$\Upsilon_h\geq 0$ by Proposition \ref{proposition:Upsilon}.
\end{proof}

For a geometric intuition of the quantities $\Upsilon_h$, $\mu_h$, and $\nu_h$,
let us expand the Weierstrass parametrization around $p_h$:
\begin{alignat*}{4}
	\re \int_{z_0}^{z} \Phi_1 &=
	\overbrace{\vphantom{\bigg(\bigg)}\re \mu_h - \cos\theta_h\log|w_h(z)|}^\text{planar terms}
	&&
	+ \overbrace{\re \left(w_h(z) \Res{\frac{\Phi_1}{w_h}}{p_h}\right)}^\text{undulation terms}
	&&+ O(|w_h(z)|^2),\\
	\re \int_{z_0}^{z} \Phi_2 &= \im \mu_h - \sin\theta_h\log|w_h(z)| &&
	+ \re \left(w_h(z) \Res{\frac{\Phi_2}{w_h}}{p_h}\right)
	&&+ O(|w_h(z)|^2),\\
	\re \int_{z_0}^{z} \Phi_3 &= \phantom{im} \nu_h     &&+\sigma_h \arg(w_h(z))
	&&+ O(|w_h(z)|).
\end{alignat*}
We then observe that
\begin{itemize}
	\item The terms under the first brace describe a vertical half-plane that
		passes through $(\mu_h, \nu_h)$ and extends in the direction $\theta_h$ as
		$w_h(z) \to 0$.

	\item The terms under the second brace describe a sinusoidal undulation in
		the vertical coordinate that decays exponentially with the horizontal
		distance from $(\mu_h, \nu_h)$; see Figure~\ref{fig:tower}.

	\item The quantity $\Upsilon_h$ describes the ``initial'' amplitude of the
		undulation on the $0$-plane in the asymptotic normal direction of the wing.
\end{itemize}

\begin{remark} \label{remark:change-of-coordinate}
 	The quantities $\Upsilon_h$ and $\mu_h$ behave as follows under change of
 	coordinate: if $\tilde w_h$ is another complex coordinate in a neighborhood
 	of $p_h$ and $\tilde\Upsilon_h$, $\tilde\mu_h$ denote the corresponding
 	quantities, we have, from the definitions,
 	\[
 		\Upsilon_h=\tilde\Upsilon_h\frac{d\tilde w_h}{d w_h}(p_h)
 		\quad \text{and} \quad
 		\tilde\mu_h-\mu_h=e^{i\theta_h}\log\left|\frac{d\tilde w_h}{d w_h}(p_h)\right|.
	\]
	Also observe that $\frac{d\tilde w_h}{d w_h}(p_h)>0$ if
 	both coordinates are adapted.
\end{remark}

\subsection{Examples}\label{sec:exampletower}

In this section, we compute the quantities $\Upsilon_h$ and $\mu_h$ for certain
saddle towers with explicitly known Weierstrass data.  The punctures $p_h$ will
be placed on the unit circle in the clockwise order, so we use the adapted
coordinate
\[
	w_h=\ii\,\frac{z-p_h}{z+p_h}.
\]
With this choice, we have for $h\in\gH$
\begin{equation} \label{eq:muh-circle}
	\mu_h=-e^{\ii\theta_h}\log 2-\sum_{j\neq h}e^{\ii\theta_j}\log|p_h-p_j|.
\end{equation}

\subsubsection{Symmetrically deformed saddle towers} \label{sec:symtower}

These examples are described in \cite[\S 2.4.1]{karcher1988}.  They have $n=2k$
ends with $k\geq 2$. Their Weierstrass data are given by
\[
	G=z^{k-1}
	\quad\text{and}\quad
	\Phi_3=\frac{2k\sin(k\varphi)}{z^k+z^{-k}-2\cos(k\varphi)}\cdot \frac{dz}{z}
\]
where $0<\varphi<\pi/k$.  Comparing to~\cite[\S 2.4.1]{karcher1988}, we
multiplied $\Phi_3$ by $2k\sin(k\varphi)$ so that the vertical period is
$2\pi$.

The punctures are at
\[
	p_h = \exp\left(-\ii\Big\lfloor \frac{h}{2} \Big\rfloor \frac{2\pi}{k} + \ii(-1)^h \varphi\right),
	\quad 1\leq h\leq n,
\]
we have $\sigma_h=(-1)^h$ and the direction of the wings are
\[
	\theta_h = \Big\lfloor \frac{h}{2} \Big\rfloor \frac{2\pi}{k} - (-1)^h \psi, \quad \psi - \frac{\pi}{n} = (k-1) \Big( \frac{\pi}{n} - \varphi \Big).
\]
The case $n=4$ gives Scherk's surfaces. The case $\psi = \varphi=\pi/n$ gives
the most symmetric saddle towers.  The saddle towers are embedded for $0 \leq
\psi \leq \pi/k$ (with strict inequalities for $n=4$).  If $n\geq 6$, the limit
cases $\psi = 0$ and $\psi=\pi/k$ give a saddle tower with $k$ pairs of
parallel wings; an example with $n=10$ is illustrated on right side of
Figure~\ref{fig:Matthias}.

Given our choice of adapted coordinate and the simple form of the Gauss map,
Proposition~\ref{proposition:Upsilon} immediately gives
\[
	\Upsilon_h = n-2.
\]
As for $\mu_h$, Equation \eqref{eq:muh-circle} does not simplify very much in
general, so we only present some special cases.  When $n=4$,
\eqref{eq:muh-circle} simplifies to
\[
	\mu_h=e^{-\ii\theta_h}\log\tan(\psi), \quad 1\leq h\leq 4,
\]
which is not collinear to $e^{\ii\theta_h}$ unless $\psi=\pi/4$ (see Remark
\ref{remark:tower}).  For arbitrary $n$, if $\psi=\pi/n$, then by symmetry,
$\mu_h$ is collinear to $e^{\ii\theta_h}$, and its norm (independent of $h$) is
tabulated below for small values of $n$.
\begin{center}
	\begin{tabular}{c|c}
		$n$ &$e^{-\ii\theta_h}\mu_h$\\
		\hline
		4 & 0\\
		6 & $\log\sqrt{3}$\\
		8 & $\sqrt{2}\log(1+\sqrt{2})$\\
		10 & $\frac{1}{4}\log 5+\frac{\sqrt{5}}{2}\log(2+\sqrt{5})$\\
		12 & $\frac{1}{2}\log 3 +\sqrt{3}\log(2+\sqrt{3})$\\
	\end{tabular}
\end{center}

\subsubsection{Isosceles saddle tower with 6 wings}

These examples with $n=6$ wings are described in \cite[\S 2.5.1]{karcher1988}.
Their Weierstrass data are given by
\[
	G=\frac{z^2+r}{1+rz^2}
	\quad\text{and}\quad
	\Phi_3=\frac{8\cos(\varphi)^2}{(1-r)^2}\cdot\frac{1+r^2+r\cdot(z^2+z^{-2})}{(z+z^{-1})\cdot(z^2+z^{-2}-2\cos(2\varphi))}\cdot\frac{dz}{z},
\]
where $r\in(-1,1)$ is the unique solution of
$$\frac{4r}{(r-1)^2}=\frac{2\sin(\varphi)-1}{\cos(\varphi)^2}.$$
Comparing
to~\cite[\S 2.5.1]{karcher1988}, we multiplied $\Phi_3$ by
$8\cos(\varphi)^2/(1-r)^2$ so that the vertical period is $2\pi$.

The punctures are at
\[
	(p_0,\cdots,p_5)=(e^{-\ii\varphi},-\ii,-e^{\ii\varphi},-e^{-\ii\varphi}, \ii,
	e^{\ii\varphi}),
\]
we have $\sigma_h=(-1)^{h+1}$ and the directions of the wings are
\[
	(\theta_0,\cdots,\theta_5)
	=(\psi, \pi/2, \pi-\psi,-\pi+\psi,-\pi/2,-\psi),
\]
where $\psi\in(0,\pi/2)$ is the solution of $\sin\psi+\sin\varphi = 1$,
so the wings are parallel to the sides of an isosceles triangle.  The most
symmetric saddle tower is recovered by $\psi = \pi/6$.  An example with $\psi =
\pi/3$ is illustrated on the left side of Figure~\ref{fig:Matthias}.  The
Jenkins--Serrin polygon is degenerate in the limit $\psi \to \pi/2$, and
special in the limit $\psi \to 0$.

We compute explicitly
\[
	r=\frac{\cos(\psi)-\cos(\varphi)}{\cos(\psi)+\cos(\varphi)}
\]
\[
	\Upsilon_h = \begin{cases}
		4\cos(\psi)/\cos(\varphi) & h = 1 \pmod{3},\\
		4\cos(\psi)/\sin(2\varphi) & h \ne 1 \pmod{3},
	\end{cases}
\]
\[
	\mu_0 = \ii\, \log\frac{\cos(\varphi)}{1-\sin(\varphi)} + \exp^{-\ii\psi} \log\cot(\varphi), \quad
	\mu_1 = 2\ii\, \sin(\psi) \log\frac{\cos(\varphi)}{1-\sin(\varphi)},
\]
and the others $\mu_h$ can be obtained by symmetry, namely
\[
	\mu_0 = -\overline{\mu_2} = -\mu_3 = \overline{\mu_5}
	\quad\text{and}\quad
	\mu_1 = -\mu_4.
\]

\subsection{Rigidity} \label{sec:rigid-tower}

By a result of Cosin-Ros \cite{cosin2001}, all saddle towers are rigid, in the
sense that the space of bounded Jacobi fields on a saddle tower is
3-dimensional and consists of translations.  This means that when the angles
$\theta_h$ are fixed, a saddle tower admits no deformation other than
translations.  In this section, we reformulate this result in term of
Weierstrass Representation, in a way that can be used in our gluing
construction.

The Weierstrass data of a saddle tower can always be written as in Equation
\eqref{eq:weierstrassz}. The equation to solve is $Q=0$, where
$Q=\Phi_1^2+\Phi_2^2+\Phi_3^2$.
Note that $Q$ has at most simples poles at the punctures $p_h$.
The angles $\theta_h$ are fixed, and the
unknowns are the poles $p_h$ for $h\in\gH$. We are given a solution, denoted by
$\cent{p}_h$, and we want to study its infinitesimal deformations. The
corresponding Weierstrass data is denoted
$(\cent{\Phi}_1,\cent{\Phi}_2,\cent{\Phi}_3)$.

We formulate the equation $Q=0$ as follows.  Without loss of generality, we may
assume by rotation that $\cos(\theta_h)\neq 0$ for all $h\in\gH$ and all zeros
of $\cent{\Phi}_1$ are simple.  Then for $p$ in a neighborhood of $\cent{p}$,
$\Phi_1$ has $n-2$ simple zeros $\zeta_1,\cdots,\zeta_{n-2}$ which depend
holomorphically on $p$. The meromorphic 1-form $Q/\Phi_1$ is holomorphic at
$p_h$, $h\in\gH$, and has (at most) simple poles at $\zeta_1, \cdots,
\zeta_{n-2}$.  We define
\[
	\Lambda(p)=\left(\Res{\frac{Q}{\Phi_1}}{\zeta_i}\right)_{1\leq i\leq n-3}.
\]
By the Residue Theorem, $Q=0$ is equivalent to $\Lambda(p)=0$.  By a M\"obius
transformation, we may fix the value of three points $p_h$, so the parameter
$p$ lies in a space of complex dimension $n-3$.
\begin{theorem} \label{thm:rigid-tower}
	The differential $D\Lambda(\cent{p})$ is an isomorphism.
\end{theorem}
This theorem is proved in Annexe A of \cite{younes2009}, which is unfortunately
not published.  So we include a proof in Appendix \ref{app:rigid-tower}.

\section{Main result} \label{sec:main}
\subsection{Vertical balance and rigidity}
We prescribe the
phase differences between adjacent saddle towers through an antisymmetric
\emph{phase function} $\cent\phi \colon \gH \to \R/2\pi\Z$ that assigns a phase
difference $\cent\phi_h$ to each half-edge $h$.  We say that $\cent\phi$ is trivial if
$\cent\phi = 0$ or $\cent\phi=\pi$ on every half-edge, which is the case in
Theorem~\ref{thm:rami}.

Define the vertical periods as the function
\begin{align*}
	\Pver\colon \cA &\to \cC\\
	\phi &\mapsto \curl(\phi).
\end{align*}
We require that the periods of the phase function are given on a canonical
cycle basis $\gC^*$ as
\begin{equation}\label{eq:Vperiod}
	\Pver_c(\cent\phi) = \sum_{h\in c}\cent\phi_h= \begin{cases}
		0, & c \in \gF;\\
		\Psi_i, & c=c_i,\quad i=1,2.
	\end{cases}
\end{equation}
for some $\Psi_1, \Psi_2 \in \R/2\pi\Z$.  We call $\Psi_1$ and $\Psi_2$ the
\emph{fundamental shifts}.  We want to construct minimal surfaces in the flat
3-torus
\[
	\T_\varepsilon^3=\R^3/\langle (\Lambda_1 + T_1\varepsilon^{-2},\Psi_1),
	(\Lambda_2 + T_2\varepsilon^{-2},\Psi_2), (0,0,2\pi)\rangle,
\]
Here $\Lambda_1$, $\Lambda_2$ are fixed complex numbers that prescribe a first
order horizontal deformation of the lattice as $\varepsilon$ varies.

\medskip

The phase function must satisfy a balancing condition, which we now explain.

For each vertex $v \in \gV$, let $\hat\C_v$ be the punctured Riemann sphere on
which the saddle tower $\sS_v$ is parametrized.  Fix an adapted local
coordinate $w_h$ in a neighborhood of the puncture $p_h \in \C_{v(h)}$ for
every $h \in \gH$.  Recall the definition of the numbers $\Upsilon_h$, $\mu_h$
in Section \ref{sec:towers} and the notation $\mu^a_h=\mu_h-\mu_{-h}$.  If the
graph $\gG$ is balanced and rigid, the system
\begin{equation}\label{eq:system}
	\begin{cases}
		D\Fhor_b(\cent{x}) \cdot \xi = 0, & b \in \gB;\\
		\Phor_c(\xi) = - \Phor_c(\mu^a), & c \in \gF^*;\\
		\Phor_{c_i}(\xi) = - \Phor_{c_i}(\mu^a) + \Lambda_i, &  i = 1,2\\
	\end{cases}
\end{equation}
has a unique solution $\xi \in \cA^2$ by Definition
\ref{definition:horizontal-rigidity} and Remark \ref{remark:basis}.

\begin{remark}
	The system~\eqref{eq:system} is invariant by horizontal translations of the
	saddle towers.  Indeed, if $\sS_v$ is translated by a horizontal vector
	$X_v$, then $X_{v(h)}$ is added to $\mu_h$, so $\curl(\grad(X)) = 0$ is added
	to the right-hand side of~\eqref{eq:system}.
\end{remark}

We define a symmetric function
\begin{equation} \label{equation:Kh}
	K_h = \Upsilon_h \Upsilon_{-h} e^{-\re(\xi_h\overline{\cent{u}_h})}.
\end{equation}

By Proposition \ref{prop:adapted}, we have $K_h>0$.  We will see in
Proposition \ref{proposition:Kh} that $K_h$ is independent of the choice of
adapted coordinates $w_h$. Both $\xi$ and $K_h$ depend on $\Lambda_1$ and
$\Lambda_2$, but the dependence is omitted for simplicity.  When the values
of $\Lambda_1$ and $\Lambda_2$ matter, but are not specified in the context,
it is implied that $\Lambda_1 = \Lambda_2 = 0$.

We define the vertical forces as the function
\begin{align*}
	\Fver \colon \cA &\to \cB_m \\
	\phi &\mapsto \mdiv\big( (K_h \sin\phi_h)_{h \in \gH}\big).
\end{align*}
More explicitly, for any cut $b\in\gB$
\[
	\Fver_b(x)=\sum_{h\in m(b)}K_h \sin(\phi_h).
\]
\begin{definition}
	The phase function $\cent\phi$ is balanced if $\Fver(\cent\phi) = 0$.
\end{definition}

\begin{remark}
	Trivial phase functions are trivially balanced.
\end{remark}

\begin{remark}
	Unlike horizontal balancing (see Remark \ref{remark:vertex-cuts}), the
	equation $\Fver(\cent\phi)=0$ is in general not equivalent to
	$\Fver_{b(v)}(\cent\phi)=0$ for $v\in\gV^*$: it is not enough to consider
	vertex cuts; see Example~\ref{ex:necessarycuts}.  This is the reason why it
	is necessary to introduce the whole cut space to define vertical balancing.
\end{remark}

\begin{remark}\label{rmk:continuous}
	In general, the vertical forces do not depend continuously on the horizontal
	periods $T_1$ and $T_2$, but they depend continuously on the deformations
	$\Lambda_1$ and $\Lambda_2$.
\end{remark}

\begin{definition}
	The phase function $\cent\phi$ is rigid if the differential
	$D\Fver(\cent\phi)$ restricted to $\cB$ is an isomorphism between $\cB$ and
	$\cB_m$.  Equivalently, the phase function is rigid if
	$(D\Fver(\cent\phi),\Pver)$ is an isomorphism.
\end{definition}

We call the pair $(\gG, \cent\phi)$ a \emph{configuration}.  We say that the
configuration is horizontally balanced (resp.\ rigid) if the graph is balanced
(resp.\ rigid), and vertically balanced (resp.\ rigid) if the phase function is
balanced (resp.\ rigid).  And we say that the configuration is balanced (resp.\
rigid) if it is both horizontally and vertically balanced (resp.\ rigid).  Our
main result for TPMSs is the following.

\begin{theorem}[TPMSs]\label{thm:main}
 	Let $(\gG, \cent\phi)$ be a configuration, where the graph $\gG$ is represented in
 	$\T^2 = \C/\langle T_1, T_2 \rangle$, and the fundamental shifts of $\cent\phi$ is
 	$\Psi_1$ and $\Psi_2$.  Assume that $\gG$ is orientable, that the
 	configuration is balanced and rigid, and that all vertices are ordinary.
 	Then for sufficiently small $\varepsilon>0$, there is a family
 	$\sM_\varepsilon$ of embedded minimal surfaces of genus $|\gF| + 1$ in the
 	flat 3-torus
	\[
		\T_\varepsilon^3=\R^3/\langle (\Lambda_1 + T_1\varepsilon^{-2},\Psi_1),
		(\Lambda_2 + T_2\varepsilon^{-2},\Psi_2), (0,0,2\pi)\rangle.
	\]
	They lift to triply periodic minimal surfaces $\widetilde\sM_\varepsilon$ in
	$\R^3$ such that
 	\begin{enumerate}
		\item $\widetilde{\sM}_{\varepsilon}$ converges, after a scaling by
			$\varepsilon^2$, to $\widetilde{\gG}\times\R$ as $\varepsilon\to 0$,
			where $\widetilde{\gG}$ is the lift of $\gG$ to $\R^2$.

		\item For each vertex $v$ of $\widetilde\gG$, there exists a horizontal
			vector $X_v(\varepsilon)$ such that
			$\widetilde\sM_\varepsilon-X_v(\varepsilon)$ converges on compact subset
			of $\R^3$ to a saddle tower $\sS_v$ as $\varepsilon \to 0$. Moreover,
			$\varepsilon^2X_v(\varepsilon) \to \widetilde\varrho(v)$ as
			$\varepsilon\to 0$, where $\widetilde\varrho$ is the lift of $\varrho$.

		\item For each half-edge $h$, the phase difference of $\sS_{v(-h)}$ over
			$\sS_{v(h)}$ is equal to $\cent\phi_h$.
	\end{enumerate}
\end{theorem}

\begin{remark}
	Unfortunately, Theorem \ref{thm:main} does not contain Theorem \ref{thm:rami}
	as a particular case. A trivial phase function is trivially balanced, but it
	is not necessarily rigid.  However, Proposition
	\ref{proposition:zero-is-rigid} below implies that the phase function is
	rigid when $\cent\phi_h$, $h \in \gH$, are all $0$ or all $\pi$.  So Theorem
	\ref{thm:rami} follows from Theorem \ref{thm:main} when all saddle towers are
	in-phase, or all adjacent saddle towers are anti-phase.
\end{remark}

\begin{remark}
	In fact, we construct a continuous family locally parameterized by
	$\varepsilon$, $\Lambda_{1,2}$, and $\Psi_{1,2}$.  By
	Proposition~\ref{prop:normalize-lambda} below, we may assume that $\Lambda_1
	= 0$ up to a scaling and a horizontal rotation.  This is therefore a
	5-parameter family up to Euclidean rotations and scalings, in correspondence
	with the deformations of the lattice.
\end{remark}

\subsection{Some auxiliary results} \label{sec:auxiliary}

\begin{proposition} \label{proposition:Kh}
	The constant $K_h$ defined in Equation \eqref{equation:Kh} is independent of
	the adapted local coordinates $w_h$.
\end{proposition}

\begin{proof}
	Consider another adapted local coordinates $\tilde w_h$ in a neighborhood of
	$p_h$ for each $h \in \gH$.  We use a tilde for all quantities associated to
	the coordinate $\tilde w_h$.  By Remark \ref{remark:change-of-coordinate}, we
	have, writing $\kappa_h = \frac{d \tilde w_h}{d w_h}(p_h)>0$, that
	\[
		\tilde \mu_h - \mu_h = \cent{u}_h \log \kappa_h
		\quad\text{and}\quad
		\tilde\Upsilon_h = \Upsilon_h/\kappa_h.
	\]
	Therefore, since $\cent{u}_{-h}=-\cent{u}_h$
	\[
		\tilde\mu_h^a-\mu_h^a=\cent{u}_h\log(\kappa_h\kappa_{-h}).
	\]
	Observe that $\tilde\xi - \xi$ is the solution of
	\[
		\begin{cases}
			D\Fhor(\cent{x}) \cdot (\tilde\xi - \xi) = 0;\\
			\Phor(\tilde\xi - \xi) =  - \Phor(\tilde\mu^a - \mu^a).
		\end{cases}
	\]
	By Equation \eqref{eq:DFhor}, we have
	\[
		D\Fhor(\cent{x})\cdot(\tilde\mu_h^a-\mu_h^a)=0.
	\]
	So the solution is trivially
	\[
		\tilde\xi_h -\xi_h = -\tilde\mu_h^a+\mu_h^a=-\cent{u}_h \log(\kappa_h \kappa_{-h}), \quad h \in \gH.
	\]
	Therefore,
	\[
		\frac{\tilde K_h}{K_h} = \frac{\tilde\Upsilon_h}{\Upsilon_h}\frac{\tilde\Upsilon_{-h}}{\Upsilon_{-h}}\exp\left(-\re\left((\tilde\xi_h-\xi_h)\overline{\cent{u}_h}\right)\right)
		=\frac{e^{\log(\kappa_h \kappa_{-h})}}{\kappa_h \kappa_{-h}} = 1.
	\]
\end{proof}

\begin{proposition} \label{prop:normalize-lambda}
	The vertical balance condition is invariant under the transform $\Lambda_i
	\mapsto \tilde\Lambda_i = \Lambda_i + \lambda T_i$, $\lambda \in \C$.
\end{proposition}

\begin{proof}
	We use a tilde for all quantities associated to $\tilde\Lambda_i$.  Using
	Equations \eqref{eq:Hperiod} and \eqref{eq:DFhor}, the solutions
	of~\eqref{eq:system} satisfy $\tilde\xi - \xi = \lambda \cent{x}$.  Hence
	$\tilde K_h = K_h \exp(-\cent\ell_h \re\lambda)$ and $\tilde{F}^\text{ver}_b
	= \Fver_b \exp(-\cent\ell_b \re\lambda)$.
\end{proof}

We conclude this section with the following vertical
rigidity result:

\begin{proposition} \label{proposition:zero-is-rigid}
	Let $\cent\phi$ be a phase function such that $\cos(\cent\phi_h)$, $h \in
	\gH$, are all positive or all negative.  Then $\cent\phi$ is rigid. In
	particular, the zero phase function is always balanced and rigid.
\end{proposition}

\begin{proof}
	It suffices to prove the case where $\cos(\cent\phi_h)$ are all positive.
	Let $\dot\phi\in\cB$ such that $D\Fver_b(\cent\phi)\cdot\dot\phi=0$.  We can write
	$\dot\phi=\grad(f)$ with $f\in\cV$. Then
	\[
		0=d\Fver_b(\cent\phi)\cdot\dot\phi=\sum_{h\in m(b)}K_h\cos(\cent\phi_h)(f_{v(-h)}-f_{v(h)}).
	\]
	Since $K_h\cos(\cent\phi_h)>0$ for all $h\in\gH$, we conclude that $\dot\phi=0$
	by the same argument as in the proof of Proposition \ref{proposition:mdiv},
	considering the maximum of $f$.
\end{proof}

\section{Examples}\label{sec:example}

In this part, examples are sketched in the form of diagrams.  Edges are
decorated with arrows to illustrate a consistent orientation $\sigma$.  Unless
otherwise specified (e.g.~Example~\ref{ex:33t}), for each $h$ such that
$\sigma(h) = 1$, we label the phase difference $\cent\phi_h$ on the edge
$e(h)$.  The fundamental parallelogram of the torus spanned by $T_1$ and $T_2$
is illustrated by dotted lines.

For all the examples presented below, the computation of $K_h$ is either
trivial because of symmetry, or not necessary (e.g. Examples~\ref{ex:aH}).

\subsection{Genus three}\label{ssec:exg3}

The genus of a TPMS is at least three, hence the graph for our construction has
at least two faces.  We notice four families of balanced configurations with
two faces.  They are illustrated in Figure~\ref{fig:g3}.
Theorem~\ref{thm:TPMSg3} in the Appendix asserts that these are the only
balanced configurations whose graphs are orientable with two faces.  Hence they
are the only possible configurations that give rise to TPMSs of genus 3.

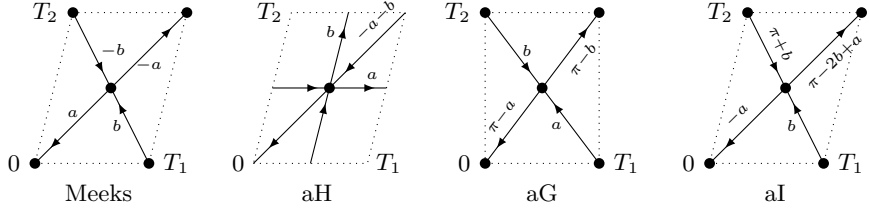
\begin{figure}[bht]
	\small
	\begin{tabular}[c]{cccc}
		\begin{tikzpicture}[scale=.5,baseline]
	\draw[dotted] (-2,-2) -- (1,-2) -- (2,2) -- (-1,2) -- cycle;
	\fill (0,0) node[circle,inner sep=0.05cm,fill=black]{};
	\fill (-2,-2) node[circle,inner sep=0.05cm,fill=black,label=left:$0$]{};
	\fill (1,-2) node[circle,inner sep=0.05cm,fill=black,label=right:$T_1$]{};
	\fill (2,2) node[circle,inner sep=0.05cm,fill=black]{};
	\fill (-1,2) node[circle,inner sep=0.05cm,fill=black,label=left:$T_2$]{};
	\draw[wing] (0,0) -- (-2,-2) node[midway,above] {$\scriptstyle a$};
	\draw[wing] (0,0) -- (2,2) node[midway,below] {$\scriptstyle -a$};
	\draw[wing] (1,-2) -- (0,0) node[midway,left] {$\scriptstyle b$};
	\draw[wing] (-1,2) -- (0,0) node[midway,right] {$\scriptstyle -b$};
\end{tikzpicture}&
		\begin{tikzpicture}[scale=0.5,baseline]
	\draw[dotted] (-2,-2) -- (1,-2) -- (2,2) -- (-1,2) -- cycle;
	\fill (0,0) node[circle,inner sep=0.05cm,fill=black]{};
	\node at (-2,-2) [left] {$0$};
	\node at (1,-2) [right] {$T_1$};
	\node at (-1,2) [left] {$T_2$};
	\draw[wing] (0,0) -- (1.5,0) node[near end,above] {$\scriptstyle a$};
	\draw[wing] (-1.5,0) -- (0,0);
	\draw[wing] (0,0) -- (0.5,2) node[near end,left] {$\scriptstyle b$};
	\draw[wing] (-0.5,-2) -- (0,0);
	\draw[wing] (0,0) -- (-2,-2);
	\draw[wing] (2,2) -- (0,0) node[near start,above,sloped] {$\scriptstyle -a-b$};
\end{tikzpicture}&
		\begin{tikzpicture}[scale=0.5,baseline]
	\draw[dotted] (-1.5,-2) -- (1.5,-2) -- (1.5,2) -- (-1.5,2) -- cycle;
	\fill (0,0) node[circle,inner sep=0.05cm,fill=black]{};
	\fill (-1.5,-2) node[circle,inner sep=0.05cm,fill=black,label=left:$0$]{};
	\fill (1.5,-2) node[circle,inner sep=0.05cm,fill=black,label=right:$T_1$]{};
	\fill (1.5,2) node[circle,inner sep=0.05cm,fill=black]{};
	\fill (-1.5,2) node[circle,inner sep=0.05cm,fill=black,label=left:$T_2$]{};
	\draw[wing] (0,0) -- (-1.5,-2) node[midway,above,sloped] {$\scriptstyle \pi-a$};
	\draw[wing] (1.5,-2) -- (0,0) node[midway,left] {$\scriptstyle a$};
	\draw[wing] (0,0) -- (1.5,2) node[midway,below,sloped] {$\scriptstyle \pi-b$};
	\draw[wing] (-1.5,2) -- (0,0) node[midway,right] {$\scriptstyle b$};
\end{tikzpicture}&
		\begin{tikzpicture}[scale=.5,baseline]
	\draw[dotted] (-2,-2) -- (1,-2) -- (2,2) -- (-1,2) -- cycle;
	\fill (0,0) node[circle,inner sep=0.05cm,fill=black]{};
	\fill (-2,-2) node[circle,inner sep=0.05cm,fill=black,label=left:$0$]{};
	\fill (1,-2) node[circle,inner sep=0.05cm,fill=black,label=right:$T_1$]{};
	\fill (2,2) node[circle,inner sep=0.05cm,fill=black]{};
	\fill (-1,2) node[circle,inner sep=0.05cm,fill=black,label=left:$T_2$]{};
	\draw[wing] (0,0) -- (-2,-2) node[midway,above,sloped] {$\scriptstyle -a$};
	\draw[wing] (0,0) -- (2,2) node[midway,below,sloped] {$\scriptstyle \pi-2b+a$};
	\draw[wing] (1,-2) -- (0,0) node[midway,left] {$\scriptstyle b$};
	\draw[wing] (-1,2) -- (0,0) node[midway,above,sloped] {$\scriptstyle \pi+b$};
\end{tikzpicture}\\
		Meeks & aH & aG & aI
	\end{tabular}
	\caption{
		The four balanced configurations that could give rise to TPMSs of genus
		three.\label{fig:g3}
	}
\end{figure}

\begin{example}[Meeks family] \label{ex:meeks}
	The first diagram illustrates a 4-parameter family (parameterized by $T_1$,
	$T_2$, $a$, and $b$).  It actually describes the Scherk limit of Meeks'
	family~\cite{meeks1990}.  To see this, note that the configurations are
	invariant under the translation $(T_1+T_2)/2$.  This implies an
	orientation-reversing translational symmetry in the corresponding TPMSs,
	which characterizes Meeks' surfaces.  See Figures~\ref{fig:gallery}(a--c) for
	examples in this family.

  Let us work out this small example explicitly. The graph is clearly balanced
  and rigid.  Let $v$ be the center vertex, and label $1$, $2$, $3$, $4$ the
  half-edges adjacent to $v$ in anti-clockwise order, starting with the
  half-edge marked $a$. Then the horizontal period condition~\eqref{eq:Hperiod}
  on any face cycle gives
  \[
  	\phi_1-\phi_2+\phi_3-\phi_4 = 0
  \]
  and the fundamental shifts are given by
  \[
  	\Psi_1= -\phi_1+\phi_2\quad\text{and}\quad
  	\Psi_2=-\phi_1+\phi_4.
  \]
  \begin{itemize}
    \item If $T_1$ and $T_2$ are orthogonal, all edges have the same length so
      $m(b(v))=\{1,2,3,4\}$ and
      \[
      	\Fver_{b(v)}(\phi)=K\left(\sin(\phi_1)+\sin(\phi_2)+\sin(\phi_3)+\sin(\phi_4)\right).
      \]
      A solution is $(\phi_1,\phi_2,\phi_3,\phi_4)=(a,-b,-a,b)$ with
      $a=\frac{-\Psi_1-\Psi_2}{2}$ and $b=\frac{\Psi_2-\Psi_1}{2}$.  Regarding
      rigidity, we have $\dot{\phi}\in\ker(\Pver,D\Fver(\phi))$ if and only if
      \[
      	\left\{\begin{array}{l}
        		\dot{\phi}_1-\dot{\phi}_2+\dot{\phi}_3-\dot{\phi}_4=0\\
        		-\dot{\phi}_1+\dot{\phi}_2=0\\
        		-\dot{\phi}_1+\dot{\phi}_4=0\\
        		\cos(a)(\dot{\phi}_1+\dot{\phi}_3)+\cos(b)(\dot{\phi}_2+\dot{\phi}_4)=0
      	\end{array}\right.
    	\]
      which gives $\dot{\phi}=0$ if $\cos(a)+\cos(b)\neq 0$, so the
      configuration is vertically rigid if $\cos(a)+\cos(b)\neq 0$.

    \item If $\arg(T_2/T_1)<\pi/2$, the edges $e(2)$ and $e(4)$ are shorter so
      $m(b(v))=\{2,4\}$ and
      \[
      	\Fver_{b(v)}(\phi)=K\left(\sin(\phi_2)+\sin(\phi_4)\right).
      \]
      The solution $(a,-b,-a,b)$ is rigid if $\cos(b)\neq 0$.

    \item If $\arg(T_2/T_1)>\pi/2$, the edges $e(1)$ and $e(3)$ are shorter, so
    	$m(b(v))=\{1,3\}$ and the solution $(a,-b,-a,b)$ is rigid if $\cos(a)\neq
    	0$.\qed
  \end{itemize}

\end{example}

\begin{example}[aH] \label{ex:aH}
	The second diagram is again a 4-parameter family.  We name it aH because,
	when $|T_1|=|T_2|$ and $a=b=0$, it gives Scherk limits of the oH
	family~\cite{chen2021a}; see Figure~\ref{fig:gallery}(d).  Another special
	case in this family is the Scherk limit of the rhombohedral deformation
	family rGL of the Gyroid~\cite{chen2019}, given by $T_2/T_1 =
	\exp(2\ii\pi/3)$ and $a=b=2\pi/3$; see Figure~\ref{fig:tGrGL} (left).  All
	configurations in the family are rigid, hence give rise to a new 5-parameter
	family of TPMSs, generalizing H and rGL.  \qed
\end{example}

Previously, we knew that both the Gyroid and the H surfaces can be continuously
deformed to Meeks surfaces~\cite{chen2019, chen2021a}.  The aH family implies a
deformation path between them that does not pass through the Meeks family.

\begin{corollary}
	The Gyroid can be continuously deformed to an H surface along a path in the
	space of TPMSs of genus 3 that stays outside the Meeks family.
\end{corollary}

\begin{example}[aG] \label{ex:aG}
	The third diagram is constrained to $\arg(T_2/T_1)=\pi/2$, hence a
	3-parameter family.  One of the fundamental shift must be $\pi$.  We name it
	aG as it includes the Scherk limit of the tetragonal deformation family tG of
	the Gyroid~\cite{chen2019}, given by $|T_1|=|T_2|$ and $a=b=\pi/2$; see
	Figure~\ref{fig:tGrGL} (right).  It intersects Meeks family when $a=b=0$.
	Moreover, when $a=0$ and $b=\pi$, we recognize alternative Scherk limits of
	the oH family~\cite{chen2021a}; see Figure~\ref{fig:gallery}(e).

	Unfortunately, configurations in this family are not vertically rigid, hence
	our construction is inconclusive for them.  To see this, note that adding a
	common constant to $a$ and $b$ does not change the fundamental shifts.  For
	the corresponding TPMSs, this seems to suggest that one can vertically slide
	one Scherk tower with respect to the other without changing the lattice.
	Numerical experiments suggest that these configurations give indeed rise to
	TPMSs, but a vertical sliding between the towers must be accompanied by a
	very slight deformation of the horizontal lattice, which becomes undetectable
	in the Scherk limit.  \qed
\end{example}

\begin{example}[aI] \label{ex:aI}
	The fourth diagram is again a 4-parameter family, but constrained to
	$\arg(T_2/T_1) \ne \pi/2$.  The fundamental shifts satisfy $\Psi_2-\Psi_1 =
	\pi$.  When $\Psi_1=\pi$, it tends to aG configurations as $\arg(T_2/T_1) \to
	\pi/2$.  Otherwise, the dependence of the vertical balancing condition on the
	horizontal lattice is not continuous; see Remark~\ref{rmk:continuous}.

	Configurations in this family are not vertically rigid: adding a common
	constant to $a$ and $b$ does not change the fundamental shifts.  Our
	construction is inconclusive and, apart from the aG limit, we are not aware
	of any known TPMS that admits this kind of Scherk limit.  But by
	Theorem~\ref{thm:rami}, the trivial phase functions given by $a,b=0$ or $\pi$
	should give rise to a family of TPMSs.  For a lack of a better name, we call
	this family aI (following the pattern of aG, aH.) \qed
\end{example}

\subsection{Technical examples with triangular lattice}

\begin{figure}[bht]
	\begin{tikzpicture}[rotate=90,scale=0.5,baseline]
	\draw[dotted] (-2,-6) -- (2,-6) -- (2,6) -- (-2,6) -- cycle;
	\fill (0,0) node[circle,inner sep=0.05cm,fill=black]{};
	\fill (-2,0) node[circle,inner sep=0.05cm,fill=black]{};
	\fill (2,0) node[circle,inner sep=0.05cm,fill=black]{};
	\fill (-1,-6) node[circle,inner sep=0.05cm,fill=black]{};
	\fill (1,-6) node[circle,inner sep=0.05cm,fill=black]{};
	\fill (1,6) node[circle,inner sep=0.05cm,fill=black]{};
	\fill (-1,6) node[circle,inner sep=0.05cm,fill=black]{};
	\draw[wing] (0,0) -- (-1,-6);
	\draw[wing] (1,-6) -- (0,0);
	\draw[wing] (-2,-6) -- (-1,-6);
	\draw[wing] (-1,-6) -- (1,-6);
	\draw[wing] (1,-6) -- (2,-6);
	\draw[wing] (-1,6) -- (-2,0);
	\draw[wing] (-1,-6) -- (-2,0);
	\draw[wing] (-2,0) -- (0,0);
	\draw[wing] (0,0) -- (2,0);
	\draw[wing] (2,0) -- (1,6);
	\draw[wing] (2,0) -- (1,-6);
	\draw[wing] (0,0) -- (-1,6);
	\draw[wing] (1,6) -- (0,0);
	\draw[wing] (-2,6) -- (-1,6);
	\draw[wing] (-1,6) -- (1,6);
	\draw[wing] (1,6) -- (2,6);
\end{tikzpicture}
	\caption{
		For vertical balancing on this graph, it is not enough to consider only
		vertex cuts. \label{fig:necessarycuts}
	}
\end{figure}
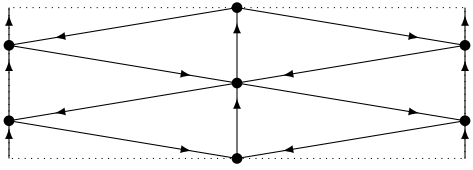

\begin{example}\label{ex:necessarycuts}
	The purpose of this example is to demonstrate the necessity to define
	vertical balancing on the whole cut space.

	Figure~\ref{fig:necessarycuts} illustrates a graph with four vertices.  If
	only vertex cuts are considered, solving the vertical balancing equation
	would only determine phase differences on the four vertical edges.  This
	gives a false illusion that one may vertically slide the saddle towers with
	respect to each other.  In fact, to determine phase differences on the
	remaining edges, we must acknowledge that they form a cut and must also be
	balanced.  \qed
\end{example}

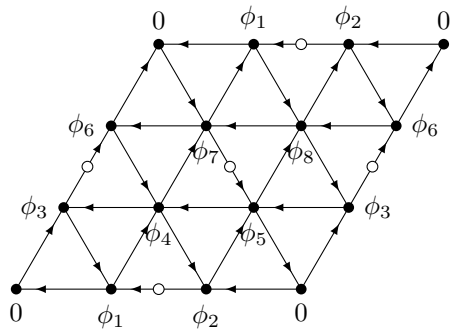
\begin{figure}[bht]
	\begin{tikzpicture}[scale=1.25,baseline]
	\fill (0,0) node[circle,inner sep=0.05cm,fill=black,label={below:$0$}]{};
	\fill (1,0) node[circle,inner sep=0.05cm,fill=black,label={below:$\phi_1$}]{};
	\fill (2,0) node[circle,inner sep=0.05cm,fill=black,label={below:$\phi_2$}]{};
	\fill (3,0) node[circle,inner sep=0.05cm,fill=black,label={below:$0$}]{};
	\fill (0,0)++(60:1) node[circle,inner sep=0.05cm,fill=black,label={left:$\phi_3$}]{};
	\fill (1,0)++(60:1) node[circle,inner sep=0.05cm,fill=black,label={below:$\phi_4$}]{};
	\fill (2,0)++(60:1) node[circle,inner sep=0.05cm,fill=black,label={below:$\phi_5$}]{};
	\fill (3,0)++(60:1) node[circle,inner sep=0.05cm,fill=black,label={right:$\phi_3$}]{};
	\fill (0,0)++(60:2) node[circle,inner sep=0.05cm,fill=black,label={left:$\phi_6$}]{};
	\fill (1,0)++(60:2) node[circle,inner sep=0.05cm,fill=black,label={below:$\phi_7$}]{};
	\fill (2,0)++(60:2) node[circle,inner sep=0.05cm,fill=black,label={below:$\phi_8$}]{};
	\fill (3,0)++(60:2) node[circle,inner sep=0.05cm,fill=black,label={right:$\phi_6$}]{};
	\fill (0,0)++(60:3) node[circle,inner sep=0.05cm,fill=black,label={above:$0$}]{};
	\fill (1,0)++(60:3) node[circle,inner sep=0.05cm,fill=black,label={above:$\phi_1$}]{};
	\fill (2,0)++(60:3) node[circle,inner sep=0.05cm,fill=black,label={above:$\phi_2$}]{};
	\fill (3,0)++(60:3) node[circle,inner sep=0.05cm,fill=black,label={above:$0$}]{};
	\draw[wing] (0,0) -- +(60:1);
	\draw[wing] (1,0) -- +(60:1);
	\draw[wing] (2,0) -- +(60:1);
	\draw[wing] (3,0) -- +(60:1);
	\draw[wing] (0,0)++(60:1) -- +(60:1);
	\draw[wing] (1,0)++(60:1) -- +(60:1);
	\draw[wing] (2,0)++(60:1) -- +(60:1);
	\draw[wing] (3,0)++(60:1) -- +(60:1);
	\draw[wing] (0,0)++(60:2) -- +(60:1);
	\draw[wing] (1,0)++(60:2) -- +(60:1);
	\draw[wing] (2,0)++(60:2) -- +(60:1);
	\draw[wing] (3,0)++(60:2) -- +(60:1);
	\draw[wing] (1,0) -- +(0:-1);
	\draw[wing] (2,0) -- +(0:-1);
	\draw[wing] (3,0) -- +(0:-1);
	\draw[wing] (1,0)++(60:1) -- +(180:1);
	\draw[wing] (2,0)++(60:1) -- +(180:1);
	\draw[wing] (3,0)++(60:1) -- +(180:1);
	\draw[wing] (1,0)++(60:2) -- +(180:1);
	\draw[wing] (2,0)++(60:2) -- +(180:1);
	\draw[wing] (3,0)++(60:2) -- +(180:1);
	\draw[wing] (1,0)++(60:3) -- +(180:1);
	\draw[wing] (2,0)++(60:3) -- +(180:1);
	\draw[wing] (3,0)++(60:3) -- +(180:1);
	\draw[wing] (0,0)++(60:1) -- +(-60:1);
	\draw[wing] (1,0)++(60:1) -- +(-60:1);
	\draw[wing] (2,0)++(60:1) -- +(-60:1);
	\draw[wing] (0,0)++(60:2) -- +(-60:1);
	\draw[wing] (1,0)++(60:2) -- +(-60:1);
	\draw[wing] (2,0)++(60:2) -- +(-60:1);
	\draw[wing] (0,0)++(60:3) -- +(-60:1);
	\draw[wing] (1,0)++(60:3) -- +(-60:1);
	\draw[wing] (2,0)++(60:3) -- +(-60:1);
	\draw (0,0)++(60:1.5) node[circle, inner sep=0.05cm, draw=black, fill=white]{};
	\draw (1.5,0)++(60:1.5) node[circle, inner sep=0.05cm, draw=black, fill=white]{};
	\draw (3,0)++(60:1.5) node[circle, inner sep=0.05cm, draw=black, fill=white]{};
	\draw (1.5,0) node[circle, inner sep=0.05cm, draw=black, fill=white]{};
	\draw (1.5,0)++(60:3) node[circle, inner sep=0.05cm, draw=black, fill=white]{};
\end{tikzpicture}

	\caption{
		A configuration with the graph of the triangular lattice.  Phase functions
		are labeled on the vertices instead of the edges.  We look for phase
		functions symmetric in the points marked with the empty circles.
		\label{fig:33t}
	}
\end{figure}

\begin{example}\label{ex:33t}
	The purpose of this example is to demonstrate the diversity of balanced phase
	functions.

	The graph in Figure~\ref{fig:33t}, which is a $3\times 3$ block of the
	triangular lattice, is trivially balanced and rigid by Theorem
	\ref{theorem:triangles}.  Unlike other diagrams in this section, the phases
	of the saddle towers are labeled on the vertices up to the addition of a
	common constant.  This is possible because the fundamental shifts are 0.

	On this small graph, we look for phase functions symmetric in the points
	marked with empty circles in the figure.  That is,
	\[
		\phi_1 = -\phi_2,\quad
		\phi_3 = -\phi_6,\quad
		\phi_4 = -\phi_8,\quad
		\phi_5 = -\phi_7.
	\]
	Under these strict restrictions, we still find two non-trivial, balanced, and
	rigid phase functions, namely
	\begin{itemize}
		\item $-\phi_1 = \phi_2 = 2\pi/3$,
			$-\phi_3 = \phi_5 = \phi_6 = -\phi_7 = \pi/3$, and
			$\phi_4 = -\phi_8 = \pi$;
		\item $\phi_1 = -\phi_2 = -\phi_3 = \phi_5 = \phi_6 = -\phi_7 = 2\arctan\sqrt{5/7}$
			and $\phi_4 = -\phi_8 = \pi$.
	\end{itemize}
	For readers who are interested in double checking: the rigidities are
	confirmed numerically by computing the determinant of an $8 \times 8$
	Jacobian matrix, which is the derivative of vertical forces on 8 vertices
	with respect to their phases.  If we assume $K_h = 1$ on all half-edges (up
	to scaling of local coordinates), the Jacobian determinant is $-3/4$ for the
	first phase function, and $-315/4$ for the second.

	We certainly did not find all non-trivial, balanced, and rigid phase
	functions on this graph.  It can be imagined that, if we take a larger block
	of triangular lattice or relax the inversion symmetry, there will be more
	balanced and rigid phase functions. \qed
\end{example}

\subsection{Generalizing some known examples}

We generalize here some interesting known examples to demonstrate the power of
our construction.  Many known families of TPMSs admit saddle tower limits.  In
addition to those discussed below, examples also include Schoen's so-called
RII, RIII, I-6, I-8, I-9\footnotemark, and GW surfaces~\cite{schoen1970,
brakke}, and many more constructed in~\cite{fischer1987,karcher1989}.  Brakke's
webpage~\cite{brakke} is a great source of examples.  We certainly do not plan
to discuss all of them.  \footnotetext{These names are coined by Brakke for
they are the 6th, 8th, and 9th surface on Page I of a note by Schoen.}

\medskip

If a graph contains a pair of parallel edges, they must be adjacent, and their
orientations and phase differences must all be opposite.  In fact, given any
graph $\gG = (\gH, \involute, \rotate, \varrho)$ with only simple edges, we may
construct a \emph{doubling graph} $\bar \gG = (\bar\gH, \bar\involute,
\bar\rotate, \bar\varrho)$ with only parallel edges as follows (we use a bar
for all objects associated to $\bar\gG$).
\begin{itemize}
	\item For each $h \in \gH$, we have two half-edges $\bar h_+, \bar h_- \in
		\bar\gH$;

	\item We define $\bar\involute$ by $\bar\involute(h_+) = (\involute(h))_-$
		and $\bar\involute(h_-) = (\involute(h))_+$;

	\item We define $\bar\rotate$ by $\bar\rotate(h_+) = h_-$ and
		$\bar\rotate(h_-) = (\rotate(h))_+$;

	\item We define $\bar\varrho$ by $\bar\varrho(\bar v(\bar h_\pm)) =
		\varrho(v(h))$ and $\bar\varrho(\bar e(\bar h_\pm)) = \varrho(e(h))$.
\end{itemize}

Clearly, the doubling graph $\bar\gG$ is always orientable, even if the
original graph $\gG$ is not.  This observation expands the power of our
constructions: Even if the graph is not orientable, as long as a graph is
balanced and rigid, its doubling would be an orientable, balanced, and rigid
graph.

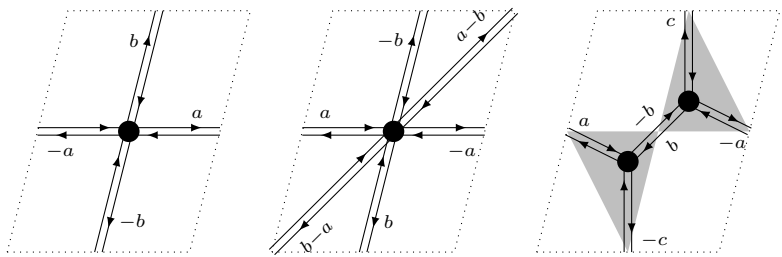
\begin{figure}[!bht]
	\begin{tikzpicture}[scale=.8,baseline] 
	\draw[dotted] (-2,-2) -- (1,-2) -- (2,2) -- (-1,2) -- cycle;
	\fill (0,0) node[circle,inner sep=0.1cm,fill=black]{};
	\draw[wing] (0,0)++(0,.062) -- ++(1.5,0) node[near end, above] {$\scriptstyle a$};
	\draw[wing] (1.5,0)++(0,-.062) -- ++(-1.5,0);
	\draw[wing] (0,0)++(-.062,0) -- ++(0.5,2) node[near end, left] {$\scriptstyle b$};
	\draw[wing] (0.5,2)++(.062,0) -- ++(-0.5,-2);
	\draw[wing] (0,0)++(0,-.062) -- ++(-1.5,0) node[near end, below] {$\scriptstyle -a$};
	\draw[wing] (-1.5,0)++(0,.062) -- ++(1.5,0);
	\draw[wing] (0,0)++(.062,0) -- ++(-0.5,-2) node[near end, right] {$\scriptstyle -b$};
	\draw[wing] (-0.5,-2)++(-.062,0) -- ++(0.5,2);
\end{tikzpicture} 
	\begin{tikzpicture}[scale=.8,baseline] 
	\draw[dotted] (-2,-2) -- (1,-2) -- (2,2) -- (-1,2) -- cycle;
	\fill (0,0) node[circle,inner sep=0.1cm,fill=black]{};
	\draw[wing] (0,0)++(0,.062) -- ++(1.5,0);
	\draw[wing] (1.5,0)++(0,-.062) -- ++(-1.5,0) node[near start, below] {$\scriptstyle -a$};
	\draw[wing] (0,0)++(-.062,0) -- ++(0.5,2) node[near end, left] {$\scriptstyle -b$};
	\draw[wing] (0.5,2)++(.062,0) -- ++(-0.5,-2);
	\draw[wing] (0,0)++(0,-.062) -- ++(-1.5,0);
	\draw[wing] (-1.5,0)++(0,.062) -- ++(1.5,0) node[near start, above] {$\scriptstyle a$};
	\draw[wing] (0,0)++(.062,0) -- ++(-0.5,-2) node[near end, right] {$\scriptstyle b$};
	\draw[wing] (-0.5,-2)++(-.062,0) -- ++(0.5,2);
	\draw[wing] (0,0)++(135:.062) -- ++(2,2) node[near end, above,sloped] {$\scriptstyle a-b$};
	\draw[wing] (2,2)++(-45:.062) -- ++(-2,-2);
	\draw[wing] (0,0)++(-45:.062) -- ++(-2,-2) node[near end, below, sloped] {$\scriptstyle b-a$};
	\draw[wing] (-2,-2)++(135:.062) -- ++(2,2);
\end{tikzpicture}
	\begin{tikzpicture}[scale=.8,baseline] 
	\draw[dotted] (-2,-2) -- (1,-2) -- (2,2) -- (-1,2) -- cycle;
	\fill[lightgray] (0,0) -- (-1.5,0) -- (-0.5,-2) -- cycle;
	\fill[lightgray] (0,0) -- (1.5,0) -- (0.5,2) -- cycle;
	\fill (-0.5,-0.5) node[circle,inner sep=0.1cm,fill=black]{};
	\fill (0.5,0.5) node[circle,inner sep=0.1cm,fill=black]{};
	\draw[wing] (-0.5,-0.5)++(135:.062) -- ++(1,1) node[midway,above,sloped] {$\scriptstyle -b$};
	\draw[wing] (0.5,0.5)++(-45:.062) -- ++(-1,-1) node[midway,below,sloped] {$\scriptstyle b$};
	\draw[wing] (-0.5,-0.5)++(0:.062) -- ++(0,-1.5) node[very near end,right] {$\scriptstyle -c$};
	\draw[wing] (-0.5,-2)++(180:.062) -- ++(0,1.5);
	\draw[wing] (0.5,0.5)++(180:.062) -- ++(0,1.5) node[very near end,left] {$\scriptstyle c$};
	\draw[wing] (0.5,2)++(0:.062) -- ++(0,-1.5);
	\draw[wing] (-0.5,-0.5)++(-120:.062) -- ++(-1,0.5);
	\draw[wing] (-1.5,0)++(60:.062) -- ++(1,-0.5) node[near start,above] {$\scriptstyle a$};
	\draw[wing] (0.5,0.5)++(60:.062) -- ++(1,-0.5);
	\draw[wing] (1.5,0)++(-120:.062) -- ++(-1,0.5) node[near start,below] {$\scriptstyle -a$};
\end{tikzpicture}
	\caption{
		Configurations doubling the graphs of parallelogram, triangle and hexagonal
		tilings.
	} \label{fig:double}
\end{figure}

\begin{example}
	In Figure~\ref{fig:double}, we illustrate three configurations that double
	the graphs of tilings of the Euclidean plane.  Each diagram describes a
	4-parameter family of balanced configurations.  They generalize,
	respectively, the Scherk limits of Schoen's S'--S'', T'--R', and H''--R
	families~\cite{schoen1970}; see Figures~\ref{fig:gallery}(g--i).  The first
	two configurations in the figure are trivially balanced and rigid, hence give
	rise to 5-parameter families of TPMSs.

	As for the last one that doubles hexagonal tilings, the graph is balanced if
	and only if the vertices are at the unique Fermat--Torricelli points of the
	gray triangles spanned by the 2-division points of the torus.  This is well
	defined only when these triangles do not have an angle $\ge 2\pi/3$.  Let
	$l_a$, $l_b$ and $l_c$ be the lengths of the edges labeled by $\pm a$, $\pm
	b$ and $\pm c$ in the figure, respectively.  We may assume that $l_a \le l_b
	\le l_c$.  Then the balance condition for the phase function is
	\[
		\left.
			\begin{array}{r}
				e^{\ii a}\\
				e^{\ii a} + e^{\ii b}\\
				e^{\ii a} + e^{\ii b} + e^{\ii c}
			\end{array}
		\right\}
		\text{being real, if}
		\left\{
			\begin{array}{l}
				l_a < l_b \le l_c;\\
				l_a = l_b < l_c;\\
				l_a = l_b = l_c.
			\end{array}
		\right.
	\]
	The phase function is rigid if the left hand side is non-zero.

	We see again that balancing conditions do not depend continuously on the
	horizontal lattice.  But the dependence on $\Lambda_1$ and $\Lambda_2$ is
	continuous, so the rigid configurations still give rise to 5-parameter
	families of TPMSs.  See Remark~\ref{rmk:continuous}.
	\qed
\end{example}

We end the section with the saddle tower limit of a known TPMS family, which is
however not horizontally rigid.

\begin{figure}[!htb]
	\begin{tikzpicture}[scale=.8,baseline]
	\draw[dotted] (-2,-2) -- (1,-2) -- (2,2) -- (-1,2) -- cycle;
	\fill (0,0) node[circle,inner sep=0.05cm,fill=black]{};
	\fill (1.5,0) node[circle,inner sep=0.05cm,fill=black]{};
	\fill (-1.5,0) node[circle,inner sep=0.05cm,fill=black]{};
	\fill (0.5,2) node[circle,inner sep=0.05cm,fill=black]{};
	\fill (-0.5,-2) node[circle,inner sep=0.05cm,fill=black]{};
	\node at (-2,-2) [left] {$0$};
	\node at (1,-2) [right] {$T_1$};
	\node at (-1,2) [left] {$T_2$};
	\draw[wing] (0,0) -- (1.5,0) node[midway,below] {$\scriptstyle a$};
	\draw[wing] (1.5,0) -- (0.5,2) node[midway,above,sloped] {$\scriptstyle b-a$};
	\draw[wing] (0.5,2) -- (0,0) node[midway,left] {$\scriptstyle -b$};
	\draw[wing] (0,0) -- (-1.5,0) node[midway,above] {$\scriptstyle -a$};
	\draw[wing] (-1.5,0) -- (-0.5,-2) node[midway,below,sloped] {$\scriptstyle a-b$};
	\draw[wing] (-0.5,-2) -- (0,0) node[midway,right] {$\scriptstyle b$};
\end{tikzpicture}
	\begin{tikzpicture}[scale=1.25,baseline]
	\draw[dotted] (-60:1) -- ++(180:2) -- ++(60:2) -- ++(0:2) -- cycle;
	\fill (0,0) node[circle,inner sep=0.05cm,fill=black]{};
	\fill (0:1) node[circle,inner sep=0.05cm,fill=black]{};
	\fill (60:1) node[circle,inner sep=0.05cm,fill=black]{};
	\fill (180:1) node[circle,inner sep=0.05cm,fill=black]{};
	\fill (240:1) node[circle,inner sep=0.05cm,fill=black]{};
	\draw[wing] (0,0) -- (0:1) node[midway,below] {$\scriptstyle 2\pi/3$};
	\draw[wing] (0:1) -- (60:1) node[midway,above,sloped] {$\scriptstyle 2\pi/3$};
	\draw[wing] (60:1) -- (0,0) node[midway,above,sloped] {$\scriptstyle 2\pi/3$};
	\draw[wing] (0,0) -- (180:1) node[midway,above] {$\scriptstyle 2\pi/3$};
	\draw[wing] (180:1) -- (240:1) node[midway,below,sloped] {$\scriptstyle 2\pi/3$};
	\draw[wing] (240:1) -- (0,0) node[midway,below,sloped] {$\scriptstyle 2\pi/3$};
\end{tikzpicture}
	\caption{
		Left: Configurations generalizing the Scherk limit of Schoen's H'--T
		family. Right: Scherk limit of the QTZ--QZD family, a chiral generalization
		of H'--T. \label{fig:HT}
	}
\end{figure}
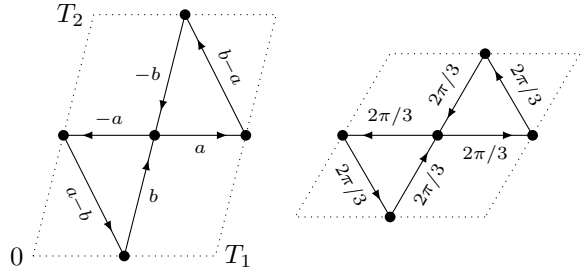

\begin{example}
	On the left of Figure~\ref{fig:HT} is a configuration that generalizes the
	Scherk limit of Schoen's H'--T family~\cite{schoen1970} (with $T_2/T_1 =
	\exp(\ii\pi/3)$ and $a=b=0$; see Figure~\ref{fig:gallery}(f)).  The graph is
	obviously balanced but, unfortunately, not rigid.  To see this, notice that
	the graph is represented as the union of three lines.  Any of the lines can
	be moved parallelly yet the graph remains balanced.  As a consequence, our
	construction does not work directly on this graph.

	But we may impose an inversion symmetry in the vertices.  Under this imposed
	symmetry, the phase function must have the form as shown in the figure.  The
	configurations form a 4-parameter family.  They are trivially balanced and
	rigid modulo symmetries.  Our construction works with little modification,
	and gives rise to a 5-parameter family of TPMSs of genus 4 generalizing the
	H'--T surfaces.  We do not plan to write down the details.  \qed
\end{example}

\begin{remark}
	Recently, a chiral variation of H'--T was discovered and named
	QTZ--QZD~\cite{markande2018}.  It arises from the same graph as H'--T, but
	has a non-trivial phase function; see the right side of Figure~\ref{fig:HT}.
\end{remark}

\section{Construction}\label{sec:construction}

In this section, we prove Theorem \ref{thm:main}.  The given phase function is
denoted $\cent\phi$.  All parameters will vary in a neighborhood of a central
value denoted with a superscript $\circ$, which depends on the given
configuration $(\gG, \cent\phi)$.  The Implicit Function Theorem will be
applied at $\varepsilon=0$ and the central value of all parameters.

Without loss of generality, we assume as in Section \ref{sec:rigid-tower} that
$\cos(\cent{\theta}_h)\neq 0$ for all $h\in\gH$ and that the zeros of $\Phi_1$
are simple for all saddle towers involved in the construction. If this is not
the case, we could always apply a horizontal rotation to the configuration.

\subsection{Opening nodes} \label{ssec:opennode}

For each vertex $v \in \gV$, we want to place a saddle tower $\sS_v$ determined
by the angles $\cent{\theta}_h$, $h \in v$, and glue them along the wings.

Recall that a saddle tower is conformally a sphere with punctures corresponding
to the wings.  Hence our initial surface at $\varepsilon=0$ is a singular
Riemann surface consisting of $|\gV|$ spheres identified at their punctures.
More specifically: to each vertex $v \in \gV$ of the oriented graph $\gG$, we
associate a Riemann sphere $\hat\C_v$.  To each half-edge $h \in v$, we
associate a complex number $\cent p_h \in \hat\C_v$, so that $\hat\C_v$
punctured at $\cent p_h$, $h \in v$, provides a conformal model for $\sS_v$.
Then we identify $\cent p_h$ and $\cent p_{-h}$.  The resulting singular
Riemann surface with nodes is denoted $\Sigma_0$.

\medskip

As $\varepsilon$ increases, we want to desingularize the nodes into necks.
This is done as follows.  For each $h \in \gH$, $p_h$ is a complex parameter in
a neighborhood of $\cent p_h$.  Consider a local coordinate $w_h$ in a
neighborhood of $p_h \in \hat\C_{v(h)}$ such that $w_h(p_h) = 0$.  The local
coordinate $w_h$ depends holomorphically on $p_h$, but the dependence is
omitted for simplicity.  Moreover, we assume that the local coordinate $\cent
w_h$ associated to $\cent p_h$ is adapted. We denote $p=(p_h)_{h\in\gH}$.

Since the graph is finite, it is possible to fix a small number $\delta > 0$
independent of $v$ such that the disks $|\cent w_h| < 2\delta$ for $h\in v$ are
disjoint in each Riemann sphere $\hat\C_v$.  Then for $p$ close enough to
$\cent{p}$, the disks $|w_h|<\delta$ for $h\in v$ are disjoint.

Consider a symmetric complex parameter $t = (t_h)_{h \in \gH}$ in the
neighborhood of $0$ with $|t_h| < \delta^2$.  For every $h \in \gH$, if
$t_h\neq 0$, we remove the disk
\[
	|w_h| < |t_h|/\delta,
\]
and identify the annuli
\[
	|t_h| / \delta \le |w_h| \le \delta
	\quad \text{and} \quad
	|t_{-h}| / \delta \le |w_{-h}| \le \delta
\]
by
\[
	w_h w_{-h} = t_h.
\]
When $t_h = 0$, $p_h$ and $p_{-h}$ are simply identified to form a node.  This
produces a Riemann surface, possibly with nodes, denoted by $\Sigma_t$.  Note
that $\Sigma_t$ also depends on the parameter $p$, but the dependence is not
written for simplicity.  When $t_h \ne 0$ for all $h \in \gH$, $\Sigma_t$ is a
regular Riemann surface of genus $g=|\gF|+1$ which provides the conformal model
for our construction.

We consider the following fixed domains in all $\Sigma_t$:
\[
	U_{v,\delta}=\{ z \in \hat\C_v \colon \forall h\in v,\; |\cent w_h| > \delta/2 \}
	\quad\text{and}\quad
	U_{\delta}=\bigsqcup_{v\in\gV} U_{v,\delta}
\]
For $v\in\gV$, we denote $\infty_v$ the point at infinity in the Riemann sphere
$\hat\C_v$ and fix an origin $O_v \in U_{v,\delta}\setminus\{\infty_v\}$ as the
starting point of the integration defining the Weierstrass parameterization of
the saddle tower $\sS_v$.  

\subsection{Regular 1-forms}

Let $\Sigma_t$ be a family of Riemann surfaces defined by opening nodes as
above.  A \emph{regular 1-form} $\omega$ on $\Sigma_t$ is a differential 1-form
that is holomorphic away from the nodes and, whenever two points $p$ and $q$
are identified to form a node, has simple poles of opposite residues at $p$ and
$q$.  Regular $1$-forms extend the notion of holomorphic 1-forms to noded
Riemann surfaces.  By~\cite[Proposition 4.1]{masur1976}, there is a basis
$\omega_{1,t},\cdots,\omega_{g,t}$ for the space of regular 1-forms on
$\Sigma_t$ which ``depends holomorphically'' on $t$ in a neighborhood of $0$.
More formally, in our case, this means that the restriction of $\omega_{j,t}$
to $U_\delta$ depends holomorphically on $z \in U_\delta$ and $t$.

One can also consider regular 1-forms with simple poles away from the nodes: by
\cite[Proposition 4.2]{masur1976}, if $p$, $q$ are two points on $\Sigma_0$
minus the nodes, there exists a unique 1-form $\omega_{p,q,t}$ on $\Sigma_t$
which has simple poles at $p$ and $q$ with residues $1$ and $-1$, is otherwise
regular in the sense above, and has suitably normalized periods (a normalized
differential of the third kind).  Moreover, $\omega_{p,q,t}$ depends
holomorphically on $t$ in a neighborhood of $0$.

More specifically, we shall use the following result.  For each half-edge
$h\in\gH$, let $A_h$ denote a small anticlockwise circle in $U_{v(h),\delta}$
around $p_h$; it is then homologous in $\Sigma_t$ to a clockwise circle in
$U_{v(-h),\delta}$ around $p_{-h}$.
\begin{proposition} \label{proposition:regular-1forms}
 	Given an antisymmetric function $(\alpha_h)_{h\in\gH}$, there exists a unique
 	regular 1-form $\omega_t$ on $\Sigma_t$, possibly with simple poles at
 	$\infty_v$, $v\in\gV$, such that
	\[
		\forall h\in\gH,\quad\int_{A_h}\omega_t=2\pi\ii\, \alpha_h.
	\]
	Moreover, the restriction of $\omega_t$ to $U_{\delta}$ depends
	holomorphically on $t$ in a neighborhood of $0$.
\end{proposition}
\begin{proof}
	this follows from \cite[Proposition 4.2]{masur1976}. See also \cite[Theorem
	8.2]{traizet2013} for a constructive proof.
\end{proof}
Note that by the Residue Theorem in $\hat\C_v$,
\[
	\Res{\omega_t}{\infty_v}=-\sum_{h\in v}\alpha_h.
\]

We will also need the following results.

\begin{lemma}[{\cite[Lemma~3]{traizet2008}}] \label{lem:derivative}
	The derivative $\partial\omega_t/\partial t_h$ at $t=0$, restricted to
	$U_\delta$, coincides with a meromorphic 1-form on $\Sigma_0$ with double
	poles at the nodes, holomorphic elsewhere, and vanishing A-periods.  In term
	of the local complex coordinates $w_h$ used to open nodes, the principal part
	at $p_h$ is
	\[
		-\frac{d\,w_h}{w_h^2}\Res{\frac{\omega_0}{w_{-h}}}{p_{-h}}.
	\]
\end{lemma}

For every half-edge $h$ and $t_h\neq 0$, let $B_h$ be the concatenation of
\begin{enumerate}
	\item a path in $U_{v(h),\delta}$ from $O_{v(h)}$ to $w_h = \delta$,
	\item the path parameterized by $w_h = \delta^{1-2s}\,t_h^s$ for $s \in
		[0,1]$, from $w_h = \delta$ to $w_h = t_h/\delta$, which is identified with
		$w_{-h}=\delta$, and
	\item a path in $U_{v(-h),\delta}$ from $w_{-h}=\delta$ to $O_{v(-h)}$.
\end{enumerate}
\begin{lemma}[{\cite[Lemma 1]{traizet2002}}] \label{lem:Bintegral}
	The difference
	\begin{equation}\label{eq:the-difference}
		\Big(\int_{B_h} \omega_t\Big) - \alpha_h\log t_h 
	\end{equation}
	extends holomorphically to
	$t_h = 0$.
	Moreover, its value at $t=0$ is equal to
	\[
		\lim_{z\to p_h}\bigg[\Big(\int_{O_{v(h)}}^{z}\omega_0\Big)-\alpha_h\log w_h(z)\bigg]
		-\lim_{z\to p_{-h}}\bigg[\Big(\int_{O_{v(-h)}}^{z}\omega_0\Big)-\alpha_{-h}\log w_{-h}(z)\bigg].
	\]
\end{lemma}

Lemma \ref{lem:Bintegral} was essentially proved in~\cite[Lemma
1]{traizet2002}.  As the lemma has been and will be used in similar
constructions, we consider it a good time to refurbish the proof in Appendix
\ref{appendix:neck}.

\subsection{Weierstrass data}

We construct a conformal minimal immersion using the Weierstrass
parameterization in the form
\[
	z\mapsto	\re \int^z (\Phi_1, \Phi_2, \Phi_3) ,
\]
where $\Phi_i$ are meromorphic 1-forms on $\Sigma_t$ satisfying the
conformality equation
\begin{equation}\label{eq:conformal}
	Q := \Phi_1^2 + \Phi_2^2 + \Phi_3^2 = 0.
\end{equation}
Observe that $Q$ is a meromorphic quadratic differential on $\Sigma_t$.

\subsubsection{A-periods}

We need to solve the following A-period problem
\[
	\re \int_{A_h} (\Phi_1, \Phi_2, \Phi_3) = (0,0, 2\pi\sigma_h), \quad \forall h \in \gH.
\]
We define $\Phi_1$, $\Phi_2$, and $\Phi_3$, using
Proposition~\ref{proposition:regular-1forms}, as the unique regular 1-forms on
$\Sigma_t$ with (at most) simple poles at $\infty_v$ for $v\in\gV$ and the
A-periods
\[
	\int_{A_h} (\Phi_1, \Phi_2, \Phi_3) =
	2 \pi \ii (\alpha_h, \beta_h, \gamma_h - \ii \sigma_h), \quad\forall h\in\gH,
\]
where $(\alpha,\beta,\gamma) \in \cA^3$ are antisymmetric parameters.
This way, the A-period problems are solved by definition.  We choose the
following central value for the parameters:
\[
	\cent\alpha_h = -\cos(\cent\theta_h),\qquad
	\cent\beta_h = -\sin(\cent\theta_h),\qquad
	\cent\gamma_h = 0.
\]
Then at $\varepsilon=0$ and the central value of all parameters, we have in $\hat\C_v$
\[
	\cent\Phi_1=\sum_{h\in v}\frac{-\cos(\cent\theta_h)}{z-\cent p_h}dz,\quad
	\cent\Phi_2=\sum_{h\in v}\frac{-\sin(\cent\theta_h)}{z-\cent p_h}dz\quad\text{and}\quad
	\cent\Phi_3=\sum_{h\in v}\frac{-\ii\sigma_h}{z-\cent p_h}dz.
\]
In other words, $(\cent\Phi_1,\cent\Phi_2,\cent\Phi_3)$ is precisely the
Weierstrass data of the saddle tower $\sS_v$ as we want.

Note that $\sigma \in \ker(\div)$, that is
\[
	\sum_{h \in v} \sigma_h = 0,
\]
so we have by Residue Theorem in $\hat\C_v$
\[
	\Res{\Phi}{\infty_v} =
	-\sum_{h\in v}(\alpha_h,\beta_h,\gamma_h)=
	-\div_{b(v)}(\alpha, \beta, \gamma).
\]

We want $\infty_v$ to be regular points, so we need to solve
\begin{equation}\label{eq:balance}
	\div_{b(v)}(\alpha, \beta, \gamma) =
	\sum_{h \in v} (\alpha_h, \beta_h, \gamma_h) = 0,
	\qquad \text{for all } v \in \gV.
\end{equation}
If the graph is balanced, then the central values solve \eqref{eq:balance} at
$\varepsilon=0$.  We call \eqref{eq:balance} the \emph{balance equations}.

\begin{remark}
	If we see the surface as a soap film, then the saddle tower at $v(h)$ is
	pulled by a surface tension force along the wing of $h$, which can be
	calculated (up to a physical coefficient) as
	\[
	 	- \im \int_{A_h} (\Phi_1, \Phi_2, \Phi_3) = - 2\pi (\alpha_h, \beta_h, \gamma_h).
	\]
\end{remark}

\subsubsection{B-periods}

For any cycle $c=(h_1, \cdots, h_n)$ of the graph, let $B_c$ denote the
concatenation $B_{h_1}*\cdots * B_{h_n}$ which is a cycle in $\Sigma_t$.
Recall from Section \ref{sec:cut-cycle} the cycle basis
$\gC^*=\gF^*\cup\{c_1,c_2\}$.  We need to solve the following B-period problem:
\begin{gather}
	\varepsilon^2 \Big( \re \int_{B_c} \Phi_1 + \ii \re \int_{B_c} \Phi_2 \Big) =
	\begin{cases}
		0 & c \in \gF^*,\\
		T_i + \varepsilon^2 \Lambda_i & c=c_i,\; i = 1,2;
	\end{cases} \label{eq:BC12}\\
	\re \int_{B_c} \Phi_3 =
	\begin{cases}
		0 \pmod{2\pi} & c \in \gF^*,\\
		\Psi_i \pmod{2\pi} & c=c_i,\; i=1,2.
	\end{cases} \label{eq:BC3}
\end{gather}

\subsubsection{Conformality}

At $\varepsilon=0$ and the central value of all parameters, $\cent\Phi_1$ has
$\deg(v)$ simple poles in $\hat\C_v$, hence $\deg(v)-2$ zeros denoted
$\cent\zeta_{v,j}$ for $1\leq j\leq \deg(v)-2$.  Recall that these zeros are
simple (see the beginning of Section \ref{sec:construction}).  We can also
assume that they are not $\infty_v$.  When the parameters are close to their
central value, $\Phi_1$ has a simple zero $\zeta_{v,j}$ close to
$\cent\zeta_{v,j}$ in $\hat\C_v$ for $1\leq j\leq \deg(v)-2$.

\begin{remark}
	If the balancing equations are not solved, $\Phi_1$ may have a pole at
	$\infty_v$ and consequently, an extra zero near $\infty_v$. We may ignore
	this zero, as it disappears when the balancing equations are solved.  See
	Proposition~\ref{prop:alt-conformal} below.
\end{remark}

We will solve the following equations:
\begin{align}
	\int_{A_h} \frac{Q}{\Phi_1} &= 0, &&h \in \gH, \label{eq:conform1} \\
	\Res{\frac{Q}{\Phi_1}}{\zeta_{v,j}} &= 0, &&1 \le j \le \deg(v)-3,\; v \in \gV. \label{eq:conform2}
\end{align}

\begin{proposition}\label{prop:alt-conformal}
	The conformality equation~\eqref{eq:conformal} is solved if the
	equations~\eqref{eq:balance}, \eqref{eq:conform1}, and
	\eqref{eq:conform2} are solved.
\end{proposition}

\begin{proof}
 	If Equation \eqref{eq:balance} is solved, then $Q$ and $\Phi_1$ are
 	holomorphic at $\infty_v$ and $\Phi_1$ has no extra zero in $\hat\C_v$.  By
 	the Residue Theorem in $\hat\C_v$
 	\[
 		\sum_{h \in v} \int_{A_h} \frac{Q}{\Phi_1} + 
 		2\pi\ii \sum_{j=1}^{\deg(v)-2} \Res{\frac{Q}{\Phi_1}}{\zeta_{v,j}}=0.
 	\]
 	Hence if Equations \eqref{eq:conform1} and \eqref{eq:conform2} are solved,
 	the residue of $Q/\Phi_1$ at the last zero $\zeta_{v,\deg(v)-2}$ must also
 	vanish.  So the 1-form $Q/\Phi_1$, being holomorphic on $\Sigma$ with
 	vanishing $A$-periods, must be $0$.
\end{proof}

\subsubsection{Dimension count}

Let us perform a dimension count before proceeding further.

We have $3|\gE|$ real parameters $(\alpha,\beta,\gamma)$.  The complex
parameters $t$ comprises $|\gE|$ complex numbers.  For each $v \in \gV$,
M\"obius transforms on $\hat\C_v$ do not change the conformal structure of
$\Sigma_t$.  So we can fix the positions of three of the punctures
$(p_h)_{h \in v}$, leaving $\deg(v)-3$ free complex parameters for each
$v \in \gV$.  So, together with the parameter $\varepsilon$ of the family, we
have $9|\gE|-6|\gV|+1$ real parameters.

Let us now count the equations.  The B-period problems is $3g=3(|\gF|+1)$ real
equations. The conformality equations contain $|\gE| + 2 |\gE| - 3 |\gV| =
3|\gF|$ complex equations.  (This is of course the dimension of the space of
holomorphic quadratic differentials on $\Sigma_t$.) Moreover, we have
$3|\gV|-3$ balancing equations~\eqref{eq:balance}.  Hence there are
$9|\gF|+3|\gV| = 9|\gE|-6|\gV|$ real equations.  This is one less than the
number of parameters, so 1-parameter families are expected out of our
construction.

\subsection{Solving the conformality problem}

\begin{proposition}\label{prop:p}
	For $(t,\alpha,\beta,\gamma)$ in a neighborhood of
	$(0,\cent\alpha,\cent\beta,0)$, there exists $p = (p_h)_{h \in \gH}$
	depending analytically on $(t,\alpha,\beta,\gamma)$, such that the
	conformality equations~\eqref{eq:conform2} are solved.  Moreover,
	$p_h(0,\cent\alpha,\cent\beta,0) = \cent p_h$.
\end{proposition}

\begin{proof}
	As in Section \ref{sec:rigid-tower}, for each $v\in\gV$, we may fix the
	position of three points $p_h$ in $\hat\C_v$ using a M\"obius transformation.
	At $(t,\alpha,\beta,\gamma)=(0,\cent\alpha,\cent\beta,0)$,
	$(\Phi_1,\Phi_2,\Phi_3)$ restricted to $\hat\C_v$ is the Weierstrass data
	considered in Section \ref{sec:rigid-tower}, depending on the parameters
	$(p_h)_{h\in v}$.  By Theorem \ref{thm:rigid-tower}, the partial differential
	of
	\[
		\left(\Res{\frac{Q}{\Phi_1}}{\zeta_{v,j}}\right)_{1\leq j\leq \deg(v)-3}
	\]
	with respect to $(p_h)_{h\in v}$ is an isomorphism. So Proposition
	\ref{prop:p} follows from the Implicit Function Theorem.
\end{proof}

From now on, we assume that $p_h$ are given by Proposition~\ref{prop:p}.  We
make the change of parameters
\[
	\alpha_h + \ii\beta_h= - \rho_h \exp(\ii \theta_h).
\]
Clearly, $\rho=(\rho_h)_{h \in \gH} \in \cS$, and $(\exp(\ii \theta_h))_{h \in
\gH} \in \cA^2$.  The central values of $\rho_h$ is $\cent\rho_h = 1$ and the
central value of $\theta_h$ is $\cent\theta_h$ with
$\cent{u}_h=e^{\ii\cent\theta_h}$.

\begin{proposition}\label{prop:rhogamma}
	For $(t,\theta)$ in a neighborhood of $(0, \cent\theta)$, there exist unique
	values of $\rho$ and $\gamma$, depending real-analytically on $(t,\theta)$,
	such that the equations~\eqref{eq:conform1} are solved. At $t_h=0$ we have,
	no matter the values of $\theta$, $(t_k)_{k\neq \pm h}$, and other
	parameters, that
	\begin{equation}\label{eq:rhogamma1}
		\rho_h = 1\quad\text{and}\quad \gamma_h =0.
	\end{equation}
	In particular, we have $|\gamma_h|\leq C|t_h|$ for a uniform constant $C$.
	Moreover, at $(t,\theta)=(0,\cent{\theta})$, we have the Wirtinger
	derivatives
	\begin{equation}\label{eq:rhogamma2}
		\frac{\partial\rho_h}{\partial t_h} = -\frac{1}{2}\Upsilon_h\Upsilon_{-h}
		\quad\text{and}\quad 
		\frac{\partial\gamma_h}{\partial t_h} = -\frac{\ii}{2}\sigma_h\Upsilon_h\Upsilon_{-h}.
	\end{equation}
\end{proposition}
\begin{proof}
	Define for $h\in\gH$
	\[
		\mathcal E_h(t,\rho,\gamma,\theta)=\frac{1}{2\pi\ii}\int_{A_h}\frac{Q}{\Phi_1}.
	\]
	Assume that $t_h = 0$. Then $\Phi_1$, $\Phi_2$, and $\Phi_3$ have a simple
	pole at $p_h$, so $Q/\Phi_1$ has a simple pole at $p_h$ and, by the Residue
	Theorem
	\[
		\mathcal E_h\mid_{t_h=0} =
		\frac{\alpha_h^2+\beta_h^2+(\gamma_h-\ii\sigma_h)^2}{\alpha_h}=
		\frac{\rho_h^2+\gamma_h^2-1-2\ii\sigma_h\gamma_h}{\alpha_h}.
	\]
	So the solution of~\eqref{eq:conform1} at $t_h = 0$ is $(\rho_h, \gamma_h) =
	(1,0)$, no matter the values of the other parameters.  This
	proves~\eqref{eq:rhogamma1}. We compute the partial derivatives of $\mathcal
	E_h$ with respect to $\rho_h$ and $\gamma_h$ at
	$(t,\rho,\gamma,\theta)=(0,1,0,\cent{\theta})$:
	\begin{equation} \label{eq:rhogamma4}
		\frac{\partial \mathcal E_h}{\partial \rho_h} = \frac{2}{\cent\alpha_h},\qquad
		\frac{\partial\mathcal E_h}{\partial \gamma_h} = \frac{-2\ii\sigma_h}{\cent\alpha_h}.
	\end{equation}
	So the existence and uniqueness statement of the proposition follows from the
	Implicit Function Theorem.  To prove the last point, we need to compute the
	partial derivative of $\mathcal E_h$ with respect	to $t_h$ at
	$(t,\rho,\gamma,\theta)=(0,1,0,\cent{\theta})$. We use the following
	elementary results: if $f=a_{-1}z^{-1}+a_0+O(z)$ and
	$g=b_{-1}z^{-1}+b_0+O(z)$ are meromorphic functions with a simple pole at
	$z=0$, then
	\[
		\Res{\frac{f}{z^2g}}{0}=\frac{a_0b_{-1}-a_{-1}b_0}{(b_{-1})^2},\;
		\Res{f^2}{0}=2a_{-1}a_0,
		\;\text{and}\;
		\Res{f^2z}{0}=(a_{-1})^2.
	\]
	We have
	\begin{align}
		\frac{\partial\mathcal E_h}{\partial t_h}
		=& \Res{
		\sum_{j=1}^3 2\frac{\cent\Phi_j}{\cent\Phi_1} \frac{\partial\Phi_j}{\partial t_h}}{\cent p_h}\nonumber\qquad\text{ by the Residue Theorem}\\
		=& -2\sum_{j=1}^3 \Res{\frac{\cent\Phi_j\,dw_h}{\cent\Phi_1w_h^2}}{\cent p_h}
		\Res{\frac{\cent\Phi_j}{w_{-h}}} {\cent p_{-h}}\qquad \text{by Lemma~\ref{lem:derivative}}\nonumber\\
		=& -\frac{2}{\cent\alpha_h}
		\sum_{j=1}^3 \Res{ \frac{\cent\Phi_j}{w_h}}{\cent p_h} \Res{\frac{\cent\Phi_j}{w_{-h}}}{\cent p_{-h}} \label{eq:rhogamma5}\\
		&+ \frac{2}{(\cent\alpha_h)^2} \Res{\frac{\cent\Phi_1}{w_h}}{\cent p_h}
		\sum_{j=1}^3 \Res{\cent\Phi_j }{\cent p_h}\Res{\frac{\cent\Phi_j}{w_{-h}}}{\cent p_{-h}}.\label{eq:rhogamma6}
	\end{align}
	Since $\cent Q=0$, we have
	\begin{align}
		\Res{\frac{w_h \cent Q}{d w_h}}{\cent p_h}
		&= \sum_{j=1}^3 \Res{\cent\Phi_j}{\cent p_h}^2 = 0, \label{eq:rhogamma7}\\
		\Res{\frac{\cent Q}{d w_h}}{\cent p_h}
		&= 2\sum_{j=1}^3 \Res{\cent\Phi_j}{\cent p_h} \Res{\frac{\cent\Phi_j}{w_h}}{\cent p_h}= 0. \label{eq:rhogamma8}
	\end{align}

	Because $\Res{\cent\Phi_j}{\cent p_{-h}}= - \Res{\cent\Phi_j}{\cent p_h}$,
	Eq.~\eqref{eq:rhogamma8} implies that the term~\eqref{eq:rhogamma6} vanishes.
	Note that
	\[
		\frac{1}{\sqrt{2}}\Res{\cent\Phi}{\cent p_h}, \quad
		\frac{1}{\sqrt{2}}\Res{\overline{\cent\Phi}}{\cent p_h}, \quad
		\cent N(\cent p_h)
	\]
	form an orthonormal basis of $\C^3$ for the standard hermitian product
	$\langle\cdot,\cdot\rangle_\text{H}$.  After recalling the definition of
	$\Upsilon_h$ in Section \ref{sec:undulation}, we decompose in this basis
	\[
		\Res{\frac{\cent\Phi}{w_h}}{\cent p_h} =
		\frac{1}{2} \left\langle \Res{\cent\Phi}{\cent p_h}, \Res{\frac{\cent\Phi}{w_h}}{\cent p_h} \right\rangle_{\!\text{H}} \Res{\cent\Phi}{\cent p_h}
		\;+\;\sigma_h\Upsilon_h\cent N(\cent p_h).
	\]
	Here, the component on $\Res{\overline{\cent\Phi}}{\cent p_h}$ vanishes
	because of~\eqref{eq:rhogamma8}.  In the same way, after recalling that
	$\cent{N}(\cent p_{-h})=\cent{N}(\cent p_h)$ and $\sigma_{-h}=-\sigma_h$,
	\[
		\Res{\frac{\cent\Phi}{w_{-h}}}{\cent p_{-h}} =
		\frac{1}{2} \left\langle \Res{\cent\Phi}{\cent p_{h}}, \Res{\frac{\cent\Phi}{w_{-h}}}{\cent p_{-h}} \right\rangle_{\!\text{H}} \Res{\cent\Phi}{\cent p_h}
		\;-\;\sigma_{h}\Upsilon_{-h}\cent N(\cent p_h).
	\]

	Hence by Equation ~\eqref{eq:rhogamma5}
	\[
		\frac{\partial\mathcal E_h}{\partial t_h}=
		\frac{-2}{\cent \alpha_h}\left\langle\overline{\Res{\frac{\cent\Phi}{w_h}}{\cent p_h}}, \Res{\frac{\cent\Phi}{w_{-h}}}{\cent p_{-h}}\right\rangle_{\!\text{H}}
		= \frac{2}{\cent\alpha_h}\Upsilon_h\Upsilon_{-h}.
	\]
	If $k\neq \pm h$, we have $\partial\mathcal E_h/\partial t_k=0$ as
	$\partial\Phi_j/\partial t_k$ is holomorphic at $\cent p_h$.

	We differentiate $\mathcal E_h(t,\rho(t,\theta),\gamma(t,\theta),\theta)=0$
	with respect to $t_h$ and obtain, using Equation \eqref{eq:rhogamma4},
	\begin{equation} \label{eq:rhogamma9}
		\frac{2}{\cent\alpha_h}\left(\Upsilon_h\Upsilon_{-h}+\frac{\partial\rho_h}{\partial t_h}-\ii\sigma_h\frac{\partial\gamma_h}{\partial t_h}\right)=0.
	\end{equation}
	Since $\mathcal E_h$ is holomorphic in $t_h$, we obtain by differentiation
	with respect to $\overline{t_h}$ and conjugation
	\begin{equation}
		\label{eq:rhogamma10}
		\frac{2}{\cent\alpha_h}\left(\frac{\partial\rho_h}{\partial t_h}+\ii\sigma_h\frac{\partial\gamma_h}{\partial t_h}\right)=0.
	\end{equation}
	Solving the system \eqref{eq:rhogamma9}, \eqref{eq:rhogamma10} gives
	Equation \eqref{eq:rhogamma2}.
\end{proof}

From now on, we assume that $\rho$ and $\gamma$ are given by
Proposition~\ref{prop:rhogamma}.

\subsection{Solving horizontal balance and period problems} \label{ssec:horizontal-balance}

We make the change of parameters
\begin{equation}\label{eq:t}
	t_h = -\exp\Big(-\ell_h\varepsilon^{-2}+\ii \sigma_h \phi_h\Big),
\end{equation}
where $\ell = (\ell_h)_{h \in \gH} \in \cS$, and $\phi = (\phi_h)_{h \in \gH}
\in \cA$.  The central value of $\ell_h$ is $\cent\ell_h$, the length of the
edge $e(h)$. The central value of $\phi_h$ is the given phase function
$\cent\phi_h$.  We combine $\ell$ and $\theta$ into
\[
	x_h = \ell_h e^{\ii\theta_h},
\]
whose central value is $\cent x_h$ as given by the graph.

\begin{proposition}\label{prop:thetaell}
	Assume that the graph $\gG$ is rigid and balanced.  For $(\varepsilon,\phi)$
	in a neighborhood of $(0,\cent\phi)$, there exist unique values for $(x_h)_{h
	\in \gH}$, depending smoothly on $(\varepsilon,\phi)$, such that $x(0,\phi) =
	\cent x$ for all $\phi$ and the horizontal B-period equations \eqref{eq:BC12}
	as well as the $\alpha$ and $\beta$ components of the balance
	equations~\eqref{eq:balance} are solved.  Moreover, $x$ is an even function
	of $\varepsilon$ and at $\varepsilon=0$, we have for all $\phi$
	\[
		\frac{\partial^2 x_h}{\partial \varepsilon^2}(0,\phi) = 2 \xi_h
	\]
	where $\xi_h$ is the solution to~\eqref{eq:system}.
\end{proposition}

\begin{proof}
	Define for $(\varepsilon,x,\phi)$ in a neighborhood of $(0,\cent
	x,\cent\phi)$ and $h\in\gH$
	\[
		\cPhor_h(\varepsilon,x,\phi)= 
		\varepsilon^2\Big(\re \int_{B_h} \Phi_1 + \ii \re \int_{B_h} \Phi_2\Big).
	\]
	By Lemma~\ref{lem:Bintegral}
	\[
		\lambda^\text{hor}_h(\varepsilon,x,\phi):=\bigg[\Big(\re\int_{B_h}\Phi_1\Big)-\alpha_h\log|t_h|\bigg]+\ii\bigg[\Big(\re\int_{B_h}\Phi_2\Big)-\beta_h\log|t_h|\bigg]
	\]
	extends at $\varepsilon=0$ to a smooth function of remaining parameters
	$(\varepsilon,x,\phi)$.

	Recall the definition of $\mu_h$ in Section \ref{sec:undulation} and that
	$z_0=O_v$ is chosen as the origin of the parameterization of the saddle tower
	$\sS_v$.  By the last statement of Lemma~\ref{lem:Bintegral}, we have at the
	central value
	\[
		\lambda^\text{hor}_h(0,\cent x,\phi)=\mu_h-\mu_{-h}.
	\]
	We then have
	\begin{align*}
		\cPhor_h(\varepsilon,x,\phi)=&\varepsilon^2\left[(\alpha_h+\ii\beta_h)\log|t_h|+\lambda^\text{hor}_h(\varepsilon,x,\phi)\right]\\
		=&-\varepsilon^2\rho_h e^{\ii\theta_h}(-\ell_h\varepsilon^{-2})+\varepsilon^2\lambda^\text{hor}_h(\varepsilon,x,\phi)\\
		=&\rho_h x_h+\varepsilon^2\lambda^\text{hor}_h(\varepsilon,x,\phi).
	\end{align*}
	Note that by definition, $\cPhor_h$ is an even function of $\varepsilon$.
	Since $\rho_h$ is an analytic function of $t_h$ with value $1$ at $t_h=0$ and
	all derivatives of $t_h$ with respect to $\varepsilon$ vanish at
	$\varepsilon=0$, we have at $\varepsilon=0$
	\begin{equation} \label{eq:Phor0}
		\cPhor_h(0,x,\phi) = x_h \quad\text{and}\quad
		\frac{\partial^2 \cPhor_h}{\partial\varepsilon^2}(0,\cent x,\phi) = 2(\mu_h-\mu_{-h}).
	\end{equation}

	Now define for $c$ in the cycle basis $\gC^*$
	\[
		\cPhor_c(\varepsilon,x,\phi) = 
		\varepsilon^2\Big(\re \int_{B_c} \Phi_1 + \ii \re \int_{B_c} \Phi_2\Big) =
		\sum_{h \in c} \cPhor_h(\varepsilon,x,\phi).
	\]
	Then, at $\varepsilon=0$, we have
	\begin{gather*}
		\cPhor_c(0,x,\phi) = \sum_{h \in c} x_h = \Phor_c(x) \quad\text{and}\\
		\frac{\partial^2 \cPhor_c}{\partial\varepsilon^2}(0,\cent x,\phi) = 2\sum_{h \in c} (\mu_h - \mu_{-h}) = 2\Phor_c(\mu^a).
	\end{gather*}

	Regarding the horizontal balance equation \eqref{eq:balance}, we define for
	$b$ in the cut-basis $\gB^*$
	\[
		\cFhor_b = - \div_b(\alpha + \ii\beta) =\div_b(\rho_he^{\ii\theta_h})
		=\div_b(\rho_h u_h(x)).
	\]
	Again, $\cFhor_b$ is an even function of $\varepsilon$ and at
	$\varepsilon=0$, we have since $\rho=1$
	\[
		\cFhor_b(0,x,\phi) =\Fhor_b(x)\quad\text{and}\quad
		\frac{\partial^2 \cFhor_b}{\partial\varepsilon^2}(0,\cent x,\phi) = 0.
	\]
	We want to solve the system
	\[
		\begin{cases}
			\cPhor_c(\varepsilon,x,\phi) = 0, &c \in \gF^*,\\
			\cPhor_{c_i}(\varepsilon,x,\phi) = T_i + \varepsilon^2 \Lambda_i, &i = 1,2,\\
			\cFhor_b(\varepsilon,x,\phi) = 0, &b \in \gB^*.
		\end{cases}
	\]
	If $\gG$ is balanced, the system is solved at $\varepsilon=0$ with
	$x(0,\phi)=\cent x$ for all $\phi$.  If $\gG$ is rigid, then by Implicit
	Function Theorem, for $\varepsilon$ in a neighborhood of $0$, there exists a
	unique solution $x(\varepsilon,\phi)$ depending smoothly on
	$(\varepsilon,\phi)$ such that $x(0,\phi)=\cent x$.

	Moreover, since the system is even in $\varepsilon$, so is the unique
	solution $x(\varepsilon,\phi)$.  Taking the second derivative of the system
	with respect to $\varepsilon$ at $\varepsilon=0$ gives
	\[
		\begin{dcases}
			2 \Phor_c(\mu^a) + \Phor_c\Big(\frac{\partial^2 x}{\partial\varepsilon^2}(0, \phi)\Big) = 0, &c \in \gF^*,\\
			2 \Phor_{c_i}(\mu^a) + \Phor_{c_i}\Big(\frac{\partial^2 x}{\partial\varepsilon^2}(0, \phi)\Big) = 2\Lambda_i, &i = 1,2,\\
			D\Fhor_b(\cent x) \cdot \frac{\partial^2 x}{\partial\varepsilon^2}(0, \phi) = 0, &b \in \gB^*.
		\end{dcases}
	\]
	Its unique solution is, by definition, $\frac{\partial^2 x}{\partial\varepsilon^2}(0, \phi) = 2\xi$.
\end{proof}

From now on, we assume that the parameters $x_h$ are given by
Proposition~\ref{prop:thetaell}.

\subsection{Solving vertical balance and period problems}

\begin{proposition}\label{prop:phi}
	Assume that the phase function $\cent\phi$ is balanced and rigid.  For
	$\varepsilon$ in a neighborhood of $0$, there exist unique values for
	$(\phi_h)_{h \in \gH}$, depending smoothly on $\varepsilon$, such that
	$\phi_h(0) = \cent\phi_h$ and the vertical B-period problems~\eqref{eq:BC3}
	as well as the $\gamma$ component of~\eqref{eq:balance} are solved.
\end{proposition}

\begin{proof}
	Similarly to the proof of Proposition \ref{prop:thetaell}, we define for $h\in\gH$
	\[
		\cPver_h(\varepsilon,\phi) = \re \int_{B_h} \Phi_3.
	\]
	By Lemma~\ref{lem:Bintegral},
	\[
		\lambda^\text{ver}_h(\varepsilon,\phi):=\re\bigg[\Big(\int_{B_h}\Phi_3\Big)-(\gamma_h-\ii\sigma_h)\log t_h \bigg]
	\]
	extends smoothly at $\varepsilon=0$. By Proposition \ref{prop:rhogamma}, we
	have $\gamma_h=O(t_h)$ so $\gamma_h\log t_h $ extends smoothly at
	$\varepsilon=0$ with value $0$. Hence by Lemma \ref{lem:Bintegral}, we have
	at $\varepsilon=0$
	\begin{align*}
		\lambda^\text{ver}_h(0,\phi)=&\re
		\lim_{z\to p_h} \bigg[
			\Big(\int_{O_{v(h)}}^z\cent\Phi_3\Big)+\ii\sigma_h\log w_h(z)
		\bigg]\\
		&-\re\lim_{z\to p_{-h}} \bigg[
			\Big(\int_{O_{v(-h)}}^z\cent\Phi_3\Big)+\ii\sigma_{-h}\log w_{-h}(z)
		\bigg]\\
		=&\nu_h-\nu_{-h}
	\end{align*}
	where $\nu_h$ is defined as in Section \ref{sec:undulation}.  We then have
	\[
		\cPver(\varepsilon,\phi)=\lambda^\text{ver}_h(\varepsilon,\phi)-\gamma_h\ell_h\varepsilon^{-2}+\phi_h+\sigma_h\pi.
	\]
	Since $\gamma_h=O(t_h)$, we have at $\varepsilon=0$
	\begin{equation}
		\label{eq:Pver0}
		\cPver_h(0,\phi)= \nu_h-\nu_{-h} + \phi_h +\sigma_h\pi.
	\end{equation}
	Now define for $c$ in the cycle-basis $\gC^*$
	\[
		\cPver_c (\varepsilon,\phi)= \re \int_{B_c} \Phi_3 = \sum_{h \in c}\cPhor_h(\varepsilon,\phi).
	\]
	Observe that by Proposition~\ref{prop:adapted}
	\[
		\sum_{h\in c} (\nu_h-\nu_{-h}+\pi)=\sum_{h\in c}(\nu_h-\nu_{\rotate(h)}+\pi)=0\pmod{2\pi}
	\]
	Hence at $\varepsilon=0$, we have
	\[
		\cPver_c(0,\phi) = \Pver_c(\phi) \pmod{2\pi}.
	\]
	Regarding the vertical balance equation \eqref{eq:balance}, we define for $b$
	in the cut basis $\gB^*$
	\[
		\cFhor_b(\varepsilon,\phi) = -\div_b(\gamma(\varepsilon,\phi))=-\sum_{h\in b}
		\gamma_h(\varepsilon,\phi).
	\]
	Using the Mean Value Inequality and Proposition~\ref{prop:rhogamma}, we have for $\varepsilon\to 0$
	\[
		\gamma_h(\varepsilon,\phi) = \sigma_h \Upsilon_h \Upsilon_{-h} \im(t_h(\varepsilon,\phi)) + o(t_h(\varepsilon,\phi)).
	\]
	By Proposition~\ref{prop:thetaell} and that $\partial x_h / \partial
	\varepsilon = 0$, we have
	\[
		\frac{\partial^2 x_h}{\partial\varepsilon^2}(0) = \exp(\ii\cent\theta_h)
		\Big(\frac{\partial^2 \ell_h}{\partial\varepsilon^2}(0) +
		\ii\ell_h\frac{\partial^2 \theta_h}{\partial\varepsilon^2}(0)\Big) = 2\xi_h,
	\]
	so
	\[
		\frac{\partial^2\ell_h}{\partial\varepsilon^2}(0) =
		2 \re(\xi_h\exp(-\ii\cent\theta_h)).
	\]
	Therefore using Taylor formula
	\begin{align*}
		\im(t_h(\varepsilon,\phi)) =& -\sigma_h \sin(\phi_h) \exp(-\ell_h(\varepsilon,\phi)\varepsilon^{-2})\\
		=& -\sigma_h \sin(\phi_h) \exp\left(-\cent\ell_h\varepsilon^{-2} - \re(\xi_h\exp(-\ii\cent\theta_h)) + O(\varepsilon^2)\right).
	\end{align*}
	Hence we have
	\[
		\cFver_b(\varepsilon,\phi)= \sum_{h\in b}\Upsilon_h \Upsilon_{-h} \sin(\phi_h) \exp\left(-\cent\ell_h\varepsilon^{-2} - \re(\xi_h\exp(-\ii\theta_h)) + O(\varepsilon^2)\right).
	\]

	Recall that $\cent\ell_b = \min\{\cent\ell_h \mid h \in b\}$ and define for
	$\varepsilon\neq 0$
	\[
		\widetilde{\cFver_b}(\varepsilon,\phi) = \exp(\cent\ell_b\varepsilon^{-2}) \cFver_b(\varepsilon,\phi).
	\]
	Then $\widetilde\cFver$ extends smoothly at $\varepsilon=0$ and in the limit
	$\varepsilon\to 0$ only the shortest edges remain:
	\[\widetilde\cFver(0,\phi) = \sum_{h\in m(b)}\Upsilon_h\Upsilon_{-h}\sin(\phi_h)\exp\left(-\re(\xi_h\exp(-\ii\theta_h))\right)=\Fver(\phi).\]

	We want to solve the system
	\[
		\begin{cases}
			\cPver_c(\varepsilon,\phi) = 0, &c \in \gF^*,\\
			\cPver_{c_i}(\varepsilon,\phi) = \Psi_i , & i = 1,2,\\
			\widetilde{\cFver_b}(\varepsilon,\phi) = 0, &b \in \gB^*.
		\end{cases}
	\]
	If the phase function $\cent\phi$ is balanced, the system is solved at
	$\varepsilon=0$ by $\phi=\cent\phi$.  If $\cent\phi$ is rigid and the cut
	basis $\gB^*$ is $\gB_m^*$ as given by Proposition \ref{proposition:mdiv},
	then by the Implicit Function Theorem, the system has a unique solution
	$\phi(\varepsilon)$, depending smoothly on $\varepsilon$ for $\varepsilon$ in
	a neighborhood of $0$, such that $\phi(0) = \cent\phi$.
\end{proof}

\subsection{Geometry and Embeddedness}

We have constructed a 1-parameter family of Weierstrass data
$(\Phi_{1,\varepsilon},\Phi_{2,\varepsilon},\Phi_{3,\varepsilon})$ depending on
$\varepsilon>0$ that solve the Conformality and Period Problems. All parameters
are now smooth functions of $\varepsilon$ and are denoted accordingly
$t_h(\varepsilon)$, $\theta_h(\varepsilon)$, etc...  We denote
\[
	f_{\varepsilon}=(f_{\varepsilon}^\text{hor},f_{\varepsilon}^\text{ver}):\Sigma_{t(\varepsilon)}\to\T^3_{\varepsilon}
\]
the immersion given by Weierstrass Representation formula, decomposed into its
horizontal and vertical components.  Note that $\Sigma_{t(\varepsilon)}$ and
$f_{\varepsilon}$ are independent of the positive number $\delta$ used to
define $\Sigma_t$ (although the smaller $\delta$, the smaller $\varepsilon$
must be in the definition of $\Sigma_t$).  We first prove that the immersion is
regular, that is:

\begin{proposition}\label{prop:regular}
	we have
	\begin{equation}\label{eq:regular}
		|\Phi_{1,\varepsilon}|^2+|\Phi_{2,\varepsilon}|^2+|\Phi_{3,\varepsilon}|^2 > 0
	\end{equation}
	on $\Sigma_{t(\varepsilon)}$ for sufficiently small $\varepsilon$.
\end{proposition}

\begin{proof}
	The statement holds in $U_\delta$ because the Saddle towers are regular
	embeddings, and the inequality~\eqref{eq:regular} remains for sufficiently
	small $\varepsilon$.  It suffices to prove~\eqref{eq:regular} on the annuli
	that we identified to form necks.  This is the case if $\Phi_{3,\varepsilon}$
	has no zero in the necks.

	Recall that $\cent\Phi_3$ has $\deg(v)$ poles hence $\deg(v)-2$ zeros in
	$\hat\C_v$.  By choosing $\delta$ sufficiently small, we may assume that all
	$2|\gE|-2|\gV|=2|\gF|=2g-2$ zeros lie in~$U_\delta$, which remains true for
	sufficiently small $\varepsilon$.  But for $\varepsilon\neq 0$,
	$\Phi_{3,\varepsilon}$, being a holomorphic function on a Riemann surface of
	genus $g$ (without nodes), has exactly $2g-2$ zeros, therefore no other
	zeros.  The statement then follows readily.
\end{proof}

It remains to prove that the image of $f_{\varepsilon}$ has the desired
geometry, as stated in Points (1), (2), (3) of Theorem \ref{thm:main}, and is
embedded.  Define
\[
	X_v(\varepsilon)=f_{\varepsilon}^\text{hor}(O_v).
\]
Then $f_{\varepsilon}-X_v(\varepsilon)$ converges smoothly on $U_{v,\delta}$ to
$f_{v,0}$ which parametrize a Saddle tower which we denote $\sS_v$.  Hence
$f_{\varepsilon}(U_{v,\delta})$ is embedded for $\varepsilon>0$ small enough.
Let $\cent{\sS}_v$ be the reference saddle tower, namely translated so that the
image of $O_v$ is the origin. We denote $\varphi(\sS)$ the phase of a saddle
tower.  Using Equation \eqref{eq:Pver0} and Proposition \ref{prop:adapted}, we
have for $h\in\gH$ \def\phase{\varphi}
\begin{align*}
	\phase(\sS_{v(-h)})-\phase(\sS_{v(h)})&=
	\phase(\cent\sS_{v(-h)})-\phase(\cent\sS_{v(h)})
	+\lim_{\varepsilon\to 0}\cPver_h(\varepsilon)\\
	&= \phase(\cent\sS_{v(-h)})-\phase(\cent\sS_{v(h)})+\cent\phi_h+\nu_h-\nu_{-h}+\sigma_h\pi\\
	&=\cent\phi_h\pmod{2\pi}.
\end{align*}
This proves Point (3).  Using Equation \eqref{eq:Phor0}, we have for $h\in\gH$
\[
	\lim_{\varepsilon\to
	0}\varepsilon^2\left(X_{v(-h)}(\varepsilon)-X_{v(h)}(\varepsilon)\right)
	=\lim_{\varepsilon\to 0}\cPhor_h(\varepsilon)=\cent x_h.
\]
This proves Point (2). Moreover, because the vertices are represented by
distinct points, the images $f_{\varepsilon}(U_{v,\delta})$ for $v\in\gV$ are
disjoint for $\varepsilon>0$ small enough.

For $\delta$ and $\varepsilon$ small enough, the image of the circle
$|w_h|=\delta$ is close to a vertical segment $s_h(\varepsilon)$ lying in
$f_{\varepsilon}(U_{v,\delta})$.  On the boundary of the annulus
$|t_h|/\delta<|w_h|<\delta$, the gauss map is close to a horizontal constant
(from the convergence to saddle towers), so by the maximum principle for
holomorphic functions, it is close to this constant on the whole annulus.
Given the boundary behavior, the image of this annulus is a graph over a
vertical plane, hence embedded.  By the convex hull property of minimal
surfaces, the image of the annulus $|t_h|/\delta<|w_h|<\delta$ stays close to
the vertical plane $B_{e(h)}(\varepsilon)$ bounded by $s_h(\varepsilon)$ and
$s_{-h}(\varepsilon)$.  After a scaling by $\varepsilon^2$, Point (1) follows.

If there are no parallel edges, by taking $\delta>0$ small enough, we may
ensure that the planes $B_e(\varepsilon)$ for $e\in\gE$ are disjoint and
conclude that the image of $f_{\varepsilon}$ is embedded because the graph is
embedded.  The situation is more subtle in the case of parallel edges.  Assume
that the edges $e(h_1),\cdots,e(h_k)$ are parallel, with $h_1,\cdots,h_k\in v$
and $\rotate(h_i)=h_{i+1}$ for $1\leq i\leq k-1$. Because the saddle tower
$\sS_v$ is embedded, the corresponding ends are asymptotic to parallel but
distinct vertical planes.  Hence, by our choice of ordering of the ends, the
abscissa of the vertical segments $s_{h_i}(\varepsilon)$ are decreasing and
separated by a distance greater than some uniform $r>0$.  On the other side,
the abscissa of the vertical segments $s_{-h_i}(\varepsilon)$ are also
decreasing and uniformly separated.  Hence for $\varepsilon$ small enough, the
planes $B_{e(h_i)}$ for $1\leq i\leq k$ are disjoint, and we conclude again that
the image of $f_{\varepsilon}$ is embedded.  This concludes the proof of
Theorem \ref{thm:main}.

\appendix

\section{Rigidity of saddle towers} \label{app:rigid-tower}

This section is dedicated to a proof of Theorem \ref{thm:rigid-tower}.  Let
$\dot{p}\in\operatorname\ker{D\Lambda(\cent{p})}$.  We want to prove that, if
$\dot{p}_h = 0$ for three $h \in \gH$, then $\dot{p}=0$.

Define $p(\epsilon)=\cent{p}+\epsilon\dot{p}$. Let $(\Phi_1(\epsilon),
\Phi_2(\epsilon), \Phi_3(\epsilon))$ be the Weierstrass data given by Equation
\eqref{eq:weierstrassz} and $X(\epsilon)$ the corresponding immersion given by
Equation \eqref{eq:weierstrass-representation}. Note that $X(\epsilon)$ is
harmonic but not minimal. Let $H(\epsilon)$ be the mean curvature and
$N(\epsilon)$ be the normal of $X(\epsilon)$.  To ease computation, the
dependence on $\epsilon$ will not be written.  Instead, we use an exponent
$\circ$ to denote the value at $\epsilon=0$ and a dot to denote the partial
derivative with respect to $\epsilon$ at $\epsilon=0$.
\begin{lemma} \label{lem:rigid-tower}
	We have $\dot{H}=0$ on the Riemann sphere minus the points $\cent{p}_h$ for
	$h\in\gH$.
\end{lemma}
\begin{proof}
	We have for all $\epsilon$
	\[
		Q=\sum_{i=1}^{n-2}\lambda_i\frac{\Phi_1\,dz}{z-\zeta_i}\quad\mbox{ with }\quad
		\lambda_i=\Res{\frac{Q}{\Phi_1}}{\zeta_i}.
	\]
	We have $\cent{\lambda}_i=0$ and $\dot{\lambda}_i=D\Lambda_i(\cent{p})\cdot
	\dot{p}=0$ so $\dot{Q}=0$.  Let $w=x_1+\ii x_2$ be a local complex coordinate
	on the Riemann sphere.  Recall the standard formula for the mean curvature
	\begin{equation}
		\label{equation:H}
		H=\frac{g_{22}b_{11}+g_{11}b_{22}-2g_{12}b_{12}}{2(g_{11}g_{22}-g_{12}^2)}
	\end{equation}
	where $g$ and $b$ are respectively the matrices of the first and second
	fundamental forms in the coordinate system $(x_1,x_2)$.  We have for all
	$\epsilon$
	\[
		Q=\left(\left\|\frac{\partial X}{\partial x_1}\right\|^2-\left\|\frac{\partial X}{\partial x_2}\right\|^2-2\ii
		\left\langle\frac{\partial X}{\partial x_1},\frac{\partial X}{\partial x_2}\right\rangle\right)dw^2
		=(g_{11}-g_{22}-2\ii g_{12})dw^2.
	\]
	Hence from $\cent{Q}=0$ and $\dot{Q}=0$ we obtain
	\begin{equation} \label{eq:dotg}
		\cent{g}_{11}=\cent{g}_{22},\quad
		\cent{g}_{12}=0,\qquad
		\dot{g}_{11}=\dot{g}_{22}
		\quad\text{and}\quad
		\dot{g}_{12}=0.
	\end{equation}
	Since $X(\epsilon)$ is harmonic for all $\epsilon$, we have $b_{11}+b_{22}=0$ for
	all $\epsilon$, so
	\[
		\cent{b}_{11}+\cent{b}_{22}=0\quad\mbox{ and }\quad \dot{b}_{11}+\dot{b}_{22}=0.
	\]
	Taking the derivative of \eqref{equation:H}, we obtain $\dot{H}=0$.
\end{proof}
By Lemma \ref{lem:rigid-tower}, $u=\langle\dot{X},\cent{N}\rangle$ is a Jacobi
field.  In a neighborhood of $\cent{p}_h$, we have
\begin{equation}
	\label{eq:dotX}
	\dot{X}=\re\left[\left(\cos(\theta_h),\sin(\theta_h),i\sigma_h\right)\frac{\dot{p}_h}{z-\cent{p}_h}\right]+\mbox{bounded terms}
\end{equation}
\[
	\cent{N}=\sigma_h(-\sin(\theta_h),\cos(\theta_h),0)+O(z-\cent{p}_h)
\]
so $u$ is a bounded Jacobi field.  Since all saddle towers are rigid, there
exists $c\in\R^3$ such that $u=\langle \cent{N}, c \rangle$.
Consider the translated immersion $Y(\epsilon) = X(\epsilon)-\epsilon c$.
Then $\langle\dot{Y},\cent{N}\rangle=0$ so
$\dot{Y}$ is a tangent vector.  Using the local complex coordinate $w=x_1+\ii
x_2$, we decompose $\dot{Y}$ in the tangent space as
\[
	\dot{Y}=\xi_1\frac{\partial \cent{X}}{\partial
	x_1}+\xi_2\frac{\partial\cent{X}}{\partial x_2}
	=\xi\frac{\partial\cent{X}}{\partial
	w}+\overline{\xi}\frac{\partial\cent{X}}{\partial\overline{w}}\quad\mbox{
	with } \quad \xi=\xi_1+\ii \xi_2.
\]
\begin{lemma} \label{lem:rigid-tower2}
	$\xi\frac{d}{dw}$ defines a holomorphic vector field on the Riemann sphere.
	Moreover, $\xi(\cent{p}_h)=-\dot{p}_h$ for $h\in\gH$.
\end{lemma}
\begin{proof}
	Since $\cent{X}$ is conformal and harmonic, we have
	\[
		\left\langle\frac{\partial\cent{X}}{\partial w},\frac{\partial\cent{X}}{\partial w}\right\rangle=0,\qquad
		\left\langle\frac{\partial\cent{X}}{\partial w},\frac{\partial^2\cent{X}}{\partial w^2}\right\rangle=0,
		\qquad\text{and}\quad
		\frac{\partial^2\cent{X}}{\partial w\partial \overline{w}}=0
	\]
	where $\langle\cdot,\cdot\rangle$ denotes the $\C$-bilinear dot product.
	For all $\epsilon$
	\[
		\left\langle\frac{\partial X}{\partial w},\frac{\partial X}{\partial w}\right\rangle
		=\frac{1}{4}(g_{11}-g_{22}-2\ii g_{12}).
	\]
	Taking the derivative with respect to $\epsilon$ and using Equation
	\eqref{eq:dotg}, we obtain
	\[
		2\left\langle\frac{\partial\cent{X}}{\partial w},\frac{\partial \dot{Y}}{\partial w}\right\rangle
		=
		2\left\langle\frac{\partial\cent{X}}{\partial w},\frac{\partial}{\partial w}
		\left(\xi\frac{\partial\cent{X}}{\partial w}+\overline{\xi}\frac{\partial\cent{X}}{\partial\overline{w}}\right)\right\rangle=
		2 \frac{\partial\overline{\xi}}{\partial w}
		\left\langle\frac{\partial \cent{X}}{\partial w},\frac{\partial\cent{X}}{\partial\overline{w}}\right\rangle=0.
	\]
	Hence $\frac{\partial\overline{\xi}}{\partial w}=0$ and $\xi$ is holomorphic.
	If $w'$ is another local complex coordinate, we can write 
	\[
		\dot{Y}=\xi'\frac{\partial\cent{X}}{\partial w'}+\overline{\xi'}\frac{\partial \cent{X}}{\partial\overline{w'}}
		=\xi'\frac{\partial\cent{X}}{\partial w}\frac{dw}{d w'}
		+\overline{\xi'}\frac{\partial\cent{X}}{\partial\overline w}\overline{\left(\frac{dw}{d w'}\right)}
	\]
	This gives
	\[
		\xi=\xi'\frac{dw}{d w'}
	\]
	so $\xi$ transforms as a holomorphic vector field under change of coordinate.
	Using the complex coordinate $w=z-\cent{p}_h$ in a neighborhood of $\cent{p}_h$, we
	have
	\[
		\frac{\partial\cent{X}}{\partial w}=\frac{1}{2w}(-\cos(\theta_h),-\sin(\theta_h),-\ii\sigma_h)+O(1).
	\]
	Hence by Equation \eqref{eq:dotX},
	\[
		\xi(\cent{p}_h)=-\dot{p}_h.
	\]
\end{proof}

Recall that we have fixed the position of three points $p_h$. Hence $\xi$ has
at least three zeros. Now a non-zero holomorphic vector field on the Riemann
sphere has two zeros, so $\xi=0$. This implies that $\dot{p}=0$ and concludes
the proof of Theorem~\ref{thm:rigid-tower}.

\section{Horizontal rigidity of ``triangulated'' graphs} \label{appendix:triangles}

This section is dedicated to a proof of Theorem~\ref{theorem:triangles}.
Consider, for $x\in\cH$, the total length
\[
	{\mathcal L}(x)=\frac{1}{2}\sum_{h\in\gH}\|x_h\|.
\]
If $\Phor(\chi)=0$, so that $\chi=\grad f$, we have
\begin{align*}
	D{\mathcal L}(\cent{x})\cdot \chi
	=&\frac{1}{2}\sum_{h\in\gH}\frac{1}{\|\cent{x}_h\|}\langle \cent{x}_h,f_{v(-h)}-f_{v(h)}\rangle\\
	=&-\sum_{h\in\gH}\frac{1}{\|\cent{x}_h\|}\langle \cent{x}_h,f_{v(h)}\rangle
	\quad\mbox{ (using $h\to -h$ for $f_{v(-h)}$)}\\
	=&-\sum_{v\in V}\sum_{h\in b(v)}\frac{1}{\|\cent{x}_h\|}\langle\cent{x}_h,f_v\rangle\\
	=&-\sum_{v\in V}\langle\Fhor_{b(v)},f_v\rangle.
\end{align*}
We have proved that
\begin{proposition}\label{prop:minlength}
	The graph is balanced if and only if $\cent{x}$ is a critical point of
	$\mathcal L$ restricted to those $x\in\cA$ such that
	$\Phor(x)=\Phor(\cent{x})$.
\end{proposition}
Assume that $\Fhor(\cent{x})=0$ and let $\chi=\grad(f)$ be in the kernel of
$(D\Fhor(\cent{x}),\Phor)$. Differentiating the above equation
\[
	D^2{\mathcal L}(\cent{x})\cdot(\chi,\chi)=
	-\sum_{v\in V}\langle D\Fhor_{b(v)}(\cent{x})\cdot \chi,f_v\rangle =0.
\]
On the other hand, a direct computation gives
\[
	D^2{\mathcal L}(\cent{x})\cdot(\chi,\chi)=\frac{1}{2}\sum_{h\in\gH}
	\frac{1}{\|\cent{x}\|^3}\left(\|\cent{x}_h\|^2\|\chi_h\|^2-\langle \cent{x}_h,\chi_h\rangle^2\right).
\]
The summands are all non-negative, hence must be all zero, which means that
$\chi_h$ is parallel to $\cent{x}_h$ for all $h\in\gH$.

If all faces have 2 or 3 edges, this implies that $\chi_h=\lambda \cent{x}_h$
for some $\lambda\in\R$. Then $\Phor_{c_1}(\chi)=\lambda T_1=0$ implies
$\lambda=0$ and $\chi=0$, which proves Theorem \ref{theorem:triangles}.

\section{Classification of balanced configurations of genus 3}\label{app:TPMSg3}

As promised in Section~\ref{ssec:exg3}, we prove the following Classification
Theorem.

\begin{theorem}\label{thm:TPMSg3}
	The Meeks, aG, aH, and aI configurations in Figure~\ref{fig:g3} are the only
	balanced configurations whose graphs are orientable with two faces.  Hence
	they are the only possible configurations that give rise to TPMSs of genus 3.
\end{theorem}

\begin{proof}
	If a configuration gives rise to a TPMS of genus three, its graph must have
	two faces.  By Euler's formula, the average degree of the graph is
	\[
		2|\gE|/|\gV|=2(|\gV|+|\gF|)/|\gV| = 2+4/|\gV|.
	\]
	By Assumption~\ref{ass:degree}, the average degree is at least $4$, hence the
	graph has at most $2$ vertices.  We discuss two cases

	\begin{description}
		\item[$|\gV|=1$, $|\gE|=3$]

			Then all edges are loops.  Since edges are represented by straight segments,
			no loop is null-homologous.  If two loops are homologous, they must be
			parallel, and represented by the same segment that cuts the torus into an
			annulus.  To form a 2-cell embedding, the remaining edge must be in a
			different homology class.  But then, it is impossible to orient the
			half-edges alternately incoming and outgoing around the vertex,
			contradicting the orientability.

			So we have three pairwise non-homologous simple loops.  They only
			intersect at the vertex, hence any two of them form a homology basis, and
			the remaining loop must be homologous to their concatenation.  So the
			graph must be homeomorphic to that of the aH.

			Any configuration with this graph is trivially balanced because all edges
			are loops.  It is also trivially rigid as there is only one vertex, so
			the cut space is trivial.

		\item[$|\gV|=2$, $|\gE|=4$]

			We first prove that such a graph has no loop.  The graph is connected,
			hence the edges can not be all loops.  If exactly one or three edges are
			loops, the degree of a vertex will be smaller than $4$, contradicting
			Assumption~\ref{ass:degree}.  If exactly two edges are loops, they must
			be adjacent to different vertices.  So they divide the torus into two
			annuli.  Since the graph is represented as a limit of 2-cell embeddings,
			the remaining edges must lie in different annuli.  But then, it is
			impossible to orient the half-edges alternately incoming and outgoing
			around each vertex, contradicting the orientability.  This proves that
			none of the edges is a loop.

			So we have four edges between two vertices.  Any two of the edges form a
			cycle.  If some of these cycles are null-homologous, we will have
			parallel edges.  If at most one edge is simple, it is not possible to
			form a 2-cell embedding.  If exactly two edges are simple, the only
			2-cell embedding does not have a proper orientation.

			So we have four simple edges between two vertices.  Then the graph must
			be as shown in Figure~\ref{fig:general}.  This graph is balanced if and
			only if the half-edges form two collinear pairs around each vertex.  So
			if one vertex is at $0$, the other vertex must lie at a 2-division point.

			By the period condition~\eqref{eq:Vperiod}, the phase function must be of
			the form
			\[
				c-\Psi_2/2, \quad c+\Psi_2/2, \quad c+\Psi_1+\Psi_2/2, \quad, c+\Psi_1-\Psi_2/2
			\]
			on the half-edges around a vertex, where $\Psi_1$ and $\Psi_2$ are the
			fundamental shifts.  See Figure~\ref{fig:general}.

			\begin{figure}[!htb]
				\begin{tikzpicture}[scale=1,baseline]
	\draw[dotted] (-2,-2) -- (1,-2) -- (2,2) -- (-1,2) -- cycle;
	\fill (0,0) node[circle,inner sep=0.05cm,fill=black]{};
	\fill (-2,-2) node[circle,inner sep=0.05cm,fill=black,label=left:$0$]{};
	\fill (1,-2) node[circle,inner sep=0.05cm,fill=black,label=right:$T_1$]{};
	\fill (2,2) node[circle,inner sep=0.05cm,fill=black]{};
	\fill (-1,2) node[circle,inner sep=0.05cm,fill=black,label=left:$T_2$]{};
	\draw[wing] (0,0) -- (-2,-2) node[midway,above,sloped] {$\scriptstyle c-\Psi_2/2$};
	\draw[wing] (0,0) -- (2,2) node[midway,below,sloped] {$\scriptstyle c+\Psi_1+\Psi_2/2$};
	\draw[wing] (1,-2) -- (0,0) node[midway,below,sloped] {$\scriptstyle -c-\Psi_1+\Psi_2/2$};
	\draw[wing] (-1,2) -- (0,0) node[midway,above,sloped] {$\scriptstyle -c-\Psi_2/2$};
\end{tikzpicture}
				\caption{\label{fig:general}}
			\end{figure}

			If $\arg(T_2/T_1) \ne \pi/2$, the shortest edges of the graph form a
			collinear pair.  Assuming $\arg(T_2/T_1) < \pi/2$ as in
			Figure~\ref{fig:general}, then the phase function is balanced if and only
			if
			\[
				\sin(c+\Psi_1/2) = 0 \quad\text{or}\quad \cos\frac{\Psi_2-\Psi_1}{2} = 0.
			\]
			The solution $\Psi_1 = -2c$ gives Meeks' configuration which is
			generically rigid.  The solution $\Psi_2 - \Psi_1 = \pi$ gives the aI
			configurations, which is not rigid because $c$ remains a free variable.

			If $\arg(T_2/T_1) = \pi/2$, the phase function is balanced if and only if
			\[
				\sin c + \sin(c+\Psi_1) = 0 \quad\text{or}\quad \cos(\Psi_2/2) = 0.
			\]
			The solution $\Psi_1=-2c$ gives again Meeks' configurations.  The
			solutions $\Psi_1 = \pi$ and $\Psi_2 = \pi$ give the aG configurations,
			which is not rigid because $c$ remains a free variable.
	\end{description}
\end{proof}

\section{Integral of a holomorphic 1-form through a neck} \label{appendix:neck}

To prove Lemma \ref{lem:Bintegral}, we first prove a general result, Theorem
\ref{theorem:neck} below. This was in fact proved in
\cite[Lemma~1]{traizet2002} and has been used in several papers, but was not
clearly stated as an independent result and the value at $t=0$ was not
explicitly given.

Fix some numbers $0<\epsilon_1<\epsilon$. For $t\in\C$ such that
$0<|t|<\epsilon^2$, let ${\mathcal A}_t\subset\C$ be the annulus
$|t|/\epsilon<|z|<\epsilon$ and $\psi_t:{\mathcal A}_t\to{\mathcal A}_t$
be the involution defined by $\psi_t(z)=t/z$.  Let $\beta_t$ be the curve from
$\epsilon_1$ to $t/\epsilon_1$ parameterized for $s\in[0,1]$ by
\[
	\beta_t(s)=\epsilon_1^{1-2s}t^s.
\]
Note that $\beta_t$ depends on the choice of the argument of $t$.  Let $\gamma$
be the circle parameterized by $\gamma(s)=\epsilon_1e^{2\pi is}$.

\begin{theorem} \label{theorem:neck}
	Let $\omega_t$ be a family of holomorphic 1-forms on ${\mathcal A}_t$,
	depending holomorphically on $t\in D^*(0,\epsilon^2)$, and let
	$\widetilde{\omega}_t=\psi_t^*\omega_t$.  Define
	\[
		\alpha_t=\frac{1}{2\pi\ii}\int_{\gamma}\omega_t=-\frac{1}{2\pi\ii}\int_{\gamma}\widetilde{\omega}_t.
	\]
	Assume that
	\[
		\lim_{t\to 0}\omega_t=\omega_0\quad\mbox{ and }\quad\lim_{t\to 0}\widetilde{\omega}_t=\widetilde{\omega}_0
	\]
	where $\omega_0$ and $\widetilde{\omega}_0$ are holomorphic in
	$D^*(0,\epsilon)$ with at most simple poles at $z=0$, and the limit is
	uniform on compact subsets of $D^*(0,\epsilon)$.

	Then $\int_{\beta_t}\omega_t-\alpha_t\log t $ is a well-defined holomorphic
	function of $t\neq 0$ which extends holomorphically at $t=0$. Moreover, its
	value at $t=0$ is
	\[
		\lim_{z\to 0}\bigg[\Big(\int_{\epsilon_1}^z \omega_0\Big)-\alpha_0\log z \bigg]
		-\lim_{z\to 0}\bigg[\Big(\int_{\epsilon_1}^z\widetilde{\omega}_0\Big)+\alpha_0\log z \bigg].
	\]
\end{theorem}

\begin{proof}
	If $\arg(t)$ is increased by $2\pi$, the homotopy class of $\beta_t$ is
	multiplied on the left by $\gamma$, so $\int_{\beta_t}\omega_t$ is increased
	by $2\pi\ii\alpha_t$. On the other hand, $\log t $ is increased by $2\pi\ii$,
	so the difference $\int_{\beta_t}\omega_t-\alpha_t\log t $ is a well-defined
	holomorphic function of $t$ in $D^*(0,\epsilon^2)$.  Using the change of
	variable rule, we write
	\begin{equation} \label{eq:neck1}
		\Big(\int_{\beta_t}\omega_t\Big)-\alpha_t\log t =
		\bigg[\Big(\int_{\epsilon_1}^{\sqrt{t}}\omega_t\Big)-\alpha_t\log\sqrt{t}\bigg]
		-\bigg[\Big(\int_{\epsilon_1}^{\sqrt{t}}\widetilde{\omega}_t\Big)+\alpha_t\log\sqrt{t}\bigg].
	\end{equation}
	To estimate the first term, we fix $\epsilon_1<\epsilon_2<\epsilon$
	and expand $\omega_t$ in Laurent series in the annulus ${\mathcal A}_t$ as
	\[
		\omega_t=\sum_{n\in\Z}a_n(t)z^{n-1}dz
	\]
	with $a_0(t)=\alpha_t$ and
	\[
		a_n(t)=\frac{1}{2\pi\ii}\int_{|z|=\epsilon_2}z^{-n}\omega_t
		=\frac{-1}{2\pi\ii}\int_{|z|=\epsilon_2}\psi_t^*(z^{-n}\omega_t)
		=\frac{-1}{2\pi\ii}\int_{|z|=\epsilon_2} t^{-n}z^n\widetilde{\omega}_t.
	\]
	Since $\omega_t$ and $\widetilde{\omega}_t$ are uniformly bounded on the
	circle $|z|=\epsilon_2$, this gives the estimates, for $n>0$ and a uniform
	constant $C$,
	\begin{equation} \label{eq:neck2}
		|a_n(t)|\leq \frac{C}{(\epsilon_2)^n}\quad\text{and}\quad
		|a_{-n}(t)|\leq \frac{C|t|^n}{(\epsilon_2)^n}.
	\end{equation}
	Then we have
	\[
		\Big(\int_{\epsilon_1}^{\sqrt{t}}\omega_t\Big)-\alpha_t\log\sqrt{t}
		=-\alpha_t\log\epsilon_1+
		\sum_{n=1}^{\infty}\bigg[
			\frac{a_n}{n}\big(t^{n/2}-(\epsilon_1)^n\big)-
			\frac{a_{-n}}{n}\big(t^{-n/2}-(\epsilon_1)^{-n}\big)
		\bigg].
	\]
	Using the estimates \eqref{eq:neck2} and $\epsilon_1<\epsilon_2$, it is
	straightforward to check that the sum is uniformly bounded with respect to
	$t$ and that we have
	\[
		\lim_{t\to 0}\bigg[\Big(\int_{\epsilon_1}^{\sqrt{t}}\omega_t\Big)-\alpha_t\log\sqrt{t}\bigg]
		=-\alpha_0\log\epsilon_1 -\sum_{n=1}^{\infty}\frac{a_n(0)}{n}\epsilon_1^n.
	\]
	On the other hand, we have
	\[
		\int_{\epsilon_1}^z\omega_0=\alpha_0\log\left(\frac{z}{\epsilon_1}\right)
		+\sum_{n=1}^{\infty}\frac{a_n(0)}{n}\left(z^n-(\epsilon_1)^n\right)
	\]
	so
	\begin{equation} \label{eq:neck3}
		\lim_{t\to 0}\bigg[\Big(\int_{\epsilon_1}^{\sqrt{t}}\omega_t\Big)-\alpha_t\log\sqrt{t}\bigg]
		=\lim_{z\to 0}\bigg[\Big(\int_{\epsilon_1}^z\omega_0\Big)-\alpha_0\log z\bigg].
	\end{equation}
	The second term in \eqref{eq:neck1} is estimated in the exact same way,
	leading to
	\begin{equation} \label{eq:neck4}
		\lim_{t\to 0}\bigg[\Big(\int_{\epsilon_1}^{\sqrt{t}}\widetilde{\omega}_t\Big)+\alpha_t\log\sqrt{t}\bigg]
		=\lim_{z\to 0}\bigg[\Big(\int_{\epsilon_1}^z\widetilde{\omega}_0\Big)+\alpha_0\log z\bigg].
	\end{equation}
	The function $\int_{\beta_t}\omega_t-\alpha_t\log t$ is bounded so extends
	holomorphically at $t=0$ by Riemann Extension Theorem.  The last point of
	Theorem \ref{theorem:neck} follows from Equations \eqref{eq:neck1},
	\eqref{eq:neck3} and \eqref{eq:neck4}.
\end{proof}

Proof of Lemma \ref{lem:Bintegral}. Recall the definition of the path $B_h$
just before Lemma \ref{lem:Bintegral}. On path number (1), $\omega_t$ depends
holomorphically on $t_h$ in a neighborhood of $0$ so
\[
	\int_{O_{v(h)}}^{w_h=\delta}\omega_t
\]
is a holomorphic function of $t_h$ in a neighborhood of $0$.  Same for path
number (3).  Regarding path number (2) we write
\begin{equation}
	\label{eq:neck5}
	\Big(\int_{w_h=\delta}^{w_h=t_h/\delta}\omega_t\Big)-\alpha_h\log t_h
	=\Big(\int_{\delta}^{t_h/\delta}(w_h^{-1})^*\omega_t\Big)-\alpha_h\log t_h
\end{equation}
and we apply Theorem \ref{theorem:neck} to the 1-form $(w_h^{-1})^*\omega_t$
with $\epsilon_1=\delta$, $t=t_h$, observing that
\[
	\psi_{t_h}^* (w_h^{-1})^*\omega_t=((\psi_{t_h}\circ w_h)^{-1})^*\omega_t
	=(w_{-h}^{-1})^*\omega_t
\]
so the hypotheses of Theorem \ref{theorem:neck} are satisfied.  Hence
\eqref{eq:neck5} extends holomorphically at $t_h=0$ and its value there is
\begin{align*}
	&\lim_{z\to 0}\bigg[\Big(\int_{\delta}^z(w_h^{-1})^*\omega_0\Big)-\alpha_h\log z \bigg]
	-\lim_{z\to 0}\bigg[\Big(\int_{\delta}^z(w_{-h}^{-1})^*\omega_0\Big)-\alpha_{-h}\log z \bigg]\\
	=&
	\lim_{z\to p_h}\bigg[\Big(\int_{w_h=\delta}^z\omega_0\Big)-\alpha_h\log w_h(z) \bigg]
	-\lim_{z\to p_{-h}}\bigg[\Big(\int_{w_{-h}=\delta}^z\omega_0\Big)-\alpha_{-h}\log w_{-h}(z) \bigg].
\end{align*}
Adding the three terms gives Lemma \ref{lem:Bintegral}.

\bibliography{References}
\bibliographystyle{alpha}

\end{document}